\documentclass[letter,12pt]{amsart}
\usepackage{amsmath,bm}
\usepackage{amsfonts}
\usepackage{amssymb}
\usepackage{graphicx}
\usepackage{rotating}
\usepackage[height=8.7in,left=1in,right=1in,bottom=1in]{geometry}
\usepackage[breaklinks=true,bookmarksopen=true,colorlinks=true,citecolor=blue]{hyperref}
\usepackage{color}
\usepackage{colortbl}
\usepackage{paralist}
\usepackage{amsmath}
\usepackage{float}

\RequirePackage{hypernat}

\def \@seccntformat#1{\csname the#1\endcsname.\quad}
\makeatother
\numberwithin{equation}{section}

\begin{document}
\title[Auto DML for Causal Effects]{Automatic Debiased Machine Learning of
Causal and Structural Effects}
\thanks{The present paper formed the basis of the Fisher-Schultz Lecture given by
Victor Chernozhukov at the 2019 European Meeting of the Econometric Society
in Manchester. This research was supported by NSF grants 1559172 and 1757140. Helpful
comments were provided by the editor G. Imbens, three referees, J. Robins,
Y. Zhu, and participants at a 2016 demand workshop at Boston College and a
2018 machine learning and statistical inference workshop at the Banff
International Research Station.}
\author{Victor Chernozhukov}
\author{Whitney K. Newey}
\author{Rahul Singh}
\date{October 7, 2021}

\begin{abstract}
Many causal and structural effects depend on regressions. Examples include
policy effects, average derivatives, regression decompositions, average
treatment effects, causal mediation, and parameters of economic structural
models. The regressions may be high dimensional, making machine learning
useful. Plugging machine learners into identifying equations can lead to
poor inference due to bias from regularization and/or model selection. This
paper gives automatic debiasing for linear and nonlinear functions of
regressions. The debiasing is automatic in using Lasso and the function of
interest without the full form of the bias correction. The debiasing can be
applied to any regression learner, including neural nets, random forests,
Lasso, boosting, and other high dimensional methods. In addition to
providing the bias correction we give standard errors that are robust to
misspecification, convergence rates for the bias correction, and primitive
conditions for asymptotic inference for estimators of a variety of
estimators of structural and causal effects. The automatic debiased machine
learning is used to estimate the average treatment effect on the treated for
the NSW job training data and to estimate demand elasticities from Nielsen
scanner data while allowing preferences to be correlated with prices and
income.

\end{abstract}

\maketitle

%\maketitle

%\tableofcontents

\section{Introduction}

Many causal and structural parameters of economic interest depend on
regressions, i.e. on conditional expectations or least squares projections.
Examples include policy effects, average derivatives, regression
decompositions, average treatment effects, causal mediation, and parameters
of economic structural models. Often, regressions may be high dimensional,
depending on many variables. There may be many covariates for policy
effects, average derivatives, and treatment effects, or many prices and
covariates in the economic demand for some commodity. This paper is about
estimating economic and causal parameters that depend on high dimensional
regressions.

Machine learning is a collection of modern, adaptive statistical learning
methods for estimating regression functions and other statistical objects.
These methods exploit structured parsimony restrictions (such as approximate
sparsity) on regressions, together with various forms of regularization and
model selection, to enable high quality prediction in high dimensional
settings. Key methods include neural nets (deep learning), random forests,
and Lasso. The goal of this paper is to deploy these methods to infer causal
and structural parameters that depend on regression functions, including
policy, derivative, decomposition, and treatment effects as well as economic
structural parameters.

Machine learning is different than other methods in ways that are useful in
high dimensional settings. For example, Lasso has good properties with very
many potential regressors (possibly many more than sample size) when
relatively few important regressors give a good approximation but the
identity of those few is not known (i.e. the regression is approximately
sparse). In contrast, series regression is based on relatively few
regressors, often many fewer than the sample size. Lasso and series
regression are similar in that they both depend on a few regressors giving a
good approximation. They differ in that series regression requires that the
identity of the important regressors is known, while with Lasso their
identity need not be known. For Lasso, the important regressors just need to
be included somewhere among the many potential regressors. This difference
is useful in high dimensional settings, where there are potentially very
many regressors needed to approximate a function of many variables.
Typically, economics and statistics provide little guidance about which
regressors are important. With Lasso, such information is not needed, since
very many terms can be included among the potential regressors. Other
machine learning methods, such as random forests and neural nets, are also
well suited to high dimensional regression.

Machine learners provide remarkably good predictions in a variety of
settings but are inherently biased. The bias arises from using
regularization and/or model selection to control the variance of the
prediction. To obtain small mean squared prediction errors, machine learners
regularize and/or select among models so that variance and squared bias are
approximately equal. Although such equality is good for prediction, it is
not good for inference. Confidence intervals based on estimators with
approximately equal variance and squared bias will tend to have poor
coverage. This inference problem can be even worse when machine learners are
plugged into a formula for a causal or structural effect. These formulae
often involve averaging over regressor values which reduces variance without
affecting as much the bias. Variance could also potentially also be a
problem but machine learners control that for prediction purposes.

For causal and structural estimators that plug-in regularized machine
learners, the squared bias can shrink slower than the variance, leading to
extremely poor confidence interval coverage and estimators that are not
root-n consistent. Chernozhukov et al. (2017, 2018) give Lasso and random
forest examples respectively and Chernozhukov et al. (2020) shows that Lasso
plug-in estimators are not root-n consistent. Model selection inherent in
machine learners also creates inference problems. Model selection creates
bias from incorrect model choice under local alternatives, making the usual
asymptotic confidence intervals invalid over local alternatives, as shown by
Leeb and Potscher (2008a,b). Estimators of parameters of interest obtained
by plugging in machine learners can inherit this problem, as pointed out by
Belloni, Chernozhukov, and Kato (2015) and Chernozhukov, Hansen, and
Spindler (2015) and shown in Chernozhukov et al. (2020).

To reduce regularization and model selection bias we use a Neyman orthogonal
moment function where there is no first-order effect of the regression on
the expected moment function. The orthogonal moment function is constructed
by adding to an identifying moment the nonparametric influence function of
the regression on the identifying moment function. This construction is
model free, nonparametric, and based on the probability limit of the
regression learner for any distribution, as in Chernozhukov et al. (2016,
2020). As a result the orthogonality property is model free, meaning that
regression learners have no first order effect on the moments for
unrestricted, possibly misspecified, nonparametric distributions.
Consequently the standard errors are robust to misspecification because they
are constructed from the orthogonal moments while ignoring the presence of
the regression learners.

The orthogonal moment function depends on another unknown function $\bar{%
\alpha}$ in addition to the regression. We develop a Lasso minimum distance
learner of $\bar{\alpha}$ that is automatic and nonparametric, in the sense
that it depends only on the identifying moment function and not on the form
of $\bar{\alpha}$. The structure of the identifying moment function is used
to approximate $\bar{\alpha}$ as a linear combination of a dictionary (i.e.
basis) of known functions. We use the Lasso learner of $\bar{\alpha}$ and a
regression learner in the orthogonal moment functions to construct an
automatic debiased machine learner (Auto-DML) of parameters of interest. We
introduce debiased machine learning estimators for a wide variety of
effects, including policy effects, average derivatives, bounds on average
equivalent variation, and any other linear function of a regression where
debiased machine learners were not previously available. We also allow for
the identifying moment functions to be nonlinear in regressions. In addition
we give novel estimators of average treatment effects, causal mediation, and
regression decomposition.

We allow any regression learner, including neural nets, random forests,
Lasso, and other high dimensional learners to be used in the orthogonal
moment function. The primary requirement of the regression learner is that
the product of mean-square convergence rates for the learner of $\bar{\alpha}
$ and the regression learner is faster than $n^{-1/2}.$ Under this condition
and a few other regularity conditions we show root-n consistency and
asymptotic normality of the estimator of the parameter of interest. We give
convergence rates for the Lasso learner of $\bar{\alpha}$ and combine them
with existing convergence rates for regressions to verify conditions for
particular estimators. A learner of $\bar{\alpha}$ and large sample theory
is given for parameters that depend nonlinearly on regressions as well as
parameters that are linear in a regression.

The large sample theory in this paper takes the probability limit of the
regression learner and $\bar{\alpha}$ to be fixed. It would be
straightforward to extend the results to allow the regression limit and $%
\bar{\alpha}$ to change with sample size. Such a change would allow us to
accommodate sparse specifications where number of nonzero coefficients in
the true regression grows with the sample size but would complicate notation
and detail. We choose to work with a fixed regression for simplicity while
accommodating high dimensional regressions via approximate sparsity.

We give an application to estimating the treatment effect on the treated of
job training from the National Supported Work Demonstration (NSW). For many
large sets of covariates, we find similar estimates based on neural net,
random forest, and Lasso regressions with the automatic bias correction for
each. We also give an application to estimating price elasticities from
scanner panel data while allowing endogeneity of prices. We estimate the
elasticities from Auto-DML of an average derivative that includes many
covariates that account for correlated random effects. We find price
elasticities that are much smaller than cross-section elasticities,
consistent with though larger t than fixed effects elasticities found in
Chernozhukov, Hausman, and Newey (2021). We also find that plug in estimates
are similar to the cross-section elasticity estimates, so that debiasing is
important in this application.

The estimators of parameters of interest use cross-fitting, as in
Chernozhukov et al. (2018), where orthogonal moment functions are averaged
over groups of observations, the regression and $\bar{\alpha}$ learners use
all observations not in the group, and each observation is included in the
average over one group. Cross-fitting removes a source of bias and
eliminates any need for Donsker conditions for the regression learner. Early
work by Bickel (1982), Schick (1986), and Klaassen (1987) used similar
sample splitting ideas.

Auto-DML for a general linear functional of a regression, convergence rates,
and asymptotic normality results for a Dantzig selector of $\bar{\alpha}$
and the regression were given in Chernozhukov, Newey, and Robins (2018).
Chernozhukov, Newey, and Singh (2018) gave Auto-DML for any regression
learner, for nonlinear functions of a regression, and convergence rates for
a Lasso learner of $\bar{\alpha}.$ The current paper is a revised version of
Chernozhukov, Newey, and Singh (2018) with a different title. Chernozhukov,
Newey, and Singh (2019) is a revised version of Chernozhukov, Newey, and
Robins (2018) and is distinguished from the current paper and previous work
in giving and analyzing Auto-DML for local (nonparametric) effects as well
as focusing on the Dantzig selector for $\bar{\alpha}$ and the regression
for global effects. All of these papers make use of model free orthogonal
moment functions for regression learners given in Chernozuhkov et al. (2016)
and the automatic debiasing in Chernozhukov et al. (2020) builds on this
paper. The combined use of cross-fitting and orthogonal moment functions for
debiased machine learning is like Chernozhukov et al. (2018). The Auto-DML
in Chernozhukov, Newey, and Robins (2018), Chernozhukov, Newey, and Singh
(2018), and here innovates by not requiring an explicit formula for the bias
correction that is required in Chernozhukov et al. (2018) and earlier papers.

This work builds upon ideas in classical semi- and nonparametric learning
theory with low-dimensional regressions using traditional smoothing methods
(Van Der Vaart, 1991; Bickel et al., 1993; Newey 1994; Robins and Rotnitzky,
1995; Van der Vaart, 1998), that do not apply to the current
high-dimensional setting. The orthogonal moment functions developed in
Chernozhukov et al. (2016) and used here build on previous work on model
free orthogonal moment functions. Hasminskii and Ibragimov (1979) and Bickel
and Ritov (1988) suggest such estimators for functionals of a density. Newey
(1994) develops such scores for densities and regressions from computation
of the semiparametric efficiency bound for regular functionals. Doubly
robust estimating equations for treatment effects as in Robins, Rotnitzky,
and Zhao (1995) and Robins and Rotnitzky (1995) constitute model based
orthogonal moment functions and have motivated much subsequent work. Newey,
Hsieh, and Robins (1998, 2004) extend model free orthogonal moment functions
to any functional of a density or distribution in a low dimensional setting.
Model free, orthogonal moments for any learner are given and their general
properties derived in Chernozhukov et al. (2016, 2020). We use those model
free, orthogonal moment functions for regressions.

This paper also builds upon and contributes to the literature on modern
orthogonal/debiased estimation and inference, including Zhang and Zhang
(2014), Belloni et al. (2012, 2014a,b), Robins et al. (2013), van der Laan
and Rose (2011), Javanmard and Montanari (2014a,b, 2015), Van de Geer et al.
(2014), Farrell (2015), Ning and Liu (2017), Chernozhukov et al. (2015),
Neykov et al. (2018), Ren et al. (2015), Jankova and Van De Geer (2015,
2016a, 2016b), Bradic and Kolar (2017), Zhu and Bradic (2017a,b). This prior
work is about regression coefficients, treatment effects, and semiparametric
likelihood models. The objects of interest we consider are different than
those analyzed in Cai and Guo (2017). The continuity properties of
functionals we consider provide additional structure that we exploit, namely
the $\bar{\alpha}\,$, an object that is not considered in Cai and Guo
(2017). 

Targeted maximum likelihood was developed by Scharfstein, Rotnitzky, Robins
(1999) and Van Der Laan and Rubin (2006). The use of machine learning for
these estimators was proposed by Van der Laan and Rose (2011) and large
sample theory given by Luedtke and Van Der Laan (2016), Toth and van der
Laan (2016), and Zheng et al. (2016). In this paper we give a targeted
version of Auto-DML with automatic debiasing that we refer to as Auto-TML.
This estimator differs from previous ones in the objects we consider and the
use of automatic debiasing in Auto-TML.

Various papers have considered direct estimation of $\bar{\alpha}$ for
treatment effects, where $\bar{\alpha}$ is a Riesz representer that depends
on inverse propensity scores. Our work is the first to present a framework
for direct estimation of the Riesz representer of a broad class of linear
and nonlinear functionals, in a high-dimensional setting, without requiring
strong Donsker class assumptions. The earliest reference of which we know is
Robins et al. (2007), which gives a linear estimator for $\bar{\alpha}$ for
only the average treatment effect. Vermeulen and Vansteelandt (2015) base
parametric propensity score and regression estimators on double robustness
conditions for the average treatment effect. We differ in using a linear
approximation to $\bar{\alpha}$, which is restrictive in a parametric
setting but is general in high dimensional and/or nonparametric settings.
Newey and Robins (2018) present and analyze estimators based on regression
splines, while we present and analyze sparse methods for the
high-dimensional setting. The Lasso minimum distance learner of $\bar{\alpha}
$ given in Chernozhukov, Newey, and Singh (2018) and here is a direct
estimator of the Riesz representer for a broad class of linear and nonlinear
functionals that can be interpreted as being based on orthogonality of the
moment functions. Chernozhukov et al. (2020) extends this learner of $\bar{%
\alpha}$ to functions of high dimensional regression quantiles and other
objects.

In independent work on treatment effects Avagyan and Vansteelandt (2017)
give a model assisted estimator based on regularized first order conditions
and Tan (2020) developed a model assisted, multistep method of doubly robust
estimation with Lasso type regression learners having standard errors that
are robust to misspecification of the regression or propensity score.
Smucler, Rotnitzky, and Robins (2019) extended that approach to the linear
functionals of a regression considered in Chernozhukov, Newey, and Singh
(2018). For treatment effects the estimator we give is single step, allows
for any regression learner (e.g. neural nets), is model free, and has
correct standard errors if either or both the regression and the propensity
score are misspecified. Farrell, Liang, and Misra. (2021) gave a neural nets
and model based estimator of the average treatment effect and Wooldridge and
Zhu (2020) give a Lasso based debiased machine learner for panel data with
correlated random effects that depend on high dimensional regressions. Our
results also allow for a neural net regression learner but are model free
with specification robust standard error.

Chernozhukov, Newey, and Robins (2018) gave Auto-DML for linear functionals
using the Dantzig selector. More recently Hirshberg and Wager (2018) gave estimators for linear functionals based on minimax estimation of sample weights that are consistent for realizations of $\bar{\alpha}$ in sample mean square error, rather than a linear approximation to the $\bar{\alpha}$ function, in the low dimensional case, using the same orthogonal moment functions considered here.
The objects considered by Chernozhukov, Newey, and Robins (2018) include
average derivatives. More recently Hirshberg and Wager (2020) gave an
average derivative estimator based on debiasing a Lasso regression learner
of a single index high dimensional regression and Rothenhausler and Yu
(2019) gave an average derivative estimator using debiased Lasso regression.
Singh and Sun (2019) extend the present work to the instrumental variable
setting and present estimators of the local average treatment effect,
average complier characteristics, and complier counter factual
distributions. Previous to the current version of this paper Farbmacher et
al. (2020) gave DML (debiased machine learning) for causal mediation. We
propose an Auto-DML for causal mediation analysis as an example in Section 5.

In summary, contributions of the paper include the construction of DML for a
wide range of interesting policy effects and structural parameters where DML
was not previously available. This construction is based on a Lasso minimum
distance learner of $\bar{\alpha}$ we propose. The debiasing and inference
is model free and robust to misspecification and carried out in a single
step, unlike previous estimators of average treatment effects. For average
treatment and other effects we construct DML for a variety of regression
learners, such as neural nets, random forests, or high dimensional methods.

In Section 2 we describe the objects of interest we consider and associated
orthogonal moment functions. In Section 3 we give the Lasso learner of $\bar{%
\alpha},$ the Auto-DML and Auto-TML estimators, and a consistent estimator
of their asymptotic variance. Section 4 derives mean square convergence
rates for the Lasso learner of $\bar{\alpha}$ and conditions for root-n
consistency and asymptotic normality of Auto-DML and Auto-TML including
primitive conditions in examples. Section 5 gives Auto-DML for nonlinear
functionals of multiple regressions and as an example develops Auto-DML for
causal mediation analysis. Section 6 gives Auto-DML for regression
decomposition and estimates the average treatment on the treated for the NSW
experiment. Section 7 gives Auto-DML estimates of price elasticities that
allow for correlated random effects in scanner panel data. Section 8 offers
some conclusions and possible extensions.

\section{Average Linear Effects and Orthogonal Moment Functions}

For expositional purposes, in this Section we first consider parameters that
depend linearly on a single conditional expectation. To describe such an
object, let $W$ denote a data observation, and consider a subvector $%
(Y,X^{\prime})^{\prime}$ where $Y$ is a scalar outcome with finite second
moment and $X$ is a covariate vector. Denote the conditional expectation of $%
Y$ given $X\in\mathcal{X}$ as%
\begin{equation*}
\gamma_{0}(x)=\mathrm{E}[Y|X=x].
\end{equation*}
Let $m(w,\gamma)$ denote a function of the function $\gamma$ (i.e. a
functional of $\gamma),$ where $\gamma$ denotes a possible conditional
expectation function $\gamma:\mathcal{X}\longrightarrow\mathbb{R}$, that
depends on a data observation $w$ and is linear in $\gamma.$ We will
consider effects of the form%
\begin{equation*}
\theta_{0}=\mathrm{E}[m(W,\gamma_{0})].
\end{equation*}
The parameter of interest $\theta_{0}$ is an expectation of some known
formula $m(W,\gamma)$ of a data observation $W$ and a regression $\gamma.$

We also give results in later Sections for important parameters having more
general forms. In Section 5 we allow $m(W,\gamma)$ to be nonlinear in
multiple regressions and propose an estimator of causal effects with
mediation. In Section 6 we give estimators of regression decompositions and
their properties. These important examples extend the framework of this
Section to parameters that are nonlinear in multiple regressions

Several important examples of linear effects are:

\bigskip

\textsc{Example 1:} (Average Policy Effect). An average effect of a counter
factual shift in the distribution of regressors from a known $F_{0}$ to
another known $F_{1}$, when $\gamma _{0}$ does not vary with the
distribution of $X$, is%
\begin{equation*}
\theta _{0}=\int \gamma _{0}(x)d\mu (x);\text{ }\mu (x)=F_{1}(x)-F_{0}(x).
\end{equation*}%
Here $m(w,\gamma )=\int \gamma (x)d\mu (x)$ which does not depend on $w.$
This policy effect builds on but is different than Stock (1989) in comparing
averages over two known distributions rather than the empirical distribution.

\bigskip

\textsc{Example 2:} (Weighted Average Derivative). Here $X=(D,Z)$ for a
continuously distributed random variable $D,$ $\gamma_{0}(x)=%
\gamma_{0}(d,z), $ $\omega(d)$ is a pdf, and%
\begin{equation*}
\theta_{0}=\mathrm{E} \left[ \int\omega(u)\frac{\partial\gamma_{0}(u,Z)}{%
\partial d}du\right] =\mathrm{E} \left[ \int S(u)\gamma_{0}(u,Z)\omega(u)du%
\right] =\mathrm{E}[S(U)\gamma_{0}(U,Z)],
\end{equation*}
where $S(u)=-\omega(u)^{-1}\partial\omega(u)/\partial u$ is the negative
score for the pdf $\omega(u),$ the second equality follows by integration by
parts, and $U$ is a random variable that is independent of $Z$ with pdf $%
\omega(u).$ This $U$ could be thought of as one simulation draw from the pdf 
$\omega(u).$ Here $m(w,\gamma)=S(u)\gamma(u,x)$ where $W$ includes $U.$

This $\theta _{0}$ can be interpreted as an average treatment effect on $Y$
of a continuous treatment $D$ in a model where $Y=Y(D)$ for a potential
outcome stochastic process $Y(d)$ that is independent of $D$ conditional on
covariates $Z.$ By conditional independence%
\begin{equation*}
\mathrm{E}[\gamma _{0}(u,Z)]=\int \mathrm{E}[Y(D)|D=u,Z=z]F_{Z}(dz)=\int 
\mathrm{E}[Y(u)|Z=z]F_{Z}(dz)=\mathrm{E}[Y(u)],
\end{equation*}%
for $\omega (u)>0$ assuming that the joint pdf of $(D,Z)$ is positive where $%
\omega (D)>0$, as in Chamberlain (1984), Wooldridge (2002), and Blundell and
Powell (2004). The $\mathrm{E}[Y(u)]$ is the average outcome at $D=u$ and is
sometimes referred to as the average structural function. Assuming that we
can interchange the order of differentiation and integration,%
\begin{equation*}
\theta _{0}=\int \omega (u)\frac{\partial \mathrm{E}[\gamma _{0}(u,Z)]}{%
\partial u}du=\int \frac{\partial \mathrm{E}[Y(u)]}{\partial u}\omega
(u)du=\int \mathrm{E}\left[ \frac{\partial Y(u)}{\partial u}\right] \omega
(u)du,
\end{equation*}%
similarly to Imbens and Newey (2009) and Rothenh{\"{a}}usler and Yu (2019),
which build on but are different than Powell, Stock, and Stoker (1989).
Regarding $\mathrm{E}[\partial Y(u)/\partial u]$ as the average treatment
effect at $u$ we see that $\theta _{0}$ is a weighted average treatment
effect. Alternatively, $\theta _{0}$ can be regarded as an average
derivative of the average structural function. The averaging over a known
pdf $\omega (u)$ helps fulfill regularity conditions for the Auto-DML
developed here that can be used to estimate $\theta _{0}$ for high
dimensional covariates $Z.$

\bigskip

\textsc{Example 3: }(Average Treatment Effect). In this example $X=(D,Z)$
and $\gamma_{0}(x)=\gamma_{0}(d,z)$, where $D\in\{0,1\}$ is the treatment
indicator and $Z$ are covariates. The object of interest is%
\begin{equation*}
\theta_{0}=\mathrm{E}[\gamma_{0}(1,Z)-\gamma_{0}(0,Z)].
\end{equation*}
If potential outcomes are mean independent of treatment $D$ conditional on
covariates $Z$, then $\theta_{0}$ is the average treatment effect (Rosenbaum
and Rubin, 1983). Here $m(w,\gamma)=\gamma(1,z)-\gamma(0,z).$

\bigskip

\textsc{Example 4:} (Average Equivalent Variation Bound). An economic
example is a bound on average equivalent variation for heterogenous demand.
Here $Y$ is the share of income spent on a commodity and $X=(P_{1},Z),$
where $P_{1}$ is the price of the commodity and $Z$ includes income $Z_{1}$,
prices of other goods, and other observable variables affecting utility. Let 
$\check{p}_{1}<\bar{p}_{1}$ be lower and upper prices over which the price
of the commodity can change, $\kappa$ a bound on the income effect, $%
\omega(z)$ some weight function, and $U$ a random variable that is uniformly
distributed over $(\check{p}_{1},\bar{p}_{1})$ and independent of $(Y,X).$ $%
U $ can be thought of as one simulation draw from a uniform distribution on $%
(\check{p}_{1},\bar{p}_{1}).$ The object of interest is%
\begin{equation*}
\theta_{0}=\mathrm{E} \left[ \Lambda(U,Z)\gamma_{0}(U,Z)\right] ,\text{ }%
\Lambda(u,z)=\omega(z)1(\check{p}_{1}<u<\bar{p}_{1})(\bar{p}_{1}-\check{p}%
_{1})\frac{z_{1}}{u}\exp(-\kappa\lbrack u-\check{p}_{1}]).
\end{equation*}
If individual heterogeneity in consumer preferences is independent of $X$
and $\kappa$ is a lower (upper) bound on the derivative of consumption with
respect to income for all individuals, then $\theta_{0}$ is an upper (lower)
bound on the weighted average over consumers of equivalent variation for a
change in the price of the first good from $\check{p}_{1}$ to $\bar{p}_{1}$;
see Hausman and Newey (2016). Here $m(w,\gamma)=\Lambda(u,z)\gamma(u,z),$
where $W$ includes $U.$

\bigskip

We focus on $m(w,\gamma)$ where there exists a function $\alpha_{0}(X)$ with 
$\mathrm{E}[\alpha_{0}(X)^{2}]<\infty$ and 
\begin{equation}
{\mathrm{E}}[m(W,\gamma)]={\mathrm{E}}[\alpha_{0}(X)\gamma(X)]\text{ \ for
all }\gamma\text{ such that }{\mathrm{E}}[\gamma(X)^{2}]<\infty.
\label{Riesz rep}
\end{equation}
By the Riesz representation theorem, existence of such a $\alpha_{0}(X)$ is
equivalent to $\mathrm{E}[m(W,\gamma)]$ being a mean-square continuous
functional of $\gamma,$ i.e. $\mathrm{E}[m(W,\gamma)]\leq C\left\Vert
\gamma\right\Vert $ for all $\gamma$, where $\left\Vert \gamma\right\Vert =%
\sqrt{\mathrm{E}[\gamma(X)^{2}]}$ and $C>0.$ We will refer to this $%
\alpha_{0}(X)$ as the Riesz representer (Rr). Existence of the Rr is
equivalent to the semiparametric variance bound for $\theta_{0}$ being
finite, as stated in Newey (1994) and shown in Hirshberg and Wager (2018)
for conditional expectations and in Chernozhukov, Newey, and Singh (2019)
more generally for least squares projections. Thus, in assuming existence of 
$\alpha_{0}(X)$ we are just assuming that $\theta_{0}$ has a finite
semiparametric variance bound.

Each of Examples 1-4 has such a Rr. Let $f(x)$ denote the pdf of $X$ in
Example 1, $f(d|z)$ the pdf of $D$ conditional on $Z$ in Example 2, $\pi
_{0}(z)=\Pr(D=1|Z=z)$ the propensity score in Example 3, and $f(p_{1}|z)$
the pdf of $P_{1}$ conditional on $Z$ in Example 4. Table~\ref{tab:RR} summarizes the
functional $m(w,\gamma)$ and the Rr in each of the examples:

\begin{table}[ptb]
\centering
\begin{tabular}{c|c|c}
\hline\hline
Effect & $m(W,\gamma)$ & Riesz Representer \\ \hline
Policy Effect & $\int\gamma(x)[f_{1}(x)-f_{0}(x)]dx$ & $%
f(X)^{-1}[f_{1}(X)-f_{0}(X)]$ \\ 
Weighted Average Derivative & $S(U)\gamma(U,Z)$ & $f(D|Z)^{-1}\omega(D)S(D)$
\\ 
Average Treatment Effect & $\gamma(1,Z)-\gamma(0,Z)$ & $\pi_{0}(Z)^{-1}D-(1-%
\pi_{0}(Z))^{-1}(1-D)$ \\ 
Equivalent Variation Bound & $\Lambda(U,Z)\gamma(U,Z)$ & $(\bar{p}_{1}-%
\check{p}_{1})^{-1}f(P_{1}|Z)^{-1}\Lambda(P_{1},Z)$ \\
\hline\hline
\end{tabular}
\caption{$m$ and Rr for Examples 1-4}
\label{tab:RR}
\end{table}

Equation (\ref{Riesz rep}) follows in Example 1 by multiplying and dividing
by $f(x)$ inside the integral, in Example 2 by integration and multiplying
and dividing by $f(d|z)$, in Example 3 in a standard way for average
treatment effects, and in Example 4 by multiplying and dividing by $%
f(p_{1}|z)$. For $\mathrm{E}[\alpha_{0}(X)^{2}]<\infty$ to hold the
denominator must not be too small relative to the numerator in each $%
\alpha_{0}(X)$, on average. For instance Example 3 must have $\mathrm{E}%
[\{\pi_{0}(Z)(1-\pi_{0}^{{}}(Z))\}^{-1}]<\infty.$

Equation (\ref{Riesz rep}) implies that the effect of interest can be
represented in three different ways, as%
\begin{equation*}
\theta_{0}=\mathrm{E}[m(W,\gamma_{0})]=\mathrm{E}[\alpha_{0}(X)%
\gamma_{0}(X)]=\mathrm{E}[\alpha _{0}(X)Y],
\end{equation*}
where the last equality follows by iterated expectations. Any of these three
expressions could be used to estimate $\theta_{0}$. We could estimate $%
\theta_{0}$ from the first expression using a learner (estimator) of $%
\gamma_{0}$. We could also estimate $\theta_{0}$ from the last expression
using a learner of $\alpha_{0}(X).$ In addition we could use learners of
both $\gamma_{0}$ and $\alpha_{0}$ to estimate $\theta_{0}$ from the middle
expression. We focus here on using a learner of $\gamma_{0}$, though $%
\alpha_{0}$ will be important for the bias correction to follow.

We rely on a regression learner (estimator) $\hat{\gamma}$ of $\gamma_{0}$
to estimate $\theta_{0}.$ The $\hat{\gamma}$ can be any of a variety of
machine learners including neural nets, random forests, Lasso, and other
high dimensional methods. All we require is that $\hat{\gamma}$ converge in
mean square at a sufficiently fast rate, as specified in Section 4.

Whatever the choice of $\hat{\gamma},$ estimating $\theta_{0}$ by plugging $%
\hat{\gamma}$ into $m(W,\gamma)$ and averaging over observations on $W$ can
lead to large biases when $\hat{\gamma}$ involves regularization and/or
model selection, as discussed in the Introduction. For that reason we use an
orthogonal moment function for $\theta_{0}$, where the regression learner $%
\hat{\gamma}$ has no first-order effect on the moments. We follow
Chernozhukov et al. (2016, 2020) in basing the orthogonal moment function on
the probability limit (plim) $\gamma(F)$ of $\hat{\gamma}$ when one
observation $W$ has CDF $F,$ where $F$ is unrestricted except for regularity
conditions. Here $\gamma(F)$ can be thought of as the plim of $\hat{\gamma}$
under general misspecification, where $\gamma(F)$ need not be the
conditional expectation $\mathrm{E}_{F}[Y|X]$.

The plim $\gamma(F)$ of $\hat{\gamma}$ depends on the learner. For example
Lasso, the Dantzig selector, boosting, and other high dimensional methods
are based on a sequence of potential regressors $X=(X_{1},X_{2},...)$. These
learners have the form%
\begin{equation*}
\hat{\gamma}(x)=\sum_{j=1}^{\infty}\hat{\beta}_{j}x_{j}\text{, }\hat{\beta }%
_{j^{\prime}}\neq0\text{ for a finite number of }j^{\prime}\text{,}
\end{equation*}
where $x=(x_{1},x_{2},...)$ denotes a possible realization of $X$. Because
each $\hat{\gamma}(X)$ is a linear combination of $X=(X_{1},X_{2},...)$ the
plim $\gamma(F)$ of $\hat{\gamma}$ will also be a linear combination of $X$,
or at least will be approximated by such a linear combination. Define $%
\Gamma $ to be the mean square closure of the set of finite linear
combinations of $X$, i.e. $\Gamma$ is the set of $\gamma(X)$ such that $%
\mathrm{E}[\gamma(X)^{2}]<\infty$ and for every $\varepsilon>0$ there exists 
$(\beta_{j}^{\varepsilon})_{j=1}^{\infty}$ such that $\beta_{j^{\prime}}^{%
\varepsilon}\neq0$ for a finite number of $j^{\prime}$ and $\mathrm{E}%
[\{\gamma(X)-\sum_{j=1}^{\infty}\beta
_{j}^{\varepsilon}X_{j}\}^{2}]<\varepsilon.$ It will be the case that $%
\gamma(F)\in\Gamma.$ Because Lasso and other high dimensional methods are
being used for least squares prediction of $Y$ it will also be the case that%
\begin{equation}
\gamma(F)=\arg\min_{\gamma\in\Gamma}\mathrm{E}_{F}[\{Y-\gamma(X)\}^{2}],
\label{BLP}
\end{equation}
This $\gamma(F)$ minimizes population least squares criteria over the (mean
square closure of) linear combinations of $X,$ i.e. it is the best linear
predictor of $Y$ by linear combinations of $X.$ Here $\gamma(F)$ is the
infinite dimensional linear regression that is nonparametrically estimated
by Lasso and other high dimensional methods.

Neural nets and random forests may have a different $\gamma(F)$. A neural
net or random forest is often a nonparametric regression estimator for a
finite (but high) dimensional $X$. In that case%
\begin{equation*}
\gamma(F)=\mathrm{E}_{F}[Y|X],
\end{equation*}
which satisfies equation (\ref{BLP}) when $\Gamma$ is the set of all
(measurable) functions of $X$ with finite second moment. The plim of Lasso
and other high dimensional methods will also be this $\gamma(F)$ if $%
X=(X_{1},X_{2},...)$ can approximate any function of a fixed set of
regressors, but otherwise will not. A third type of learner $\hat{\gamma}$
is one that imposes additivity restrictions on $\hat{\gamma}$, such as $\hat{%
\gamma}(X)=\hat{\gamma}_{1}(X_{1})+\hat{\gamma}_{2}(X_{2})$, allowing for
nonparametric learners $\hat{\gamma}_{1}(X_{1})$ and $\hat{\gamma}%
_{2}(X_{2}).$ In that case $\gamma(F)$ will be satisfy equation (\ref{BLP})
where $\Gamma$ is the mean square closure of functions that are additive in $%
X_{1}$ and $X_{2}$.

We use the orthogonal moment function from Chernozhukov et al. (2016, 2020)
for a regression learner $\hat{\gamma}$ having plim $\gamma(F)$ satisfying
equation (\ref{BLP}) for any linear, closed $\Gamma.$ The orthogonal moment
function is constructed by adding to the identifying moment function $%
m(w,\gamma)-\theta$ the nonparametric influence function of of $\mathrm{E}%
[m(W,\gamma (F))].$ As shown in Newey (1994) the nonparametric influence
function of $\mathrm{E}[m(W,\gamma(F))]$ is%
\begin{equation*}
\bar{\alpha}(X)[Y-\bar{\gamma}(X)],
\end{equation*}
where $\bar{\gamma}(X)$ is the solution to equation (\ref{BLP}) for $F=F_{0}$
and $\bar{\alpha}\in\Gamma$ satisfies $\mathrm{E}[m(W,\gamma)]=\mathrm{E}[%
\bar{\alpha}(X)\gamma(X)]$ for all $\gamma\in\Gamma.$ As in Chernozhukov,
Newey, and Singh (2019),%
\begin{equation}
\bar{\alpha}=\arg\min_{\alpha\in\Gamma}\mathrm{E}[\{\alpha_{0}(X)-\alpha(X)%
\}^{2}].  \label{BLP alpha}
\end{equation}
This $\bar{\alpha}$ can be thought of as the Riesz representer for the
linear functional $\mathrm{E}[m(W,\gamma)]$ with domain $\Gamma.$ Evaluating
the nonparametric influence function at possible values $\gamma$ and $\alpha$
of $\bar{\gamma}$ and $\bar{\alpha}$ and adding it to the the identifying
moment function gives the orthogonal moment function%
\begin{equation}
\psi(w,\theta,\gamma,\alpha)=m(w,\gamma)-\theta+\alpha(x)[y-\gamma(x)].
\label{Debiased mom}
\end{equation}

The moment function $\psi(w,\theta,\gamma,\alpha)$ depends on a possible
value $\alpha$ of the unknown function $\bar{\alpha}$ as well as a possible
value $\gamma$ of the plim $\bar{\gamma}$ of the regression learner. A
learner $\hat{\alpha}$ of $\bar{\alpha}$ is needed to use this orthogonal
moment function to estimate $\theta_{0}.$ In Section 3 we will describe how
to construct $\hat{\alpha}.$ In Chernozhukov et al. (2016, 2020) $\psi
(w,\theta,\gamma,\alpha)$ is shown to be orthogonal without being specific
about the form of $\hat{\alpha}.$ For exposition we repeat that
demonstration here. Consider any $\gamma,\alpha\in\Gamma$, representing
possible realizations of learners $\hat{\gamma}$ and $\hat{\alpha}$ that are
in $\Gamma.$ The well known necessary and sufficient conditions for equation
(\ref{BLP}) with $F=F_{0}$ are that $\mathrm{E}[\alpha(X)\{Y-\bar{\gamma}%
(X)\}]=0$ for all $\alpha\in\Gamma.$ Therefore%
\begin{align}
\mathrm{E}[\psi(W,\theta,\gamma,\alpha)-\psi(W,\theta,\bar{\gamma},\bar{%
\alpha})] & =\mathrm{E}[m(W,\gamma)]-\mathrm{E}[m(W,\bar{\gamma})]+\mathrm{E}%
[\alpha(X)\{Y-\gamma (X)\}]  \label{2nd order} \\
& =\mathrm{E}[\alpha_{0}(X)\{\gamma(X)-\bar{\gamma}(X)\}]+\mathrm{E}%
[\alpha(X)\{Y-\gamma (X)\}]  \notag \\
& =\mathrm{E}[\bar{\alpha}(X)\{\gamma(X)-\bar{\gamma}(X)\}]+\mathrm{E}%
[\alpha(X)\{\bar{\gamma }(X)-\gamma(X)\}]  \notag \\
& =-\mathrm{E}[\{\alpha(X)-\bar{\alpha}(X)\}\{\gamma(X)-\bar{\gamma}(X)\}], 
\notag
\end{align}
where the second equality follows by equation (\ref{Riesz rep}) and the
third equality by the necessary and sufficient condition for equation (\ref%
{BLP alpha}) that $\mathrm{E}[\{\alpha_{0}(X)-\bar{\alpha}(X)\}\gamma(X)]=0$
for all $\gamma\in\Gamma.$ Here we see that $\psi(w,\theta,\gamma,\bar{\alpha%
})$ "partials out" $\gamma$ in the sense that%
\begin{equation*}
\mathrm{E}[m(W,\gamma)+\bar{\alpha}(X)\{Y-\gamma(X)\}]=\mathrm{E}[m(W,\bar{%
\gamma})]
\end{equation*}
does not depend on $\gamma$. Also equation (\ref{2nd order}) gives an
explicit formula showing that the effect of $\gamma$ and $\alpha$ on $%
\mathrm{E}[\psi (W,\theta,\gamma,\alpha)]$ is second order and hence $%
\psi(W,\theta ,\gamma,\alpha)$ is orthogonal.

The orthogonality property of $\psi (W,\theta ,\gamma ,\alpha )$ only
depends on $\gamma ,$ $\alpha \in \Gamma $ and $\bar{\gamma}$ satisfying
equation (\ref{BLP}). In particular orthogonality does not depend on either $%
\bar{\gamma}$ being $\mathrm{E}[Y|X]$ or on $\bar{\alpha}=\alpha _{0}.$ In
this sense orthogonality of $\psi (W,\theta ,\gamma ,\alpha )$ is model
free, i.e. nonparametric. Consequently the estimator of $\theta $ will be
asymptotically normal and standard errors consistent even if either $\bar{%
\gamma}\neq \gamma _{0}$ or $\bar{\alpha}\neq \alpha _{0}$ or both, which is
possible when neither $\gamma _{0}(X)=\mathrm{E}[Y|X]$ nor $\alpha _{0}(X)$
satisfying equation (\ref{Riesz rep}) is an element of $\Gamma .$ This
robustness of the standard errors results from the orthogonality of the
moments only depending on the $\bar{\gamma}$ limit of the regression
estimator, so that the sample average of the estimated orthogonal moment
function will be asymptotically equivalent to the sample average at the
truth, without any model assumptions.

The orthogonal moment function could also be viewed as the efficient
influence function of $\mathrm{E}[m(W,\bar{\gamma})]$ which clarifies that
the Auto-DML is an efficient semiparametric estimator of $\mathrm{E}[m(W,%
\bar{\gamma})]$. Viewing $\psi(w,\theta,\gamma,\alpha)$ in this way is not
useful for debiasing because the results of Chernozhukov et. al. (2016,
2020) already imply model free orthogonality.

The moment function $\psi(w,\theta,\gamma,\alpha)$ is doubly robust for
estimation of the true parameter $\theta_{0}.$ Evaluating at $\theta_{0},%
\bar{\gamma},\bar{\alpha}$ and taking the expectation gives%
\begin{align}
\mathrm{E}[\psi(W,\theta_{0},\bar{\gamma},\bar{\alpha})] & =\mathrm{E}[m(W,%
\bar{\gamma })]-\theta_{0}+\mathrm{E}[\bar{\alpha}(X)\{Y-\bar{\gamma}(X)\}]
\label{Double rob} \\
& =\mathrm{E}[\alpha_{0}(X)\{\bar{\gamma}(X)-\gamma_{0}(X)\}]+\mathrm{E}[%
\bar{\alpha }(X)\{\gamma_{0}(X)-\bar{\gamma}(X)\}]  \notag \\
& =-\mathrm{E}[\{\bar{\alpha}(X)-\alpha_{0}(X)\}\{\bar{\gamma}(X)-\gamma
_{0}(X)\}],  \notag
\end{align}
which is zero for $\bar{\gamma}=\gamma_{0}$ or $\bar{\alpha}=\alpha_{0}.$
Thus $\mathrm{E}[\psi(W,\theta_{0},\bar{\gamma},\bar{\alpha})]=0$, so that
the orthogonal moment condition identifies $\theta_{0},$ when either $\bar{%
\gamma}(X)=\mathrm{E}[Y|X]$ or $\alpha_{0}(X)\in\Gamma.$ These conditions
both hold when the regression learner is nonparametric so that $\Gamma$ is
the set of all functions of $X$ with finite second moment. For high
dimensional regressions where $\Gamma$ is the closed linear span of $%
X=(X_{1},X_{2},...)$ the plim of the learner $\hat{\gamma}$ may not be $%
\mathrm{E}[Y|X]$ but the orthogonal moment function still identifies $%
\theta_{0}$ when $\alpha_{0}(X)\in\Gamma.$ That is, $\theta_{0}$ is
identified when $\alpha_{0}(X)$ can be approximated arbitrarily well in mean
square by a linear combination of $X.$ This robustness condition can be
interpreted in each of Examples 1-4:

\bigskip\ 

\textsc{Example 1:} For high dimensional $\hat{\gamma},$ where $\Gamma$ is
the mean square closure of linear combinations of $X,$ $\mathrm{E}%
[\psi(W,\theta_{0},\bar{\gamma},\bar{\alpha})]=0$ even when $\bar{\gamma}%
(X)\neq \mathrm{E}[Y|X]$ if $\alpha_{0}(X)=[f_{1}(X)-f_{0}(X)]/f(X)\in%
\Gamma. $

\bigskip

\textsc{Example 2:} For high dimensional $\hat{\gamma},$ where $\Gamma$ is
the mean square closure of linear combinations of $X,$ $\mathrm{E}%
[\psi(W,\theta_{0},\bar{\gamma},\bar{\alpha})]=0$ even when $\bar{\gamma}%
(X)\neq \mathrm{E}[Y|X]$ if $\alpha_{0}(X)=f(D|Z)^{-1}\omega(D)S(D)\in%
\Gamma. $

\bigskip

\textsc{Example 3: }For the average treatment effect where $\Gamma$ is
nonparametric, so that $\bar{\gamma}(X)=\mathrm{E}[Y|X]$ and $\bar{\alpha}%
(X)=\alpha_{0}(X),$ the orthogonal moment function in equation (\ref%
{Debiased mom}) corresponds to the seminal doubly robust moment function of
Robins, Rotnitzky, and Zhao (1995). When $\hat{\gamma}$ is high dimensional,
with say $X=(DZ,(1-D)\tilde{Z})$ for sequences $Z=(Z_{1},Z_{2},...)$ and $%
\tilde{Z}=(\tilde{Z}_{1},\tilde{Z}_{2},...)$, with each $\tilde{Z}_{j}$ a
function of $Z,$ the orthogonal moment function is%
\begin{equation*}
\psi(W,\theta,\bar{\gamma},\bar{\alpha})=\bar{\gamma}(1,Z)-\bar{\gamma }%
(1,0)-\theta+\bar{\alpha}(X)[Y-\bar{\gamma}(X)].
\end{equation*}
This orthogonal moment function is different than those previously
considered in $\bar{\alpha}(X)$ being the projection of $\alpha_{0}(X)$ on $%
\Gamma$ rather than $\alpha_{0}(X)$. Here $\mathrm{E}[\psi(W,\theta_{0},\bar{%
\gamma},\bar{\alpha})]=0$ if linear combinations of $Z$ and $\tilde{Z}$ can
approximate abitrarily well $\pi_{0}(Z)^{-1}$ and $[1-\pi_{0}(Z)]^{-1}$
respectively, even when $\bar{\gamma}(X)\neq \mathrm{E}[Y|X].$

\bigskip

For brevity we omit further discussion of Example 4 from the paper and refer
the interested reader to Chernozhukov, Hausman, and Newey (2021).

\bigskip

\section{Estimation}

To estimate (learn) $\theta_{0}$ we use cross-fitting where the orthogonal
moment function $\psi(w,\gamma,\alpha,\theta)$ is averaged over observations
different than used to estimate $\bar{\gamma}$ and $\bar{\alpha}.$ We assume
that the data $W_{i},$ $(i=1,...,n)$ are i.i.d.. Let $I_{\ell},$ $%
(\ell=1,...,L)$, be a partition of the observation index set $\{1,...,n\}$
into $L$ distinct subsets of about equal size. In practice $L=5$ (5-fold) or 
$L=10$ (10-fold) cross-fitting is often used. Let $\hat{\gamma}_{\ell}$ and $%
\hat{\alpha}_{\ell}$ be estimators constructed from the observations that
are \textit{not} in $I_{\ell}.$ We construct the estimator $\hat{\theta}$ by
setting the sample average of $\psi(W_{i},\theta,\hat{\gamma}_{\ell},\hat{%
\alpha}_{\ell})$ to zero and solving for $\theta.$ This $\hat{\theta}$ and
an associated asymptotic variance estimator $\hat{V}$ have explicit forms%
\begin{align}
\hat{\theta} & =\frac{1}{n}\sum_{\ell=1}^{L}\sum_{i\in I_{\ell}}\{m(W_{i},%
\hat{\gamma}_{\ell})+\hat{\alpha}_{\ell}(X_{i})[Y_{i}-\hat{\gamma }%
_{\ell}(X_{i})]\},  \label{Estimator} \\
\hat{V} & =\frac{1}{n}\sum_{\ell=1}^{L}\sum_{i\in I_{\ell}}\hat{\psi}%
_{i\ell}^{2},\text{ }\hat{\psi}_{i\ell}=m(W_{i},\hat{\gamma}_{\ell})-\hat{%
\theta}+\hat{\alpha}_{\ell}(X_{i})[Y_{i}-\hat{\gamma}_{\ell}(X_{i})],  \notag
\end{align}

Any regression learner $\hat{\gamma}_{\ell }$ can be used here as long as
its mean-square convergence rate is a power of $1/n,$ as assumed in Section
4. Such a convergence rate is available for neural nets (Chen and White,
1999, Schmidt-Heiber, 2020, Farrell, Liang, and Misra, 2021), random forests
(Syrgkanis and Zampetakis, 2020), Lasso (Bickel, Ritov, and Tsybakov, 2009),
boosting (Luo and Spindler, 2016), and other high dimensional methods. As a
result any of these regression learners can be used to construct an Auto-DML 
$\hat{\theta}$ from equation (\ref{Estimator}), in conjunction with a
learner $\hat{\alpha}_{\ell }$ of $\bar{\alpha}.$

The correctness of $\hat{V}$ relies on consistency of the regression learner 
$\hat{\gamma}_{\ell }.$ It would be interesting to investigate whether the
finite sample approximation could be improved by using a variance estimator
that allowed $\hat{\gamma}_{\ell }$ to not be consistent because the
dimension of the regression grows as fast as the sample size, e.g. as in
Cattaneo, Jansson, and Newey (2018).

An alternative estimator of $\theta _{0}$ can be constructed that extends
the targeted maximum likelihood approach of Scharfstein, Rotnitzky, and J.M.
Robins (1999) and van der Laan and Rubin (2006) to the objects we consider.
This Auto-TML estimator is a plug-in estimator based on a regression learner
that has been debiased in a direction specific to the object of interest.
This estimator is given by 
\begin{equation}
\tilde{\theta}=\frac{1}{n}\sum_{\ell =1}^{L}\sum_{i\in I_{\ell }}m(W_{i},%
\tilde{\gamma}_{\ell }),\text{ }\tilde{\gamma}_{\ell }(x)=\hat{\gamma}_{\ell
}(x)+\frac{\sum_{i\in I_{\ell }}\hat{\alpha}_{\ell }(X_{i})[Y_{i}-\hat{\gamma%
}_{\ell }(X_{i})]}{\sum_{i\in I_{\ell }}\hat{\alpha}_{\ell }(X_{i})^{2}}\hat{%
\alpha}_{\ell }(x).  \label{TMLE}
\end{equation}%
As with other targeted estimators the plug-in form of Auto-TML allows
imposition of constraints through $m(W,\gamma )$. In Section 4 we show that
this estimator is asymptotically equivalent to $\hat{\theta}$.

To describe $\hat{\alpha}_{\ell}$ let $b(x)=(b_{1}(x),...,b_{p}(x))$ be a $%
p\times1$ dictionary of functions of $x,$ where $p$ can be large, with each $%
b_{j}(x)$ standardized to have mean $0$ and standard deviation $1,$ to be
further discussed in this Section. For convenience we ignore dependence of $%
b(x)$ on the data in the notation. The learner $\hat{\alpha}_{\ell}$ given
here is%
\begin{align}
\hat{\alpha}_{\ell}(x) & =\frac{1}{n-n_{\ell}}\sum_{i\notin
I_{\ell}}m(W_{i},1)+b(x)^{\prime}\hat{\rho}_{\ell},\text{ }\hat{\rho}%
_{\ell}=\arg \min_{\rho}\{-2\hat{M}_{\ell}^{\prime}\rho+\rho^{\prime}\hat{G}%
_{\ell}\rho+2r\sum_{j=1}^{J}\left\vert \rho_{j}\right\vert \},
\label{CFLasso} \\
\hat{M}_{\ell} & =\frac{1}{n-n_{\ell}}\sum_{i\notin I_{\ell}}m(W_{i},b),%
\text{ }\hat{G}_{\ell}=\frac{1}{n-n_{\ell}}\sum_{i\notin
I_{\ell}}b(X_{i})b(X_{i})^{\prime},  \notag
\end{align}
where $n_{\ell}$ is the number of observations in $I_{\ell}$ and $r>0$ is a
positive scalar. This $\hat{\alpha}_{\ell}$ is used in equation (\ref%
{Estimator}) to construct $\hat{\theta}$ and $\hat{V}$.

To explain and motivate $\hat{\alpha}_{\ell}$ it is notationally convenient
to drop the $\ell$ subscript, with the understanding that $\hat{\alpha}%
_{\ell}$ is computed using only observations not in $I_{\ell}$ for each $%
\ell,$ as in equation (\ref{CFLasso}). It is also notationally convenient to
drop the $0$ mean normalization of $b(x)$ and consider $\hat{\alpha}$ having
the form%
\begin{equation}
\hat{\alpha}(x)=b(x)^{\prime}\hat{\rho}\text{,}  \label{Riesz est}
\end{equation}
where $\hat{\rho}$ is a vector of estimated coefficients.

The $\hat{\alpha}$ depends on the choice of dictionary $b(x)$ and penalty
degree $r.$ For the dictionary we require that each $b_{j}(x)$ belongs to
the set $\Gamma$ of possible plims of $\hat{\gamma}(x)$ discussed in Section
2 and that linear combinations of the dictionary "span" $\Gamma.$

\bigskip

\textsc{Assumption 1:} $b(x)=(b_{1}(x),...,b_{p}(x))^{\prime}$ where i) $%
b_{j}\in\Gamma$ for all $j$ and ii) for any $\alpha\in\Gamma$ and $%
\varepsilon>0$ there is $p$ and $\rho\in%
%TCIMACRO{\U{211d} }%
%BeginExpansion
\mathbb{R}
%EndExpansion
^{p}$ such that $\mathrm{E}[\{\alpha(X)-b(X)^{\prime}\rho\}^{2}]<%
\varepsilon. $

\bigskip

One key feature of this condition is that each $b_{j}\in\Gamma.$ This
feature allows us to use $m(w,\gamma)$ to construct $\hat{\alpha}$ and will
guarantee that $\hat{\alpha}\in\Gamma$, as required for the orthogonality
shown in equation (\ref{2nd order}). Another key feature is that linear
combinations of $b(x)$ can approximate anything that belongs to $\Gamma.$
This feature will lead to $\hat{\alpha}$ estimating $\bar{\alpha}.$ The link
imposed by Assumption 1, between the regression learner $\hat{\gamma}$ and
the dictionary $b(x)$ used to construct $\hat{\alpha},$ is important for the
orthogonality property of $\psi(w,\gamma,\alpha,\theta)$ and hence for $\hat{%
\theta}$ to be asymptotically normal and $\hat{V}$ to be a consistent
estimator of the asymptotic variance under general misspecification.

Assumption 1 requires that linear combinations of $b(x)$ must be able to
approximate any $\gamma$ in the set of possible plims of $\hat{\gamma}$ and
that each $b_{j}$ must be a possible plim of $\hat{\gamma}$. For Lasso and
other high dimensional regression learners where $X=(X_{1},X_{2},...)$
Assumption 1 will be satisfied for%
\begin{equation}
b(x)=(x_{1},...,x_{p})^{\prime}.  \label{high dimensional dict}
\end{equation}
Evidently each element $b_{j}(X)=X_{j}$ is an element of $\Gamma$ and the
spanning condition is satisfied because any linear combination of $X$ with a
finite number of nonzero coefficients will also be a linear combination of $%
b(x)$ for $p$ large.

We emphasize that $b(X)$ is required to approximate only the projection $%
\bar{\alpha}(X)$ and not $\alpha _{0}(X).$ For instance, in the average
treatment effect example $\bar{\alpha}(X)$ is the projection of the
difference of inverse propensity scores on the space spanned by $%
X=(X_{1},X_{2},...)$ which is naturally approximated by linear combinations
of $X=(X_{1},...,X_{p}).$ Assumption 1 does not require that this $b(X)$
approximate the inverse propensity score.

For neural nets, random forests, and other learners that nonparametrically
estimate $\mathrm{E}[Y|X],$ Assumption 1 will require that a linear
combination of $b(X)$ can approximate any function of $X$ for large enough $%
p.$ Such a $b(x)$ can be formed from low order multivariate powers of
components of $x$, with a full set of approximating functions included as $p$
grows. In applications one may use a variety of nonlinear functions
including powers of transformations of $X.$

The learner $\hat{\alpha}$ also depends on the choice of penalty degree $r.$
An important, useful feature of Lasso is that $r=A\sqrt{\ln(p)/n}$ for a
constant $A$ gives the fastest possible mean square convergence rate for
Lasso, that optimally trades off bias and variance. In Appendix~\ref%
{sec:computing}, we describe cross-validation and theoretical methods for
choosing the choosing $r$ based on data that have proven stable across
several different applications. We also provide R code, available upon
request, for the construction of $\hat{\alpha}(x)$ and $\hat{\theta}$.

We can motivate $\hat{\rho}$ in $\hat{\alpha}(x)=b(x)^{\prime}\hat{\rho}$ as
being based on the Riesz representation in equation (\ref{Riesz rep}) and $%
\bar{\alpha}$ satisfying equation (\ref{BLP alpha}), which imply that for $%
m(w,b)=(m(w,b_{1}),...,m(w,b_{p}))^{\prime}$,%
\begin{equation}
M:=\mathrm{E}[m(W,b)]=\mathrm{E}[\alpha_{0}(X)b(X)]=\mathrm{E}[\bar{\alpha}%
(X)b(X)],  \label{Rr basis}
\end{equation}
where the last equality is satisfied by $b_{j}\in\Gamma,$ which implies $%
\mathrm{E}[b_{j}(X)\{\alpha_{0}(X)-\bar{\alpha}(X)\}]=0$ for each $j$. We
see that the cross moments $M$ between the true, unknown $\bar{\alpha}(x)$
and the dictionary $b(x)$ are equal to the expectation of the known vector
of functions $m(w,b).$ Also, the second moment matrix $G=\mathrm{E}%
[b(X)b(X)^{\prime}]$ of the dictionary is an expectation of a known function
of the data. Estimating $M$ and $G$ enables learning coefficients $\rho$ of
the least squares regression of $\bar{\alpha}(X)$ on $b(X),$ satisfying $%
M=G\rho.$ We learn $\rho$ using a Lasso minimum distance objective function
to allow for large $p$. Let%
\begin{equation*}
\hat{M}=\frac{1}{n}\sum_{i=1}^{n}m(W_{i},b),\text{ }\hat{G}=\frac{1}{n}%
\sum_{i=1}^{n}b(X_{i})b(X_{i})^{\prime},
\end{equation*}
be unbiased estimators of $M$ and $G.$ The coefficient estimator is given by 
\begin{equation}
\hat{\rho}=\arg\min_{\rho}\{-2\hat{M}^{\prime}\rho+\rho^{\prime}\hat{G}%
\rho+2r\left\Vert \rho\right\Vert _{1}\},\text{ }\left\Vert \rho\right\Vert
_{1}=\sum_{j=1}^{p}|\rho_{j}|.  \label{RRLasso}
\end{equation}

The estimator $\hat{\rho}$ can be interpreted as a minimum distance version
of Lasso. Here $\hat{M}$ is analogous to $\sum_{i=1}^{n}Y_{i}b(X_{i})/n$ in
Lasso. The objective function in equation (\ref{RRLasso}) can be thought of
as the Lasso objective with $\sum_{i=1}^{n}Y_{i}b(X_{i})/n$ replaced by $%
\hat{M}$ and $\sum_{i=1}^{n}Y_{i}^{2}/n$ dropped. In this way the objective
function is a penalized approximation to the least squares regression of $%
\alpha_{0}(x)$ on $b(x),$ where $2r\left\Vert \rho\right\Vert _{1}$ is the
penalty. We refer to this as minimum distance Lasso because $\hat{M}$ does
not have the product form of Lasso regression.

The learner $\hat{\alpha}(x)$ of $\bar{\alpha}(x)$ is automatic in being
based on $\hat{M}$ and $\hat{G}$, neither of which requires knowledge of the
form of $\bar{\alpha}.$ In particular, $\hat{\alpha}(x)=b(x)^{\prime }\hat{%
\rho}$ does not depend on plugging in nonparametric estimates of components
of $\bar{\alpha}(x).$ Instead, $b(x)^{\prime }\hat{\rho}$ is linear in the
dictionary $b(x)$ and uses the known functional $m(w,\gamma )$ in the
construction of $\hat{M}$ to obtain the learner $\hat{\rho}.$ This automatic
nature of $\hat{\alpha}(x)$ is especially useful for Lasso and other high
dimensional regression learners where $b(x)$ can be taken to be the first $p$
elements of $x=(x_{1},x_{2},...),$ and where $\bar{\alpha}(x)$ is a least
squares projection of $\alpha _{0}(X)$ on $\Gamma ,$ as in Section 2. The
projection $\bar{\alpha}(x)$ will generally not have a simple form that can
be learned by plugging in nonparametric learners to an explicit formula. For
instance, in the average treatment effect example the projection of the
inverse propensity score on the high dimensional regressors $%
(X_{1},X_{2},...)$ does not have a closed form but is naturally approximated
by a linear combination of the first $p$ regressors where $%
b(X)=(X_{1},...,X_{p})^{\prime }.$

The learner $\hat{\alpha}(x)=b(x)^{\prime }\hat{\rho}$ also avoids inverting
a learner of a conditional probability or pdf. The finite sample properties
of methods that rely on inverses of learners can be poor; see Singh and Sun
(2019) for recent examples. Instead, $\hat{\alpha}$ approximates and learns $%
\bar{\alpha}$ by a linear combination of functions. In this way the $\hat{%
\alpha}$ that we propose here avoids potential instability from inverting a
high dimensional estimator. The inverse of a conditional probability or
density is present in $\alpha _{0}(x)$ in all of the examples in this paper.
We anticipate that this feature is present quite generally for causal and
structural models involving shifts in regressors, because the Rr equation (%
\ref{Riesz rep}) involves an expectation with respect to the data
distribution rather than the shifted distribution. Thus absence of an
inverse of a machine learner in $\hat{\alpha}$ may prove to be widely
useful. In some economic structural models the linearity of $\hat{\alpha}$
in $b(x)$ may not be quite as appealing, because inverse densities can have
a parametric form and so mitigate the problem of inverting a high
dimensional learner. An example is the dynamic discrete choice learner of
Chernozhukov et al. (2016, 2020). Also there is more work to be done to see
whether this approach has better properties than previously proposed ones in
practical settings.

This learner $\hat{\alpha}(x)$ can be thought of as being based on
orthogonality of the moment function with respect to $\gamma.$ Let $\tau$
denote a scalar and $b_{j}(x)$ an element of $b(x)$. Then by equation (\ref%
{Rr basis})%
\begin{equation*}
\frac{\partial}{\partial\tau}\mathrm{E}[\psi(W,\theta,\gamma+\tau b_{j},\bar{%
\alpha })]=\mathrm{E}[m(W,b_{j})-\bar{\alpha}(X)b_{j}(X)]=0,\text{ }%
(j=1,...,p).
\end{equation*}
Replacing the expectation by a sample average and $\bar{\alpha}(X)$ by $%
b(X)^{\prime}\rho$ gives 
\begin{equation*}
\frac{1}{n}\sum_{i=1}^{n}\{m(W_{i},b_{j})-[b(X_{i})^{\prime}%
\rho]b_{j}(X_{i})\}=e_{j}^{\prime}(\hat{M}-\hat{G}\rho),
\end{equation*}
where $e_{j}$ is the jth columin of a $p$ dimensional identity matrix. This
sample average is a scaled version of the derivative of objective function
in equation (\ref{RRLasso}) without the penalty term. The first-order
conditions for equation (\ref{RRLasso}) will set $\hat{\rho}$ so that this
object is close to zero, subject to the penalty, i.e. will solve penalized
versions of a moment equation. Thus, the Lasso minimum distance learner can
be thought of as a method that uses orthogonality of $\psi(W,\theta,\gamma,%
\alpha)$ with respect to $\gamma$ to learn $\bar{\alpha}$ while penalizing
to facilitate high dimensional estimation. In Section 6 we use an extension
of this approach to construct an Auto-DML when $m(W,\gamma)$ is nonlinear in 
$\gamma.$

To illustrate $\hat{\alpha}$ we consider the choice of dictionary and the
form of $\hat{\alpha}$ for Examples 1-3.

\bigskip

\textsc{Example 1:} If the regression learner $\hat{\gamma}$ is
nonparametric the dictionary $b(X)$ should also be nonparametric while if $%
\hat{\gamma}$ is a high dimensional regression the dictionary should be
chosen as in equation (\ref{high dimensional dict}). Here $m(w,b)=\int
b(x)[f_{1}(x)-f_{0}(x)]dx$ does not depend on the data observation $w$ and
the first order conditions for $\hat{\rho}$ imply that for each $j$,

\begin{equation*}
\left\vert \int b_{j}(x)[f_{1}(x)-f_{0}(x)]dx-\frac{1}{n-n_{\ell}}%
\sum_{i\notin I_{\ell}}b_{j}(X_{i})\hat{\alpha}_{\ell}(X_{i})\right\vert
\leq r.
\end{equation*}
Here $\hat{\alpha}_{\ell}(X_{i})$ acts to approximately re-weight so that
the integral of the basis function $b_{j}(x)$ over the policy shift is
approximately equal to the sample average of the re-weighted basis function $%
b_{j}(X_{i})\hat{\alpha}_{\ell}(X_{i}).$

\bigskip

\textsc{Example 2:} The dictionary $b(X)$ should be chosen as in Example 1.
Also by $m(w,b)=S(u)\gamma(u,z)$ the first order conditions for $\hat{\rho}$
imply that for each $j$,

\begin{equation*}
\left\vert \frac{1}{n-n_{\ell}}\sum_{i\notin
I_{\ell}}\{S(U_{i})b_{j}(U_{i},Z_{i})-b_{j}(X_{i})\hat{\alpha}%
_{\ell}(X_{i})\}\right\vert \leq r.
\end{equation*}
Here $\hat{\alpha}_{\ell}(X_{i})$ acts approximately as a re-weighting
scheme, making the sample average of the score $S(U_{i})$ times the basis
function $b_{j}(U_{i},Z_{i})$ be approximately equal to the sample average
of the re-weighted basis function $b_{j}(X_{i})\hat{\alpha}_{\ell}(X_{i}).$

\bigskip

\textsc{Example 3: }The dictionary should be chosen similarly to Example 1.
For instance suppose that $X=(DZ,(1-D)Z)$, where $Z=(Z_{1},Z_{2},...)$ is a
sequence or possible covariates. Then the dictionary%
\begin{equation}
b(x)=(dq(z)^{\prime},(1-d)q(z)^{\prime})^{\prime},\text{ }%
q(z)=(z_{1},...,z_{p/2})^{\prime},  \label{ATE dic}
\end{equation}
would satisfy Assumption 1. The estimator $\hat{\alpha}_{\ell}$ has an
interesting form for this dictionary. Note that $%
m(w,b)=b(1,z)-b(0,z)=(q(z)^{\prime},0^{\prime})^{\prime}-(0^{\prime},q\left(
z\right) ^{\prime})^{\prime}=(q(z)^{\prime},-q(z)^{\prime})$. Then%
\begin{equation*}
\hat{M}_{\ell}=\left( 
\begin{array}{c}
\bar{q}_{\ell} \\ 
-\bar{q}_{\ell}%
\end{array}
\right) ,\text{ }\bar{q}_{\ell}=\frac{1}{n-n_{\ell}}\sum_{i\notin
I_{\ell}}q(Z_{i}).
\end{equation*}
Let $\hat{\rho}_{\ell}^{1}$ be the estimated coefficients of $dq(z)$ and $%
\hat{\rho}_{\ell}^{0}$ be the estimated coefficients of $(1-d)q(z)$. Then
the learner of $\bar{\alpha}(X_{i})$ is 
\begin{equation*}
\hat{\alpha}_{\ell}(X_{i})=D_{i}\hat{\omega}_{\ell i}^{1}-(1-D_{i})\hat {%
\omega}_{\ell i}^{0},\text{ }\hat{\omega}_{\ell i}^{1}=q(Z_{i})^{\prime}\hat{%
\rho}_{\ell}^{1},\text{ }\hat{\omega}_{\ell i}^{0}=-q(Z_{i})^{\prime}\hat{%
\rho}_{\ell}^{0},
\end{equation*}
where $\hat{\omega}_{\ell i}^{1}$ and $\hat{\omega}_{\ell i}^{0}$ might be
thought of as \textquotedblleft weights.\textquotedblright\ These weights
sum to one if $q(z)$ includes a constant but may be negative. The first
order conditions for $\hat{\alpha}$ are that for each $j,$%
\begin{equation}
\left\vert \frac{1}{n-n_{\ell}}\sum_{i\notin I_{\ell}}q_{j}(Z_{i})[1-D_{i}%
\hat{\omega}_{\ell i}^{1}]\right\vert \leq r,\text{ }\left\vert \frac {1}{%
n-n_{\ell}}\sum_{i\notin I_{\ell}}q_{j}(Z_{i})[1+(1-D_{i})\omega_{\ell
i}^{0}]\right\vert \leq r.  \label{ATE balance}
\end{equation}
Here $\hat{\rho}_{\ell}$ sets the weights $\hat{\omega}_{\ell i}^{1}$ and $%
\hat{\omega}_{\ell i}^{0}$ to approximately \textquotedblleft
balance\textquotedblright\ the overall sample average with the treated and
untreated averages for each element of the dictionary $q(z).$ The
constraints of equation (\ref{ATE balance}) are like the balancing
conditions of Zubizarreta (2015) and Athey, Imbens, and Wager (2018). The
source of these constraints is regularized least squares approximation of $%
\bar{\alpha }(x)=proj(\pi_{0}(z)^{-1}d-[1-\pi_{0}(z)]^{-1}(1-d)|Z)$ by a
linear combination of the dictionary $b(x)$. The approach of this paper
shows that this type of balancing is sufficient to debias any regression
learner under regularity conditions in Section 4.

\section{Large Sample Inference}

In this Section, we give mean square convergence rates for the Lasso minimum
distance learner of $\hat{\alpha}$ and root-n consistency and asymptotic
normality results for the learner $\hat{\theta}$ of the object of interest
and its asymptotic variance estimator $\hat{V}$. Let $\varepsilon_{n}$
denote a sequence that converges to zero no faster than $\sqrt{\ln(p)/n}$
and for a random variable $a(W)$ let $\left\Vert a\right\Vert =\sqrt{\mathrm{%
E}[a(W)^{2}]}$

\bigskip

\textsc{Assumption 2:} \textit{There exists }$C>1,$ $\xi>0$\textit{\ such
that for each positive integer }$s\leq C\varepsilon_{n}^{-2/(2\xi+1)}$%
\textit{\ there is }$\bar{\rho}$\textit{\ with }$s$\textit{\ nonzero
elements such that}%
\begin{equation*}
\left\Vert \bar{\alpha}-b^{\prime}\bar{\rho}\right\Vert \leq C(s)^{-\xi}.
\end{equation*}

\bigskip

Here $\left\Vert \bar{\alpha}-b^{\prime}\bar{\rho}\right\Vert $ is the mean
square approximation error from using the linear combination $b^{\prime}\bar{%
\rho}$ to approximate $\bar{\alpha}.$ This approximate sparsity condition
specifies that there is a sparse $\bar{\rho}$, having only $s$ nonzero
elements, so that the approximation error is bounded by $C(s)^{-\xi}.$ Note
that it is not required that $\bar{\alpha}$ be equal to linear combination
of $s$ terms, i.e. it is not required that $\bar{\alpha}$ be strictly
sparse. Assumption 2 does allow unknown identity of the elements of $b(x)$
that give the approximation rate $s^{-\xi}$. In this way this condition
allows for high dimensional $x$ where statistics and economics do not
provide much guidance on which elements of $b(x)$ are important.

The $\varepsilon_{n}$ in this condition represents a convergence rate for $%
\hat{M}$ and $\hat{G}$ that will be no faster than $\sqrt{\ln(p)/n}$ under
the conditions given in the rest of this Section. When $s$ is chosen to be
approximately $C\varepsilon_{n}^{-2/(2\xi+1)},$ which is the largest $s$
allowed by Assumption 2, $s$ will grow no faster than $(\sqrt{n/\ln (p)}%
)^{2/(2\xi+1)}\leq n^{1/(2\xi+1)},$ which grows slower than $n.$ Because $%
p\geq s$ is implicitly required by this condition, Assumption 2 puts a quite
a weak restriction on $p.$ An important feature of Assumption 2 is that the
sparse approximation is based on functions included in the $p\times1$
dictionary $b(x)$. Thus larger values of $p$ give more flexibility and will
help Assumption 2 to be satisfied.

Our results will require a convergence rate for $\hat{\alpha}$ that is
faster than some power of $n.$ Assumption 2 is a natural condition that
leads to such a rate. Sufficient conditions for Assumption 2 are well known
from the approximation literature when $\bar{\alpha}(x)$ belongs to a Besov
or Holder class of function and linear combinations of $b(x)$ can
approximate any function of $x$.

We will also make use of a sparse eigenvalue condition as considered in much
of the Lasso literature. Let $\rho$ denote a $p\times1$ vector, $\rho_{J}$ a 
$J\times1$ subvector of $\rho,$ and $\rho_{J^{c}}$ the vector consisting of
components of $\rho$ that are not in $\rho_{J}$. Also for a matrix $A$ let $%
\left\Vert A\right\Vert _{1}=\sum_{i,j}\left\vert a_{ij}\right\vert .$

\bigskip

\textsc{Assumption 3: }$G=\mathrm{E}[b(X)b(X)^{\prime}]$ \textit{has largest
eigenvalue bounded uniformly in }$n$\textit{\ and there is }$C,c>0$\textit{\
such that for all }$s\approx C\varepsilon_{n}^{-2}$\textit{\ with
probability approaching one} 
\begin{equation*}
\min_{J\leq s}\min_{\left\Vert \rho_{J^{c}}\right\Vert _{1}\leq3\left\Vert
\rho_{J}\right\Vert _{1}}\frac{\rho^{\prime}\hat{G}\rho}{\rho_{J}^{\prime}%
\rho_{J}}\geq c
\end{equation*}

\bigskip

This is a sparse eigenvalue condition that is familiar from the Lasso
literature, including Bickel, Ritov, Tsybakov (2009), Belloni and
Chernozhukov (2013), and Rudelson and Zhou (2013).

We will work with a dictionary $b(X)$ with elements that are uniformly
bounded.

\bigskip

\textsc{Assumption 4}: \textit{There is }$C>0$\textit{\ such that with
probability one }$\sup_{j}|b_{j}(X)|\leq C.$

\bigskip

This condition implies a convergence rate of $\sqrt{\ln(p)/n}$ for $%
\left\Vert \hat{G}-G\right\Vert _{\infty},$ where $\left\Vert A\right\Vert
_{\infty}=\max_{i,j}\left\vert a_{ij}\right\vert $ for a matrix $A=[a_{ij}]$.

Lasso mean square convergence rates are often stated in terms of finite
sample bounds. Because the focus of this paper is root-n consistency for $%
\hat {\theta}$ and for that we only need convergence at certain powers of $n$
we can simplify the statement of convergence rates without affecting the
conditions for $\hat{\theta}$ by allowing the Lasso regularization value $r$
to shrink slightly slower than $\varepsilon_{n}.$ This does lead to
approximate sparseness conditions that are strict inequalities on the size
of $\xi$ but Bradic et al. (2019) have shown that strict inequalities are
necessary for root-n consistent estimation, meaning that there is no loss of
generality in these conditions. We also limit the growth of $p$ to be slower
than some power of $n.$

\bigskip

\textsc{Assumption 5:} $\varepsilon_{n}=o(r),$ $r=o(n^{c}\varepsilon_{n})$
for all $c>0$, and there exists $C>0$ such that $p\leq Cn^{C}.$

\bigskip

We also hypothesize a convergence rate for $\hat{M}.$

\bigskip

\textsc{Assumption 6}: $\left\Vert \hat{M}-M\right\Vert
_{\infty}=O_{p}(\varepsilon_{n})$ for $\varepsilon_{n}\longrightarrow0.$

\bigskip

We use this condition to accommodate $\hat{M}$ that can depend on the
regression learner $\hat{\gamma}$ as needed for Section 5.

\bigskip

\textsc{Theorem 1:} \textit{If Assumptions 1 - 6 are satisfied then for all }%
$c>0,$%
\begin{equation*}
\Vert\hat{\alpha}-\bar{\alpha}\Vert=o_{p}(n^{c}\varepsilon_{n}^{2\xi/(2\xi
+1)}).
\end{equation*}

\bigskip

This theorem is based on extending Lemmas of Bradic et al. (2019) to allow $%
\varepsilon_{n}$ to shrink slower than $\sqrt{\ln(p)/n}.$ The extension will
be used in Section 5 to obtain convergence rates when $\hat{M}$ depends on a
nonparametric estimator.

The sparse eigenvalue condition of Assumption 3 seems strong in some
settings. It is possible to drop Assumption 3 and Assumption 2 if the
following condition is satisfied:

\bigskip

\textsc{Assumption 7: }$\bar{\alpha}(X)=\sum_{j=1}^{\infty}\rho_{j0}b_{j}(X)$%
\textit{, }$\sum_{j=1}^{\infty}\left\vert \rho_{j0}\right\vert <\infty $%
\textit{, and for }$C>0$\textit{\ and }$\bar{s}=C\sqrt{n}$\textit{\ the }$%
b_{j}(x)$\textit{\ corresponding to the largest }$\bar{s}$\textit{\ values
of }$\left\vert \rho_{j0}\right\vert $\textit{\ are included in }$b(x).$

\bigskip

This condition allows us to drop Assumption 2 because absolute summability
of the coefficients $\rho_{0j}$ implies a sparse approximation rate of $%
\xi=1/2.$ It also allows $\hat{G}$ to converge at a rate slower $%
\varepsilon_{n}$ in order to accommodate nonparametric estimation in $\hat{G}%
.$

\bigskip

\textsc{Theorem 2:} \textit{If Assumptions 1 and 5-7 are satisfied and }$%
\left\Vert \hat{G}-G\right\Vert _{\infty}=O_{p}(\varepsilon_{n})$\textit{\
then for all }$c>0,$%
\begin{equation*}
\Vert\hat{\alpha}-\bar{\alpha}\Vert=o_{p}(n^{c}\sqrt{\varepsilon_{n}}).
\end{equation*}

\bigskip

This result extends Chatterjee and Javarov (2015) to allow $\varepsilon_{n}$
to shrink slower than $\sqrt{\ln(p)/n}.$ When $\varepsilon_{n}=\sqrt{\ln
(p)/n}$ in Assumption 6 this result gives a mean square convergence rate for 
$\hat{\alpha}$ that is faster than $n^{-1/4+c}$ for all $c>0,$ without a
sparse eigenvalue condition.

We now use these results to obtain root-n consistency and asymptotic
normality for the Auto-DML $\hat{\theta}$ and consistency of its asymptotic
variance estimator $\hat{V}.$ We impose some additional regularity
conditions.

\bigskip

\textsc{Assumption 8}: \textit{There is }$C>0$\textit{\ such that with
probability one }$\max_{j\leq p}|m(W,b_{j}))|\leq C.$

\bigskip

Under this condition Assumption 6 will be satisfied with $\varepsilon _{n}=%
\sqrt{\ln(p)/n}.$ This condition will be satisfied under by Assumption 4 in
each of Examples 1-3 under conditions of Corollaries 4-6 to follow.

\bigskip

\textsc{Assumption 9:}\textit{\ }$\mathrm{E}[\{Y-\bar{\gamma}(X)\}^{2}|X]$%
\textit{\ and }$\bar{\alpha}(X)$\textit{\ are bounded.}

\bigskip

We impose this condition for simplicity; it could be weakened. We also
impose the following condition.

\bigskip

\textsc{Assumption 10:} $\mathrm{E}[m(W,\gamma_{0})^{2}]<\infty$ \textit{and 
}$\int[m(w,\hat{\gamma})-m(w,\bar{\gamma})]^{2}F_{W}(dw)\overset{p}{%
\longrightarrow}0.$

\bigskip

This condition will be implied by existence of $C>0$ with $\left\vert 
\mathrm{E}[m(W,\gamma)^{2}]\right\vert \leq C\left\Vert \gamma\right\Vert
^{2}$ for all $\gamma$, which will be satisfied in the examples we consider
under regularity conditions to be specified.

\bigskip

\textsc{Assumption 11: }\textit{With probability approaching one }$\hat {%
\gamma}_{\ell}\in\Gamma$\textit{\ and there is }$d_{\gamma}>0$\textit{\ such
that }$\Vert\hat{\gamma}-\bar{\gamma}\Vert=O_{p}(n^{-d_{\gamma}})$\textit{\
and either Assumptions 2 and 3 are satisfied with}%
\begin{equation}
\frac{\xi}{2\xi+1}+d_{\gamma}>\frac{1}{2},  \label{rate dr}
\end{equation}
\textit{or Assumption 7 is satisfied and }$d_{\gamma}>1/4.$

\bigskip

This assumption allows $\hat{\gamma}$ to be any learner that converges in
mean square at a rate that is some power of $n.$ By Theorem 1, the mean
square convergence rate for $\hat{\alpha}$ is as close as desired to $%
n^{-\xi /(2\xi+1)}.$ Thus Assumption 11 requires that the product of
convergence rates for $\hat{\alpha}$ and $\hat{\gamma}$ must go to zero
faster than $1/\sqrt {n}.$ This is a rate double robustness condition that
appears in earlier low dimensional and high dimensional literatures cited in
the introduction. Under Assumptions 2 and 3 a full trade-off in rates
between $\hat{\alpha}$ and $\hat{\gamma}$ is permitted, since Assumption 11
is satisfied for any $\xi$ if $d_{\gamma}$ is large enough and for any $%
d_{\gamma}$ if $\xi$ is large enough. Under Assumption 7 this trade-off is
not present, since $d_{\gamma }>1/4$ is required by Assumption 11. Assumption 11 can be dropped if $\alpha_0(X)$ is known and is used in place of $\hat\alpha(X)$ in the construction of $\hat\theta$ in equation (3.1). In that case only mean square consistency of $\hat\gamma$ will be required for root-n consistency and asymptotic normality of $\hat\theta$.

The following gives the large sample inference results for $\hat{\theta}$
and $\hat{V}.$ Define%
\begin{equation*}
\bar{\theta}=\mathrm{E}[m(W,\bar{\gamma})],\text{ }\psi(w)=m(w,\bar{\gamma})-%
\bar{\theta}+\bar{\alpha}(x)[y-\bar{\gamma}(x)],\text{ }V=\mathrm{E}%
[\psi(W)^{2}].
\end{equation*}
Here $\bar{\theta}$ will be the object estimated by $\hat{\theta}$ when
neither of the double robustness conditions $\bar{\gamma}(X)=\mathrm{E}[Y|X]$
nor $\bar{\alpha}(X)\in\Gamma$ is satisfied.

\bigskip

\textsc{Theorem 3}: \textit{If Assumptions 1-5, and 8-11 are satisfied then }%
$\sqrt{n}(\hat{\theta}-\bar{\theta})\overset{d}{\longrightarrow}N(0,V).$ 
\textit{If in addition Assumption 7 is satisfied then }$\hat{V}\overset{p}{%
\longrightarrow}V$.

\bigskip

It is possible to construct a consistent estimator of $V$ without Assumption
7 by using a trimmed version of $\hat{\alpha}_{\ell }(x)$ but we omit that
demonstration to avoid further complicating $\hat{V}$. The conclusion of
Theorem 3 implies that asymptotic test statistics and confidence intervals
can be formed in the usual manner from $\hat{\theta}$ and $\hat{V}.$ Theorem
3 is proven by using the convergence rate results of Theorem 1 and Theorem 2
to show that the hypotheses of Lemma 15 of Chernozhukuv et al. (2020) are
satisfied.

The asymptotic variance $V$ is fixed rather than varying with $n$ because we
have chosen to work with i.i.d. data and an approximately sparse regression
for simplicity. It would be straightforward to extend the results to allow
the regression to change with sample size in order to accomodate sparse
regressions and corresponding variances that change with $n$.

Under similar conditions as Theorem 3 Auto-TML is also consistent and
asymptotically normal.

\bigskip 

\textsc{Corollary 4:} \textit{If Assumptions 1-5, and 8-11 are satisfied, }$%
\mathrm{E}[m(W,\gamma )^{2}]\leq C\left\Vert \gamma \right\Vert^{2} $ \textit{for all 
}$\gamma \in \Gamma ,$\textit{\ and }$\bar{\alpha}(X)\neq 0$\textit{\ then }$%
\sqrt{n}(\tilde{\theta}-\bar{\theta})\overset{d}{\longrightarrow }N(0,V).$

\bigskip 

Most of the conditions of Theorem 3 are quite general, with only Assumptions
8 and 10 pertaining to a particular $m(w,\gamma )$. It is straightforward to
specify conditions under which Assumptions 8 and 10 are satisfied for
Examples 1-3.

\bigskip

\textsc{Corollary 5 (Example 1):} \textit{If Assumptions 1-5, 9, and 11 are
satisfied and there is }$C>0$\textit{\ such that }$\left\vert
[f_{1}(x)-f_{0}(x)]/f(x)\right\vert \leq C$ \textit{then }$\sqrt{n}(\hat{%
\theta}-\bar{\theta})\overset{d}{\longrightarrow }N(0,V).$\textit{\ If in
addition Assumption 7 is satisfied then }$\hat{V}\overset{p}{\longrightarrow 
}V$\textit{.}

\bigskip

The specific regularity condition for the policy effect in Corollary 5 is
that the Rr $\alpha _{0}(X)=[f_{1}(X)-f_{0}(X)]/f(x)$ be bounded.

\bigskip

\textsc{Corollary 6 (Example 2):} \textit{If Assumptions 1-5, 9, and 11 are
satisfied and there is }$C>0$\textit{\ such that }$\left\vert
S(u)\right\vert \leq C$, $f(D|Z)^{-1}\omega (D)\leq C$\textit{\ then }$\sqrt{%
n}(\hat{\theta}-\bar{\theta})\overset{d}{\longrightarrow }N(0,V).$\textit{\
If in addition Assumption 7 is satisfied then }$\hat{V}\overset{p}{%
\longrightarrow }V$\textit{.}

\bigskip

The regularity conditions for the weighted average derivative in Corollary 6
are that the score $S(u)$ is bounded and the Rr $\alpha
_{0}(X)=f(D|Z)^{-1}\omega (D)S(D)$ is also bounded.

\bigskip

\textsc{Corollary 7 (Example 3): }\textit{If Assumptions 1, 4-5, 9, and 11
are satisfied and there is }$C>0$\textit{\ with }$\pi _{0}(Z)\in \lbrack
C,1-C]$\textit{\ then }$\sqrt{n}(\hat{\theta}-\bar{\theta})\overset{d}{%
\longrightarrow }N(0,V).$\textit{\ If in addition Assumption 7 is satisfied
then }$\hat{V}\overset{p}{\longrightarrow }V$\textit{.}

\bigskip

The additional condition in Corollary 7 is that the propensity score is
bounded away from $0$ and $1$, an overlap condition that is common in
asymptotic theory for estimators of the average treatment effect. Together
Corollaries 5--7 demonstrate how simple primitive conditions involving $%
m(w,\gamma )$ can be specified so that the Auto-DML $\hat{\theta}$ of an
object of interest will be asymptotically normal and the asymptotic variance
estimator $\hat{V}$ consistent.

\section{Nonlinear Effects of Multiple Regressions}

Some important effects of interest are expectations of nonlinear functions
of multiple regressions. Causal mediation analysis is an important example
that we consider in this Section. The regression decomposition in Section 6
is another important example. In this Section we give Auto-DML for such
effects. Such effects have the form $\theta_{0}=\mathrm{E}[m(W,\gamma_{0})]$
where $m(w,\gamma)$ is nonlinear in a possible value $\gamma$ of multiple
regressions $(\gamma _{1}(X_{1}),...,\gamma_{K}(X_{K}))^{\prime}$ with
regressors $X_{k}$ specific to each regression $\gamma_{k}(X_{k})$. The
corresponding orthogonal moment functions are like those discussed in
Section 3 except that the bias correction is a sum of $K$ terms with the $%
k^{th}$ term being the bias correction for the learner of $\gamma_{k}$, as
in Newey (1994, p. 1357). The estimated bias corrections are like those of
Section 4 with the $k^{th}$ term being the product of a Lasso learner $\hat{%
\alpha}_{k\ell}(X_{k})$ and the residual $Y_{k}-\hat{\gamma}_{k\ell}(X_{k}).$
Each $\hat{\alpha}_{k\ell}(X_{k})$ differs from Section 3 in the
corresponding $\hat{M}_{k\ell}$ being a derivative evaluated at a
preliminary estimator of $\bar{\gamma}$. Because the construction of $\hat{%
\theta}$ is so closely related to that in Section 3 we proceed immediately
with its description here and fill in details concerning the orthogonal
moment function below.

The Auto-DML of a nonlinear effect is similar to equation (\ref{Estimator}).
Specifically it is%
\begin{align}
\hat{\theta}& =\frac{1}{n}\sum_{\ell =1}^{L}\sum_{i\in I_{\ell }}\{m(W_{i},%
\hat{\gamma}_{\ell })+\sum_{k=1}^{K}\hat{\alpha}_{k\ell }(X_{ki})[Y_{ki}-%
\hat{\gamma}_{k\ell }(X_{ki})]\},  \label{nonlin est} \\
\hat{V}& =\frac{1}{n}\sum_{\ell =1}^{L}\sum_{i\in I_{\ell }}\hat{\psi}%
_{i\ell }^{2},\text{ }\hat{\psi}_{i\ell }=m(W_{i},\hat{\gamma}_{\ell })-\hat{%
\theta}+\sum_{k=1}^{K}\hat{\alpha}_{k\ell }(X_{ki})[Y_{ki}-\hat{\gamma}%
_{k\ell }(X_{ki})],  \notag
\end{align}%
where each $\hat{\alpha}_{k\ell }(X_{ki})$ is obtained as follows: For each $%
k$ let $b_{k}(x_{k})=(b_{k1}(x_{k}),....,b_{kp}(x_{k}))^{\prime }$ be a $%
p\times 1$ dictionary vector specific to the $k^{th}$ regression $\gamma
_{k}(x_{k})$ and let $\hat{\gamma}_{\ell ,\ell ^{\prime }}$ be the vector of
regressions computed from all observations not in either $I_{\ell }$ or $%
I_{\ell ^{\prime }}$. Also let $\tau $ denote a scalar, and $e_{k}$ the $%
k^{th}$ column of the $K$ dimensional identity matrix. Then%
\begin{align}
\hat{\alpha}_{k\ell }(X_{ki})& =b_{k}(X_{ki})^{\prime }\hat{\rho}_{k\ell },%
\text{ }\hat{\rho}_{k\ell }=\arg \min_{\rho }\{-2\hat{M}_{k\ell }^{\prime
}\rho +\rho ^{\prime }\hat{G}_{k\ell }\rho +2r_{k}\left\Vert \rho
\right\Vert _{1}\},\text{ }\left\Vert \rho \right\Vert
_{1}=\sum_{j=1}^{p}|\rho _{j}|,  \label{nonlin Rr} \\
\hat{M}_{k\ell }& =(\hat{M}_{k\ell 1},...,\hat{M}_{k\ell p})^{\prime },\text{
}\hat{G}_{k\ell }=\left( \frac{1}{n-n_{\ell }}\right) \sum_{i\notin I_{\ell
}}b_{k}(X_{ki})b_{k}(X_{ki})^{\prime },  \notag \\
\hat{M}_{k\ell j}& =\left. \frac{d}{d\tau }\left( \frac{1}{n-n_{\ell }}%
\right) \sum_{\ell ^{\prime }\neq \ell }\sum_{i\in I_{\ell ^{\prime
}}}m(W_{i},\hat{\gamma}_{\ell ,\ell ^{\prime }}+\tau e_{k}b_{kj})\right\vert
_{\tau =0},\text{ }(j=1,...,p).  \notag
\end{align}%
where $b_{kj}$ denotes the $j^{th}$ element of the dictionary $b_{k}(x_{k})$
as a function of $x_{k}.$ Thus the $\hat{\alpha}_{k\ell }(X_{i})$ in
equation (\ref{nonlin est}) is a Lasso minimum distance estimator like that
of Section 3 that is specific to $\hat{\gamma}_{k}$ and uses the $\hat{M}%
_{k\ell }$ from equation (\ref{nonlin Rr}) rather than the one in equation (%
\ref{CFLasso}).

The $\hat{M}_{k\ell j}$ given here generalizes equation (\ref{CFLasso}) to
allow for nonlinearity of $m(w,\gamma)$ in $\gamma.$ The derivative with
respect to the scalar $\tau$ in $\hat{M}_{k\ell j}$ is generally simple to
compute analytically using the chain rule of calculus, as we will illustrate
for causal mediation analysis. When $m(w,\gamma)$ is linear in a single $%
\gamma$ this derivative just evaluates $m(W_{i},\gamma)$ at $\gamma=b_{j}$,
giving the $\hat{M}_{\ell j}$ of equation (\ref{CFLasso}). As with linear $%
m(w,\gamma)$ the $\hat{M}_{k\ell j}$ and the rest of the $\hat{\theta}$
depends just on $m(w,\gamma)$ and the first step. Thus the $\hat{\theta}$ in
equation (\ref{nonlin est}) is automatic, in the same way as the estimator
of equation (\ref{Estimator}), in only requiring $m(w,\gamma)$ and the
regression residuals $Y_{{}}$for its construction.

The $\hat{M}_{k\ell j}$ given here does depend on a cross-fit regression
learner $\hat{\gamma}_{\ell ,\ell ^{\prime }}$ in order to allow for the
nonlinearity of $m(w,\gamma )$ in $\gamma .$ The cross-fitting will make the
sample average used in the construction of $\hat{M}_{k\ell j}$ independent
of the regression learner $\hat{\gamma}_{\ell ,\ell ^{\prime }}$ used in its
construction. This independence helps $\hat{M}_{k\ell j}$ to be uniformly
consistent over $j=1,...,p$ for large $p$ with only mean square convergence
convergence rates for $\hat{\gamma}_{\ell ,\ell ^{\prime }}.$ This feature
of the theory helps $\hat{\theta}$ to be root-n consistent and
asymptotically normal for a wide variety of regression learners $\hat{\gamma}%
_{\ell ,\ell ^{\prime }}.$ This $\hat{M}_{k\ell j}$ was given in
Chernozhukov, Newey, and Singh (2018, p. 17). Multiple cross-fitting has
also been used in Newey and Robins (2018) and Kennedy (2020).

The dictionary $b_{k}(x_{k})$ used in the construction of $\hat{\alpha}%
_{k\ell}(x_{k})$ should be chosen analogously to the $b(x)$ in Section 3.
Each $b_{kj}$ should be an element of the set $\Gamma_{k}$ of possible
plim's of $\hat{\gamma}_{k}$. Also linear combinations of $b_{k}(x_{k})$
should be able to approximate any element of $\Gamma_{k}$ arbitrarily well
in mean square. That is, Assumption 1 should be satisfied with $\Gamma_{k}$
and $b_{k}(x)$ replacing $\Gamma$ and $b(x)$ respectively. In particular if $%
\hat{\gamma}_{k}$ is a high dimensional regression then $%
b(x)=(x_{k1},...,x_{kp})^{\prime }$ will do. If $\hat{\gamma}_{k}$ is a
nonparametric estimator then $b_{k}(x_{k})$ should be chosen so that linear
combinations can approximate any function of $x_{k}$.

An important difference between the Lasso minimum distance learner in
Section 3 and each $\hat{\alpha}_{k\ell}(x_{k})$ here is that the penalty
size $r_{k}$ must be chosen to be larger than $\sqrt{\ln(p)/n}$ when $%
m(w,\gamma)$ depends nonlinearly on $\gamma.$ The reason for larger $r_{k}$
is that $\hat{M}_{k\ell}$ depends on the machine learner $\hat{\gamma}%
_{\ell,\ell^{\prime}}$ and so will converge at a slower rate, leading to a
requirement that $r_{k}$ converge to zero slightly slower than the mean
square convergence rate of $\hat{\gamma}_{\ell,\ell^{\prime}}.$ A choice of $%
r_{k}$ proportional to $n^{-1/4}$ will generally suffice for this purpose,
since $\hat{\gamma}_{\ell,\ell^{\prime}}$ will be required to converge
faster than $n^{-1/4}$.

This estimator will not be doubly robust due to the nonlinearity of $%
m(w,\gamma)$ in $\gamma$; see Chernozhukov et al. (2016). Nevertheless it
will have zero first order bias and so be root-n consistent and
asymptotically normal under sufficient regularity conditions. It has zero
first order bias because $\hat{\alpha}_{k\ell}(x_{k})$ will consistently
estimate $\bar{\alpha }_{k}(x_{k})$ such that $\sum_{k=1}^{K}\bar{\alpha}%
_{k}(x)[y_{k}-\bar{\gamma }_{k}(x_{k})]$ is the influence function for $%
\mathrm{E}[m(W,\gamma(F))]$ at $\gamma(F)=\bar{\gamma}$ where $\gamma(F)=$%
plim$(\hat{\gamma})$.

\bigskip

\textsc{Example 5:} (Causal Mediation Analysis) Causal mediation analysis
provides an interesting example of a nonlinear function of multiple
regressions. This effect allows for intermediate variables, called
mediators, that lie between treatment and outcome. In this example there is
an outcome variable $Y$, a treatment indicator $D\in\{0,1\},$ and covariates 
$Z$ similar to the average treatment effect in Example 3. In addition there
is a mediation variable that we will denote by $Q,$ where we assume that $%
Q\in\{1,...,K-1)$ for an integer $K\geq3.$ Let 
\begin{equation*}
\gamma_{K0}(D,Q,Z)=\mathrm{E}[Y|D,Q,Z]\text{, }\gamma_{k0}(D,Z)=\Pr
(Q=k|D,Z)=\mathrm{E}[1(Q=k)|D,Z],\text{ }(k=1,...,K-1).
\end{equation*}
The causal mediation effect of Imai, Keele, and Tingley (2010, Theorem 1) is 
\begin{equation*}
\theta_{0}(d,d^{\prime})=\mathrm{E}[\sum_{k=1}^{K-1}\gamma_{K0}(d,k,Z)\gamma
_{k0}(d^{\prime},Z)].
\end{equation*}
This effect, or parameter, has the form $\theta_{0}(d,d^{\prime})=\mathrm{E}%
[m(W,\gamma)]$ for $W=(Y,D,Q,Z)$ and 
\begin{equation*}
m(W,\gamma)=\sum_{k=1}^{K-1}\gamma_{K}(d,k,Z)\gamma_{k}(d^{\prime},Z).
\end{equation*}

In this example we have $X_{k}=(D,Z)$, $(k=1,...,K-1)$ and $X_{K}=(D,Q,Z).$
To construct the Auto-DML $\hat{\theta}$ we need to choose the dictionaries $%
b_{k}(X_{k})$ for each $k$. We choose%
\begin{equation*}
b_{K}(X_{K})=(b_{K1}(D,Q,Z),...,b_{Kp}(D,Q,Z))^{\prime}
\end{equation*}
to be a nonparametric dictionary if $\hat{\gamma}_{K}$ is a nonparametric
estimator such as a neural net or random forest or choose $b_{K}(D,Q,Z)$ to
be the leading $p$ regressors used in a high dimension regression learner $%
\hat{\gamma}_{K}.$ For $k\leq K-1$ we choose the same dictionary $%
b_{k}(X_{k})=b_{1}(D,Z)$ with%
\begin{equation*}
b_{1}(D,Z)=(b_{11}(D,Z),...,b_{1p}(D,Z))^{\prime},
\end{equation*}
for each $k\leq K-1.$ We specify $b_{1}(D,Z)$ to be a nonparametric
dictionary if each $\hat{\gamma}_{k}$ is a nonparametric estimator such as a
neural net or random forest or choose $b_{1}(D,Z)$ to be the leading $p$
regressors used in a high dimension regression learner for each $\hat{\gamma}%
_{k}.$

It is straightforward to compute each $\hat{M}_{k\ell j}.$ Note that for $%
k\leq K-1,$%
\begin{align*}
\left. \frac{d}{d\tau }m(W,\gamma +\tau e_{k}b_{kj})\right\vert _{\tau =0}&
=\left. \frac{d}{d\tau }\gamma _{K}(d,k,Z)\{\gamma _{k}(d^{\prime },Z)+\tau
b_{1j}(d^{\prime },Z)\}\right\vert _{\tau =0}=\gamma
_{K}(d,k,Z)b_{1j}(d^{\prime },Z), \\
\left. \frac{d}{d\tau }m(W,\gamma +\tau e_{K}b_{Kj})\right\vert _{\tau =0}&
=\left. \frac{d}{d\tau }\{\sum_{k=1}^{K-1}\{\gamma _{K}(d,k,Z)+\tau
b_{Kj}(d,k,Z)\}\gamma _{k}(d^{\prime },Z)]\right\vert _{\tau =0} \\
& =\sum_{k=1}^{K-1}b_{Kj}(d,k,Z)\gamma _{k}(d^{\prime },Z).
\end{align*}%
Then we have

\begin{align*}
\hat{M}_{k\ell j} & =\frac{1}{n-n_{\ell}}\sum_{\ell^{\prime}\neq\ell}\sum_{i%
\in I_{\ell^{\prime}}}\hat{\gamma}_{K\ell,\ell^{%
\prime}}(d,k,Z_{i})b_{1j}(d^{\prime},Z_{i}),\text{ }(k=1,...,K-1), \\
\hat{M}_{K\ell j} & =\frac{1}{n-n_{\ell}}\sum_{\ell^{\prime}\neq\ell}\sum_{i%
\in I_{\ell^{\prime}}}\sum_{k=1}^{K-1}b_{Kj}(d,k,Z_{i})\hat{\gamma }%
_{k\ell,\ell^{\prime}}(d^{\prime},Z_{i}),\text{ }(j=1,...,p).
\end{align*}
We can then compute $\hat{\alpha}_{k\ell}(x)$ as in equation (\ref{nonlin Rr}%
) and $\hat{\theta}$ for $Y_{ki}=1(Q_{i}=k),$ $(k=1,...,K-1)$ and $%
Y_{Ki}=Y_{i}$ as in equation (\ref{nonlin est}).

The orthogonal moment function corresponding to this estimator is%
\begin{align*}
\psi (W,\gamma ,\alpha ,\theta )& =\sum_{k=1}^{K-1}\gamma _{K}(d,k,Z)\gamma
_{k}(d^{\prime },Z)-\theta +\alpha _{K}(D,Q,Z)[Y-\gamma _{K}(D,Q,Z)] \\
+\sum_{k=1}^{K-1}\alpha _{k}(D,Z)[1(Q& =k)-\gamma _{k}(D,Z)],\text{ }\gamma
_{K},\alpha _{K}\in \Gamma _{K},\text{ }\gamma _{k},\alpha _{k}\in \Gamma
_{1},\text{ }(k\leq K-1).
\end{align*}%
where $\Gamma _{K}$ is the set of possible plims of $\hat{\gamma}_{K}$ and $%
\Gamma _{1}$ is the set of plims of $\hat{\gamma}_{k}$ for $k\leq K-1$. This
moment function differs from the multiply robust moment function of Tchetgen
Tchetgen and Shipster (2012) in imposing the constraint that each $\gamma
_{k}$ and $\alpha _{k}$ are contained in the set $\Gamma _{k}$ of possible
plim's of $\hat{\gamma}_{k}.$ For example, when $\hat{\gamma}_{K}$ is a high
dimensional regression estimator $\gamma _{K}$ and $\alpha _{K}$ must be
elements of the mean square span of $(X_{1},X_{2},...)$ similarly to Section
2. It has the multiple robustness feature that for $\bar{\theta}=\mathrm{E}%
[m(W,\bar{\gamma})]$ and any $\alpha =(\alpha _{1},...,\alpha _{K})\in \Pi
_{k=1}^{K}\Gamma _{k},$%
\begin{equation*}
\mathrm{E}[\psi (W,\bar{\gamma},\alpha ,\bar{\theta})]=0,
\end{equation*}%
shown in Chernozhukov et al. (2020) to be a general feature of orthogonal
moment functions constructed from the influence function of $\mathrm{E}%
[m(W,\gamma (F))]$. It also has other multiple robustness features. For $%
\alpha _{k0},$ $(k=1,...,K)$ given in the proof of Corollary 10 in the
Appendix, when $\alpha _{k0}\in \Gamma _{1},$ $(k\leq K-1)$ and $\alpha
_{K0}\in \Gamma _{K},$%
\begin{equation*}
\mathrm{E}[\psi (W,\gamma _{10},...,\gamma _{K-1,0},\gamma _{K},\alpha
_{0},\theta _{0})]=0,\text{ }\mathrm{E}[\psi (W,\gamma _{1},...,\gamma
_{K-1},\gamma _{K0},\alpha _{0},\theta _{0})]=0,\text{ }
\end{equation*}%
for any $\gamma _{K}\in \Gamma _{K}$ and $\gamma _{k}\in \Gamma _{1},$ $%
(k\leq K-1).$

\bigskip

We now return to the general learner $\hat{\theta}$ and give regularity
conditions for asymptotic normality and consistent estimation of the
asymptotic variance of $\hat{\theta}$. For $\tilde{\gamma}=(\tilde{\gamma}%
_{1},...,\tilde{\gamma}_{K})^{\prime}\in\Pi_{k=1}^{K}\Gamma_{k}$ and $%
\gamma_{k}\in\Gamma_{k}$ let 
\begin{equation*}
D_{k}(W,\gamma_{k},\tilde{\gamma}):=\left. \frac{\partial m(W,\tilde{\gamma }%
+e_{k}\tau\gamma_{k})}{\partial\tau}\right\vert _{\tau=0}
\end{equation*}
be the Gateaux derivative of $m(W,\gamma)$ with respect to $\gamma_{k}$ when
it exists. Comparing this definition with equation (\ref{nonlin Rr}) we see
that each $\hat{M}_{kj\ell}$ is an average of values of this Gateaux
derivative. We impose the following condition on these derivatives.

\bigskip

\textsc{Assumption 12: }\textit{There are }$C,$\textit{\ }$\varepsilon >0,$%
\textit{\ }$a_{kj}(w),$\textit{\ and }$A_{k}(w,\gamma)$\textit{\ such that
for all }$\gamma$\textit{\ with }$\Vert\gamma-\bar{\gamma}\Vert\leq
\varepsilon,$\textit{\ }$D_{k}(W,b_{kj},\gamma)$\textit{\ exists and for }$%
k=1,...,K$ 
\begin{align*}
D_{k}(W,b_{kj},\gamma) & =a_{kj}(W)A_{k}(W,\gamma),\text{ }\max_{j\leq
p}\left\vert \mathrm{E}[a_{kj}(W)\{A_{k}(W,\gamma)-A_{k}(W,\bar{\gamma}%
)\}]\right\vert \leq C\left\Vert \gamma-\bar{\gamma}\right\Vert , \\
\max_{j\leq p}\left\vert a_{kj}(W)\right\vert & \leq C,\text{ }\mathrm{E}%
[A_{k}(W,\gamma)^{2}]\leq C.
\end{align*}

\bigskip

This condition and the use of the cross-fit $\hat{\gamma}_{\ell,\ell^{%
\prime}}$ in $\hat{M}_{k\ell}$ lead to a convergence rate for $\hat{M}%
_{k\ell}.$ Let $M_{kj}=\mathrm{E}[D_{k}(W,b_{kj},\bar{\gamma})]\,$\ and $%
M_{k}=(M_{k1},...,M_{kp}),$ $(j=1,...,p;k=1,...,K).$

\bigskip

\textsc{Lemma 8:} \textit{If there is }$0<d_{\gamma }<1/2$ such that $%
\left\Vert \hat{\gamma}_{k\ell ,\ell ^{\prime }}-\bar{\gamma}_{k\ell ,\ell
^{\prime }}\right\Vert =O_{p}(n^{-d_{\gamma }}),$ $(k=1,...,K;\ell ,\ell
^{\prime }=1,...L)$, \textit{and Assumption 12 is satisfied then}%
\begin{equation*}
\left\Vert \hat{M}_{k\ell }-M_{k}\right\Vert _{\infty }=O_{p}(n^{-d_{\gamma
}}).
\end{equation*}

\bigskip

This result can be utilized to obtain mean square convergence rates for $%
\hat{\alpha}_{k}$ from Theorems 1 and 2. As for linear functionals the limit 
$\bar{\alpha}_{k}$ of the estimators $\hat{\alpha}_{k}$ are important for
the properties of $\hat{\theta}$. Here the $\bar{\alpha}_{k}$ are associated
with the Gateaux derivatives $D_{k}(W,\gamma_{k},\bar{\gamma}),$ $%
(k=1,...,K).$ The following condition specifies each $\bar{\alpha}_{k}$ and
specifies the size of the remainder in a linearization using the Gateaux
derivatives.

\bigskip

\textsc{Assumption 13: }\textit{i) For }$(k=1,...,K)$\textit{\ there is }$%
\bar{\alpha}_{k}\in\Gamma_{k}$\textit{\ such that for all }$%
\gamma_{k}\in\Gamma_{k}$\textit{, }$\mathrm{E}[D_{k}(W,\gamma_{k},\bar{\gamma%
})]=\mathrm{E}[\bar{\alpha }_{k}(X_{k})\gamma_{k}(X_{k})];$\textit{\ ii) }$%
\bar{\alpha}_{k}(X_{k})$ \textit{and }$\mathrm{E}[\{Y_{k}-\bar{\gamma}%
_{k}(X_{k})\}^{2}|X_{k}]$\textit{\ are bounded;} \textit{iii) there are }$%
\varepsilon,$\textit{\ }$C>0$\textit{\ such that for all }$%
\gamma\in\Pi_{k=1}^{K}\Gamma_{k}$\textit{\ with }$\Vert \gamma-\bar{\gamma}%
\Vert<\varepsilon,$%
\begin{equation*}
|\mathrm{E}[m(W,\gamma)-m(W,\bar{\gamma})-\sum_{k=1}^{K}D_{k}(W,\gamma_{k}-%
\bar{\gamma },\bar{\gamma})]|\leq C\Vert\gamma-\bar{\gamma}\Vert^{2}.
\end{equation*}

\bigskip

Here each $\bar{\alpha}_{k}$ is specified as the Riesz representer for the
linear functional $\mathrm{E}[D_{k}(W,\gamma_{k},\bar{\gamma})]$ on $%
\gamma_{k}\in\Gamma_{k}$ as in Newey (1994, equation 4.4). Here the
linearization $\mathrm{E}[D_{k}(W,\gamma_{k},\bar{\gamma})]$ has the role
that was fulfilled by the linear functional $\mathrm{E}[m(W,\gamma)]$
earlier. Indeed when $m(W,\gamma)$ is linear then $m(W,\gamma)$ will be its
Gateaux derivative.

From Lemma 8 we see that the convergence rate for each $\hat{M}_{k\ell }$ is
the convergence rate $n^{-d_{\gamma }}$ of $\hat{\gamma}$ rather than $\sqrt{%
\ln (p)/n}.$ Consquently conditions for root-n consistency are different in
the nonlinear $m(W,\gamma )$ case than in the linear one. The following
condition imposes the rate conditions for a nonlinear functional.

\bigskip

\textsc{Assumption 14:} \textit{There is }$1/4<d_{\gamma}<1/2$\textit{\ such
that }$\left\Vert \hat{\gamma}_{k}-\bar{\gamma}_{k}\right\Vert
=O_{p}(n^{-d_{\gamma}}),$\textit{\ }$(k=1,...,K)$\textit{\ and for }$\bar{%
\alpha}=\bar{\alpha}_{k}$\textit{\ and }$b(x)=b_{k}(x_{k})$\textit{, either
i) Assumptions 2 and 3 are satisfied and }$d_{\gamma}(1+4\xi)/(1+2\xi )>1/2$%
\textit{\ or ii) Assumption 7 is satisfied and }$d_{\gamma}>1/3.$\textit{\ }

\bigskip

The requirement $d_{\gamma}>1/4$ given here is familiar for estimators that
depend nonlinearly on unknown functions, e.g. Newey (1994).. Condition i)
allows $d_{\gamma}$ to be any rate greater than $1/4$ if $\xi$ is large
enough. Condition ii), which drops the sparse eigenvalue assumption but
requires absolute summability of the coefficients of each $\bar{\alpha}_{k},$
requires $d_{\gamma}>1/3.$

The following gives the large sample inference results for $\hat{\theta}$
and $\hat{V}.$ Define%
\begin{equation*}
\bar{\theta}=\mathrm{E}[m(W,\bar{\gamma})],\text{ }\psi(w)=m(w,\bar{\gamma})-%
\bar{\theta}+\sum_{k=1}^{K}\bar{\alpha}_{k}(x_{k})[y_{k}-\bar{\gamma}%
_{k}(x)],\text{ }V=\mathrm{E}[\psi(W)^{2}].
\end{equation*}
Here $\bar{\theta}$ will be the object estimated by $\hat{\theta}$ for $\bar{%
\gamma}=$plim$(\hat{\gamma}).$

\bigskip

\textsc{Theorem 9}: \textit{If for }$\Gamma =\Gamma _{k}$\textit{, }$%
b(x)=b_{k}(x_{k})$\textit{, }$r=r_{k}$\textit{\ for }$(k=1,...,K)$\textit{\
and }$\varepsilon _{n}=n^{-d_{\gamma }}$ Assumptions 1, 4, 5, 10, and 12-14
are satisfied \textit{then }$\sqrt{n}(\hat{\theta}-\bar{\theta})\overset{d}{%
\longrightarrow }N(0,V).$ \textit{If in addition Assumption 7 is satisfied
for }$\bar{\alpha}=\bar{\alpha}_{k}$\textit{\ and }$b=b_{k}$ \textit{for
each }$(k=1,...,K)$\textit{\ then }$\hat{V}\overset{p}{\longrightarrow }V$.

\bigskip

\textsc{Example 6: }It is straightforward to specify regularity conditions
for causal mediation that are sufficient for the conditions of Theorem 9 to
hold.

\bigskip

\textsc{Assumption 15:} $\bar{\gamma}_{k}(X_{k})$ \textit{is bounded }$%
(k=1,...,K),$\textit{\ there is }$C>0$\textit{\ such that }$\Pr(D=d,Q=q|Z)>C$%
\textit{\ for all }$d\in\{0,1\},$\textit{\ }$q\in\{1,...,K-1\},$ \textit{and 
}$\mathrm{E}[\{Y-\bar{\gamma}_{K}(D,Q,Z)\}^{2}|D,Q,Z]\leq C.$

\bigskip

This condition is used to guarantee that $\bar{\alpha}_{k}(X_{k})$ is
bounded for each $k.$ For brevity the form of $\bar{\alpha}_{k}(X_{k})$ and $%
\psi(w)$ is given in the Appendix

\bigskip

\textsc{Corollary 10:} \textit{If for }$\Gamma =\Gamma _{k}$\textit{, }$%
b(x)=b_{k}(x_{k})$\textit{, }$r=r_{k}$\textit{\ for }$(k=1,...,K)$\textit{\
and }$\varepsilon _{n}=n^{-d_{\gamma }}$\textit{\ Assumptions 1, 4, 5, 14,
and 15 are satisfied and there is }$C>0$\textit{\ such that }$\left\vert 
\hat{\gamma}_{k}(x_{k})\right\vert \leq C$\textit{\ for all }$x_{k}$ \textit{%
then }$\sqrt{n}(\hat{\theta}-\bar{\theta})\overset{d}{\longrightarrow }%
N(0,V).$\textit{\ If in addition Assumption 7 is satisfied for }$\bar{\alpha}%
=\bar{\alpha}_{k}$\textit{\ and }$b=b_{k}$ \textit{for each }$(k=1,...,K)$%
\textit{\ then }$\hat{V}\overset{p}{\longrightarrow }V$.

\bigskip

The conditions of this result are simple relative to the general regularity
conditions in Assumptions 12 and 13. This simplicity is facilitated by $%
m(W,\gamma )$ being quadratic in $\gamma .$ The condition that $\left\vert 
\hat{\gamma}_{k}(x_{k})\right\vert \leq C$ is not strong for $k=1,...,K-1$
because $Y_{ki}\in \{0,1\}$. For $k=K$ this restriction could be imposed by
truncating $\hat{\gamma}_{k}(x)$ for some $C$ larger than a known bound on $%
\gamma _{K}(X_{k})$ without affecting Assumption 14. In this way Corollary
10 provides a quite simple set of conditions for Auto-DML\ of causal
mediation effects.

\section{Regression Decomposition and the Average Treatment Effect on the
Treated}

In this Section we consider regression decompositions and the average
treatment effect on the treated (ATET). We also give an empirical
application of the ATET using Auto-DML.

\bigskip

\textsc{Example 6:} (Regression Decomposition and ATET): The effect of some
dummy variable $D\in\{0,1\}$ on an outcome variable $Y$ is often of
interest. Regression analysis can be used to decompose the unconditional
effect into an effect conditional on covariates and an effect from a shift
in the covariate distribution when $D$ shifts. One such decomposition takes
the form%
\begin{align*}
\mathrm{E}[Y|D & =1]-\mathrm{E}[Y|D=0]=\Delta_{response}+%
\Delta_{composition}, \\
\Delta_{response} & =\mathrm{E}[Y|D=1]-\frac{\mathrm{E}[D\gamma_{0}(0,Z)]}{%
\Pr(D=1)},\text{ }\Delta_{composition}=\frac{\mathrm{E}[D\gamma_{0}(0,Z)]}{%
\Pr(D=1)}-\mathrm{E}[Y|D=0],
\end{align*}
where $\gamma_{0}(D,Z)=\mathrm{E}[Y|D,Z]$. We will focus here on the
response effect%
\begin{equation*}
\theta_{0}=\Delta_{response}=\frac{\mathrm{E}[D\gamma_{0}(1,Z)]-\mathrm{E}%
[D\gamma_{0}(0,Z)]}{\Pr(D=1)}=\frac{\mathrm{E}[D\{\gamma_{0}(1,Z)-%
\gamma_{0}(0,Z)\}]}{\Pr(D=1)}.
\end{equation*}
This $\theta_{0}$ is the average effect of changing $D$ on the outcome $Y$
conditional on $Z,$ averaged over the subpopulation with $D=1.$ One could
also consider a corresponding effect on the subpopulation with $D=0.$ That
could also be estimated using Auto-DML similarly to $\theta_{0}$ but for
brevity we omit this discussion.

This $\theta_{0}$ is also the ATET when $D$ is a treatment indicator and
potential outcomes are mean independent of treatment conditional on
covariates $Z$. Thus the estimator $\hat{\theta}$ and the asymptotic
variance estimator $\hat{V}$ we give could be applied for inference for the
ATET. We do so in the application given later in this Section.

The key regression functional of interest for $\theta_{0}$ is 
\begin{align}
\mathrm{E}[D\gamma_{0}(0,Z)] & =\mathrm{E}[\pi_{0}(Z)\gamma_{0}(0,Z)]=%
\mathrm{E}[\pi_{0}(Z)\frac {1-D}{1-\pi_{0}(Z)}\gamma_{0}(0,Z)]
\label{Rr ATET} \\
& =\mathrm{E}[\alpha_{0}(X)\gamma_{0}(X)],\text{ }\alpha_{0}(X)=\frac{%
(1-D)\pi_{0}(Z)}{1-\pi_{0}(Z)}.  \notag
\end{align}
Here $\alpha_{0}(X)$ is the Rr of a linear effect as in Section 2 with $%
m(w,\gamma)=d\gamma(0,z).$ The condition $\mathrm{E}[\alpha_{0}(X)^{2}]<%
\infty$ for a finite semiparametric variance bound is $\mathrm{E}%
[1/\{1-\pi_{0}(Z)\}]<\infty.$

\bigskip

The effect $\theta_{0}=\Delta_{response}=ATET$ is a special case of the
nonlinear effect in Section 5 where $\gamma=(\gamma_{1},\gamma_{2})$, $%
Y_{1}=Y,$ $X_{1}=(D,Z),$ $Y_{2}=D,$ $X_{2}=1$, and 
\begin{equation*}
m(w,\gamma)=\frac{y_{2}}{\gamma_{2}}[y_{1}-\gamma_{1}(0,z)].
\end{equation*}
The orthogonal moment function for this object is 
\begin{equation*}
\psi(w,\gamma,\gamma_{2},\alpha,\theta)=\frac{1}{\gamma_{2}}\{d[y-\gamma
(0,z)-\theta]-\alpha(x)[y-\gamma(x)]\},
\end{equation*}
where for notational convenience we let $y_{1}=y,$ $y_{2}=d,$ and $\gamma
_{1}=\gamma$. Similarly to Section 2 this moment function is doubly robust
in that 
\begin{equation*}
\mathrm{E}[\psi(W,\bar{\gamma},\gamma_{20},\bar{\alpha},\theta_{0})]=0
\end{equation*}
if either $\bar{\gamma}(X)=\mathrm{E}[Y|X]$ or $\alpha_{0}(X)\in\Gamma$.

An Auto-DML is given by%
\begin{align}
\hat{\theta} & =\frac{1}{n_{D}}\{\sum_{\ell=1}^{L}\sum_{i\in
I_{\ell}}\{D_{i}[Y_{i}-\hat{\gamma}_{\ell}(0,Z_{i})]-\hat{\alpha}%
_{\ell}(X_{i})[Y_{i}-\hat{\gamma}_{\ell}(X_{i})]\}\}, \\
\hat{V} & =\frac{1}{n}\sum_{\ell=1}^{L}\sum_{i\in I_{\ell}}\hat{\psi}%
_{i\ell}^{2},\text{ }\hat{\psi}_{i\ell}=(\frac{n}{n_{D}})\{D_{i}[Y_{i}-\hat{%
\gamma}_{\ell}(0,Z_{i})-\hat{\theta}]-\hat{\alpha}_{\ell}(X_{i})[Y_{i}-\hat{%
\gamma}_{\ell}(X_{i})]\},  \notag
\end{align}
where $n_{D}$ is the number of treated observations and $\hat{\alpha}_{\ell
}(x)$ is the Lasso learner of the Rr for $m(w,\gamma)=d\gamma(0,z).$
Similarly to the ATE\ in Example 3 we specify the dictionary to be $%
b(x)=[dq(z)^{\prime },(1-d)q(z)^{\prime}]^{\prime},$ where $%
q(z)=(z_{1},...,z_{p/2})^{\prime}$ when $\hat{\gamma}_{\ell}$ is high
dimensional and $q(z)$ is a vector of approximating functions when $\hat{%
\gamma}_{\ell}$ is nonparametric. Then $m(w,b_{j})=d\cdot
b_{j}(0,z)=d\cdot1(j>p/2)q_{j-p/2}(z),$ so that%
\begin{equation*}
\hat{M}_{\ell}=\frac{1}{n-n_{\ell}}\sum_{i\notin I_{\ell}}m(W_{i},b)=\left( 
\begin{array}{c}
0 \\ 
\bar{q}_{\ell}%
\end{array}
\right) ,\text{ }\bar{q}_{\ell}=\frac{1}{n-n_{\ell}}\sum_{i\notin
I_{\ell}}D_{i}q(Z_{i}).
\end{equation*}
Then by block diagonality of $\hat{G}_{\ell}$ and the first block of $\hat {M%
}_{\ell}$ being zero%
\begin{align*}
\hat{\alpha}_{\ell}(x) & =(1-d)q(z)^{\prime}\hat{\rho}_{\ell2},\text{ }\hat{%
\rho}_{\ell2}=\arg\min_{\rho}\{-2\bar{q}_{\ell}^{\prime}\rho_{2}+\rho
_{2}^{\prime}\hat{G}_{\ell2}\rho_{2}+2r\left\Vert \rho_{2}\right\Vert _{1}\},
\\
\hat{G}_{\ell2} & =\frac{1}{n-n_{\ell}}\sum_{i\notin
I_{\ell}}(1-D_{i})q(Z_{i})q(Z_{i})^{\prime}.
\end{align*}
The first order conditions for the Lasso coefficients $\hat{\rho}_{\ell2}$
are%
\begin{equation}
\left\vert \frac{1}{n-n_{\ell}}\sum_{i\notin
I_{\ell}}q_{j}(Z_{i})[D_{i}-(1-D_{i})\hat{\omega}_{\ell i}]\right\vert \leq
r,\text{ }\hat{\omega }_{\ell i}=q(Z_{i})^{\prime}\hat{\rho}_{\ell2},\text{ }%
(j=1,...,p/2).
\end{equation}
The $\hat{\alpha}_{\ell}$ learner sets the \textquotedblleft
weights\textquotedblright\ $\hat{\omega}_{\ell i}$ to approximately
\textquotedblleft balance\textquotedblright\ the treated and untreated
averages for each element of $q(z).$

\bigskip

\textsc{Corollary 11}: \textit{If i) there is }$C>0$ \textit{with }$\pi
_{0}(Z)<1-C$\textit{\ and ii) Assumptions 1, 4, 5, and 9, 11 are satisfied
then for }$\bar{\theta}=\mathrm{E}[D\{Y-\bar{\gamma}(0,Z)\}]/\Pr (D=1)$%
\textit{\ and }$\psi (W)=\Pr (D=1)^{-1}\{D[Y-\bar{\gamma}(0,Z)-\bar{\theta}]-%
\bar{\alpha}(X)[Y-\bar{\gamma}(X)]\},$%
\begin{equation*}
\sqrt{n}(\hat{\theta}-\theta _{0})\overset{d}{\longrightarrow }N(0,V),\text{ 
}V=\mathrm{E}[\psi (W)^{2}].
\end{equation*}%
\textit{If Assumption 7 is also satisfied then} $\hat{V}\overset{p}{%
\longrightarrow }V.$

\bigskip

As an empirical application, we use the Auto-DML of the ATET to estimate the
effect of job training in the National Supported Work Demonstration (NSW), a
job training program for disadvantaged workers that operated in the
mid-1970s. We follow the empirical strategy of LaLonde (1986) and Dehijia
and Wahba (1999), who compare the difference-in-means estimator applied to
an experimental data set with various econometric estimators applied to
\textquotedblleft quasi-experimental\textquotedblright\ data sets. The
experimental data set consists of the treatment and control groups from a
field experiment. A quasi-experimental data set consists of the treatment
group from a field experiment and a comparison group from an unrelated
national survey.

We use sample selection and variable construction as in Dehijia and Wahba
(1999) and Farrell (2015). The outcome $Y$ is earnings in 1978. The
treatment $D$ is an indicator of participation in job training. We consider
three specifications of covariates $Z$. We impose common support of the
propensity score for the treated and untreated groups based on covariates $Z$
as in Farrell (2015). Specifically, we calculate the range of propensity
scores for the treated group, and drop observations in the untreated group
whose propensity scores lie outside this range. We implement this procedure
for each of the three specifications (inducing three different propensity
scores), and ultimately keep the untreated observations that pass all three
tests. In estimation, we consider the fully-interacted dictionary $%
b(D,Z)=(1,D,Z,DZ)$ for all three specifications of $Z$.

The covariate specifications are as follows.

\begin{enumerate}
\item Demographics and earnings, with quadratic terms of continuous
variables. In particular, the covariates are: age, education, black
indicator, Hispanic indicator, married indicator, 1974 earnings, 1975
earnings, age squared, education squared, 1974 earnings squared, and 1975
earnings squared. This specification is moderately flexible. It is one that
an analyst may reasonably implement without knowing the experimental
benchmark ex ante. Here $dim(Z)=11$ and $p=dim(b(D,Z))=24$.

\item Demographics and earnings, with quadratic terms of continuous
variables and constructed indicators. In particular, the covariates are:
those in specification 1; unemployed in 1974 indicator, unemployed in 1975
indicator, and no degree indicator. This specification includes some domain
knowledge about which signals employers may respond to while making hiring
decisions. Note that it does not include conveniently hand-crafted basis
functions to get closer to the experimental benchmark. Here $dim(Z)=14$ and $%
p=dim(b(D,Z))=30$.

\item A high dimensional specification where the covariates are: those in
specification 2; all possible first order interactions, and all polynomials
up to order five of the continuous variables (age, education, 1974 earnings,
1975 earnings). This specification was introduced by Farrell (2015). Here $%
dim(Z)=171$ and $p=dim(b(D,X))=344$.
\end{enumerate}

We estimate the Rr with Lasso minimum distance, and the regression with
Lasso minimum distance, random forests (RF), or neural networks (NN). For
Lasso minimum distance, we use the tuning procedure described in Appendix~%
\ref{sec:computing}. We use the same settings of random forest as
Chernozhukov et al. (2018). We implement a neural network with two hidden
layers of eight units each and linear activation. We use $L=5$ folds in
cross-fitting.

Tables~\ref{tab:ATT_nsw2},~\ref{tab:ATT_psid2}, and~\ref{tab:ATT_cps2}
summarize results for the NSW, PSID, and CPS data sets, respectively. For
comparison, LaLonde (1986) reports $1794$ $(633)$ by difference-in-means
applied to the NSW data, which is the experimental benchmark. Farrell (2015)
reports $1737$ $(869)$ by group Lasso applied to the PSID data using
specification 3. Our corresponding estimate is $1763$ $(1026)$, and our
other results are broadly consistent. To validate the robustness of our
results with respect to the choice of tuning procedure, we report analogous
tables using cross validated regularization in Appendix~\ref%
{sec:additional_empirics}.

\begin{table}[ptb]
\centering
\begin{tabular}{c|c|c|c|c|c|c|c|c}
\hline\hline
spec. & treated & untreated & Lasso ATET & Lasso SE & RF ATET & RF SE & NN
ATET & NN SE \\ \hline
1 & 185 & 172 & 3022.84 & 1278.54 & 3106.55 & 1327.02 & 2585.19 & 1183.85 \\ 
2 & 185 & 172 & 2959.72 & 1253.13 & 3077.26 & 1318.67 & 2606.43 & 1020.56 \\ 
3 & 185 & 172 & 2289.65 & 836.19 & 2785.13 & 819.17 & 2504.83 & 770.24 \\ 
\hline\hline
\end{tabular}
\caption{ATET using NSW treatment and NSW control, by Auto-DML}
\label{tab:ATT_nsw2}
\end{table}

\begin{table}[ptb]
\centering
\begin{tabular}{c|c|c|c|c|c|c|c|c}
\hline\hline
spec. & treated & untreated & Lasso ATET & Lasso SE & RF ATET & RF SE & NN
ATET & NN SE \\ \hline
1 & 185 & 727 & 900.58 & 873.62 & 1521.92 & 977.08 & 197.53 & 946.30 \\ 
2 & 185 & 727 & 1466.35 & 882.67 & 1336.66 & 956.22 & 1447.65 & 980.73 \\ 
3 & 185 & 727 & 1763.20 & 1026.09 & 2010.53 & 987.73 & 2698.55 & 1036.24 \\ 
\hline\hline
\end{tabular}
\caption{ATET using NSW treatment and PSID comparison, by Auto-DML}
\label{tab:ATT_psid2}
\end{table}

\begin{table}[ptb]
\centering
\begin{tabular}{c|c|c|c|c|c|c|c|c}
\hline\hline
spec. & treated & untreated & Lasso ATET & Lasso SE & RF ATET & RF SE & NN
ATET & NN SE \\ \hline
1 & 185 & 5904 & 703.21 & 583.23 & 1639.95 & 616.08 & 1686.77 & 611.81 \\ 
2 & 185 & 5904 & 971.46 & 583.48 & 1584.12 & 616.33 & 1094.86 & 590.31 \\ 
3 & 185 & 5904 & 1358.46 & 614.56 & 1906.62 & 651.77 & 2235.09 & 742.86 \\ 
\hline\hline
\end{tabular}
\caption{ATET using NSW treatment and CPS comparison, by Auto-DML}
\label{tab:ATT_cps2}
\end{table}

\section{Panel Average Derivative and Demand Elasticities}

In this Section, we apply Auto-DML to estimating demand elasticities while
allowing for individual preferences that are correlated with prices and
total expenditure. Specifically, we estimate own-price elasticity in a panel
data model with correlated random slopes. We apply this approach to Nielsen
scanner data.

A panel data model requires double indexing. Let $Y_{it}$, $%
(t=1,...,T_{i},i=1,...,n)$, denote the share of total expenditure on some
good for household $i$ in time period $t$. Let $X_{it}$ be a vector of log
prices, log expenditure, and covariates. Let $\tilde{X}_{i}=(X_{i1}^{%
\prime},...,X_{i,T_{i}}^{\prime})^{\prime}$ collect observations over all
time periods for individual $i$ into one vector. We allow for an unbalanced
panel where different households may have different numbers of observations $%
T_{i}$ as in Wooldridge (2019).

Consider the demand model of Chernozhukov, Hausman, and Newey (2021) given
by 
\begin{equation}
\mathrm{E}[Y_{it}|\tilde{X}_{i},B_{it}]=b_{1}(X_{it})^{\prime }B_{it}.
\label{eq:demand}
\end{equation}%
The $K$-dimensional dictionary $b_{1}(X_{it})$ is a vector of functions of $%
X_{it}$ that includes a constant and, for example, powers of log price and
log expenditure. $B_{it}$ represents household specific preferences that may
vary over time and that may be correlated with regressors from each time
period. We assume the conditional mean of $B_{it}$ is time stationary with%
\begin{equation}
\mathrm{E}[B_{it}|\tilde{X}_{i}]=[I_{K}\otimes \tilde{H}_{i}]^{\prime }\pi
_{0},\quad \pi _{0}=(\pi _{10}^{\prime },...,\pi _{K,0}^{\prime })^{\prime },
\label{eq:stationary}
\end{equation}%
where $I_{K}$ is a $K$-dimensional identity matrix. $\tilde{H}_{i}$ is a
vector of functions of $\tilde{X}_{i}$ with length that does not depend on $%
T_{i}$. This panel model is like that of Chamberlain (1982, 1992),
Chernozhukov et al. (2013b), Graham and Powell (2012), and Wooldridge
(2019), as further discussed in Chernozhukov, Hausman, and Newey (2021).

We will consider identifying and estimating transformations of $\beta _{0}=%
\mathrm{E}[B_{it}]$. $\beta_{0}$ is interpretable as the average marginal
effect of changing $b_{1}(X_{it})$. The transformations we consider will be
interpretable as average income, own-price, and cross-price elasticities. By
law of iterated expectations, our model implies 
\begin{equation}
\beta_{0}=\mathrm{E}[B_{it}]=[I_{K}\otimes \mathrm{E}[\tilde{H}%
_{i}]]^{\prime}\pi_{0}.  \label{eq:LIE}
\end{equation}

Combining~(\ref{eq:demand}),~(\ref{eq:stationary}), and~(\ref{eq:LIE}), we
summarize the correlated random effects model as follows. 
\begin{align}
\gamma_{0}(\tilde{X}_{i}) & =\mathrm{E}[Y_{it}|\tilde{X}_{i}]  \notag \\
& =b_{1}(X_{it})^{\prime}\{\beta_{0}+\mathrm{E}[B_{it}|\tilde{X}_{i}]-\beta
_{0}\}  \notag \\
& =b_{1}(X_{it})^{\prime}\{\beta_{0}+[I_{K}\otimes\tilde{H}%
_{i}]^{\prime}\pi_{0}-[I_{K}\otimes \mathrm{E}[\tilde{H}_{i}]]^{\prime}%
\pi_{0}\}  \notag \\
& =b_{1}(X_{it})^{\prime}\beta_{0}+[b_{1}(X_{it})\otimes(\tilde{H}_{i}-%
\mathrm{E}[\tilde{H}_{i}])]^{\prime}\pi_{0}.  \label{eq:corrRE}
\end{align}
In summary, the choice of $K$-dimensional dictionary $b_{1}(X_{it})$ in the
demand model~(\ref{eq:demand}) induces a $p$-dimensional dictionary $%
b_{it}=b(\tilde{X}_{i})=(b_{1}(X_{it})^{\prime},[b_{1}(X_{it})\otimes (%
\tilde{H}_{i}-\mathrm{E}[\tilde{H}_{i}])]^{\prime})^{\prime}$ in the
correlated random effects model~(\ref{eq:corrRE}). In practice, we replace $%
\mathrm{E}[\tilde{H}_{i}]$ with $\frac{1}{n}\sum_{i=1}^{n} \tilde{H}_{i}$
and set $\tilde{H}_{i}=\frac {1}{T_{i}}\sum_{t=1}^{T_{i}} b_{1}(X_{it})$.

\bigskip

\textsc{Example 10:} Demand elasticities. Denote $X_{it}=(D_{it},Z_{it})$
where $D_{it}$ is log own price. By the derivation in Chernozhukov et al.
(2019) for budget share regressions, an average own-price elasticity is 
\begin{equation*}
\theta_{0}^{\ast}=\frac{\theta_{0}}{\mathrm{E}[Y]}-1,\quad\theta_{0}=\mathrm{%
E} \left[ \frac{\partial\gamma_{0}(\tilde{X}_{i})}{\partial d}\right] .
\end{equation*}
Own-price elasticity $\theta_{0}^{\ast}$ is a smooth transformation of a
linear effect $\theta_{0}$, which in this case is average derivative.
Auto-DML of own-price elasticity is then given by 
\begin{equation*}
\hat{\theta}^{\ast}=\frac{\hat{\theta}}{\frac{1}{n\sum_{i=1}^{n}T_{i}}%
\sum_{i=1}^{n}\sum_{t=1}^{T_{i}}Y_{it}}-1
\end{equation*}
where $\hat{\theta}$ is the Auto-DML of average derivative from Example 4.
Income elasticity and cross-price elasticity have a similar structure; see
Appendix~\ref{sec:delta} for details.

For completeness, we present $\hat{M}_{\ell }$ for average derivative using
the panel data dictionary $b_{it}$. 
\begin{equation*}
\hat{M}_{\ell }=\frac{1}{\sum_{i=1}^{n}T_{i}-\sum_{i\in I_{\ell }}T_{i}}%
\sum_{i\not\in I_{\ell }}\sum_{t=1}^{T_{i}}\frac{\partial b_{it}}{\partial d}%
=\frac{1}{\sum_{i=1}^{n}T_{i}-\sum_{i\in I_{\ell }}T_{i}}\sum_{i\not\in
I_{\ell }}\sum_{t=1}^{T_{i}}\left( 
\begin{array}{c}
\frac{\partial b_{1}(X_{it})}{\partial d} \\ 
0%
\end{array}%
\right) .
\end{equation*}%
Recall Theorem 4 provides consistency and asymptotic normality guarantees
for Auto-DML $\hat{\theta}$. A more sophisticated estimator $\hat{V}$ of the
asymptotic variance of $\hat{\theta}$ is required that accounts for
clustering of observations by household. See the Appendix~\ref{sec:delta}
for details. Importantly, the cluster structure is also preserved in
cross-fitting. Clustering methods for DML were previously used by Chiang et
al. (2019) and Chernozhukov, Hausman, and Newey (2021). The consistency of
the own-price elasticity $\hat{\theta}^{\ast }$ follows from the continuous
mapping theorem, and the asymptotic normality of $\hat{\theta}^{\ast }$
follows from delta method.

As an empirical application, we apply Auto-DML to estimate own-price
elasticity of milk and soda with Nielsen scanner data. The empirical work
here is the researchers' own analyses calculated (or derived) based in part
on data from Nielsen Consumer LLC and marketing databases provided through
the NielsenIQ Datasets at the Kilts Center for Marketing Data Center at The
University of Chicago Booth School of Business. The conclusions drawn from
the NielsenIQ data are those of the researchers and do not reflect the views
of NielsenIQ. NielsenIQ is not responsible for, had no role in, and was not
involved in analyzing and preparing the results reported herein.

The data we use are a subset of the Nielsen Homescan Panel as in Burda,
Harding, and Hausman (2008, 2012). The data include 1483 households from the
Houston-area zip codes for the years 2004-2006. The number of monthly
observations for each household ranges from 12 to 36, with some households
being added and taken away throughout the three years covered. 609
households are included the entire time. Expenditures include all purchases
of the household in each month. The original data had time stamps for
purchases. If a household purchased a good more than once in a month, the
\textquotedblleft monthly price\textquotedblright\ is the average price that
the household paid (i.e. total amount spent on good/total quantity
purchased). We include observations with zero expenditure share as justified
in Chernozhukov, Hausman, and Newey (2021). For those observations, $%
Y_{it}=0 $ and own price is imputed in the ways described in Chernozhukov,
Hausman, and Newey (2021).

We consider 15 groups of goods: bread, butter, cereal, chips, coffee,
cookies, eggs, ice cream, milk, orange juice, salad, soda, soup, water, and
yogurt. As in Burda, Harding, and Hausman (2008, 2012), we choose these
groups because they make up a relatively large proportion of total food
expenditure. We consider budget share regressions for two of these goods:
milk and soda. $Y_{it}$ is share of expenditure spent on milk (soda) by
household $i$ in month $t$. We take as $b_{1}(X_{it})$ the concatenation of
the following variables: fourth order polynomial of log expenditure; fourth
order polynomial of log price for milk (soda); up to fourth order
interactions thereof; and log price of other goods. For $\tilde{H}_{i}$, we
use the time averages of $b_{1}(X_{it})$. Note that $K=dim(b_{1}(X_{it}))=42$
and $p=1521$.

We estimate own-price elasticity according to the procedure outlined
previously in this Section. We estimate both the Rr and the regression with
Lasso minimum distance. For Lasso minimum distance, we use the tuning
procedure described in Appendix~\ref{sec:computing}. We use $L=5$ folds in
cross-fitting. We calculate clustered standard errors by delta method, as
described in Appendix~\ref{sec:delta}.

Table~\ref{tab:dml_elasticity} summarizes results for the milk and soda
own-price elasticities using Auto-DML. For comparison, the cross sectional
estimates for milk and soda elasticities are $-1.27$ ($0.0163$) and $-0.859$
($0.00485$), respectively (Table 1 of Chernozhukov, Hausman, and Newey 2021)
and the corresponding fixed effects estimates are $-.739$ $(.0197)$ and $%
-.853$ ($.00517$). Our results show that allowing for correlated random
coefficients lowers these elasticity estimates, especially the milk
elasticity. These results confirm the finding in Table 5 of Chernozhukov,
Hausman, and Newey (2019), that panel elasticity estimates allowing for
correlation of preferences with prices and total expenditure are much
smaller than cross-section estimates for milk. Our own-price elasticity
estimates are not as small as their slope fixed effect estimates, which for
milk are between $-0.626$ $(0.00849)$ and $-0.496$ $(0.0479)$ and for soda
are between $-0.805$ $(0.00830)$ and $-0.780$ $(0.0235)$ depending on choice
of regularization parameter.

\begin{table}[ptb]
\centering
\begin{tabular}{c|c|c}
\hline\hline
good & elasticity & SE \\ \hline
milk & -.645 & .00649 \\ 
soda & -.826 & .00379 \\ 
\hline\hline
\end{tabular}
\vspace{8pt}
\caption{Average own-price elasticity, by Auto-DML}
\label{tab:dml_elasticity}
\end{table}

For further comparison, we report results from the plug-in approach in Table~%
\ref{tab:dml_elasticity_bias}. The plug-in elasticity estimates are much
closer to the cross-section estimates than the Auto-DML estimates. The
results of this table confirm the importance of debiasing in this
application, with debiased estimates differing from plug-in estimates by
much more than the associated standard errors.

\begin{table}[ptb]
\centering
\begin{tabular}{c|c|c}
\hline\hline
good & elasticity & SE \\ \hline
milk & -.863 & .00255 \\ 
soda & -.863 & .00305 \\ 
\hline\hline
\end{tabular}
\vspace{8pt}
\caption{Average own-price elasticity, by plug-in}
\label{tab:dml_elasticity_bias}
\end{table}

\section{Conclusions}

In this paper we have given an automatic method of debiasing a machine
learner of a parameter of interest that depends on a high dimensional and/or
nonparametric regression. The method only requires the form of the object of
interest. The regression learners are allowed to be anything that converges
in mean square at a fast enough rate. We have shown root-n consistency and
asymptotic normality and given a consistent asymptotic variance estimator
for a wide variety of causal and structural estimators, including nonlinear
functionals of regression. We have applied these methods to estimate the
average treatment effect on the treated in a job training experiment and
have found similar results for Lasso, neural nets, and random forests
regressions. We also have also estimated a correlated random slopes
specification for consumer demand from scanner data and found estimates that
are similar to fixed slope effect elasticities.\newpage

\appendix

\section{Computing Auto-DML}

\label{sec:computing}

\bigskip

\subsection{Tuning}

\label{sec:tuning}

\subsubsection{Theoretical Procedure}

The estimating equation~(\ref{RRLasso}) takes as given the value of
regularization parameter $r_{L}$. For practical use, we provide an iterative
tuning procedure to empirically determine $r_{L}$. Due to its iterative
nature, the tuning procedure is most clearly stated as a replacement for
equation~(\ref{RRLasso}).

Recall that the inputs to equation~(\ref{RRLasso}) are observations in $%
I^{c}_{\ell}$, i.e. excluding fold $\ell$. The analyst must also specify the 
$p$-dimensional dictionary $b$. For notational convenience, we assume $b$
includes the intercept in its first component: $b_{1}(x)=1$. In this tuning
procedure, the analyst must further specify a low-dimensional sub-dictionary 
$b^{\text{low}}$ of $b$. As in equation~(\ref{RRLasso}), the output of the
tuning procedure is $\hat{\rho}_{\ell}$, an estimator of the Rr coefficient
trained only on observations in $I^{c}_{\ell}$.

The tuning procedure is as follows. For observations in $I^{c}_{\ell}$

\begin{enumerate}
\item Initialize $\hat{\rho}_{\ell}$ using $b^{\text{low}}$ 
\begin{align*}
\hat{G}^{\text{low}}_{\ell} & =\frac{1}{n-n_{\ell}}\sum_{i\not \in I_{\ell}}
b^{\text{low}}(X_{i})b^{\text{low}}(X_{i})^{\prime} \\
\hat{M}^{\text{low}}_{\ell} & =\frac{1}{n-n_{\ell}}\sum_{i\not \in I_{\ell}}
m(W_{i},b^{\text{low}}) \\
\hat{\rho}_{\ell} & =%
\begin{bmatrix}
\left( \hat{G}^{\text{low}}_{\ell}\right) ^{-1}\hat{M}^{\text{low}}_{\ell}
\\ 
0%
\end{bmatrix}%
\end{align*}

\item Calculate moments 
\begin{align*}
\hat{G}_{\ell} & =\frac{1}{n-n_{\ell}}\sum_{i\not \in I_{\ell}}
b(X_{i})b(X_{i})^{\prime} \\
\hat{M}_{\ell} & =\frac{1}{n-n_{\ell}}\sum_{i\not \in I_{\ell}} m(W_{i},b)
\end{align*}

\item While $\hat{\rho}_{\ell}$ has not converged

\begin{enumerate}
\item Update normalization 
\begin{equation*}
\hat{D}_{\ell}=\left[ diag\left( \frac{1}{n-n_{\ell}}\sum_{i\not \in
I_{\ell}} [b(X_{i}) b(X_{i})^{\prime}\hat{\rho}_{\ell}- m(W_{i},b)]^{2}
\right) \right] ^{\frac{1}{2}}
\end{equation*}

\item Update $(r_{L},\hat{\rho}_{\ell})$ 
\begin{align*}
r_{L} & =\frac{c_{1}}{\sqrt{n-n_{\ell}}}\Phi^{-1}\left( 1-\frac{c_{2}}{2p}%
\right) \\
\hat{\rho}_{\ell} & =\arg\min_{\rho} \rho^{\prime}\hat{G}_{\ell}\rho
-2\rho^{\prime}\hat{M}_{\ell}+2r_{L} c_{3} |\hat{D}_{\ell,1}\cdot\rho
_{1}|+2r_{L} \sum_{j=2}^{p} |\hat{D}_{\ell,j}\cdot\rho_{j}|
\end{align*}
where $\rho_{j}$ is the $j$-th coordinate of $\rho$ and $\hat{D}_{\ell,j}$
is the $j$-th diagonal entry of $\hat{D}_{\ell}$.
\end{enumerate}
\end{enumerate}

In step 1, $b^{low}$ is sufficiently low-dimensional that $\hat{G}_{\ell }^{%
\text{low}}$ is invertible. In practice, we take $dim(b^{low})=dim(b)/40$.

In step 3, $(c_{1},c_{2},c_{3})$ are hyper-parameters taken as $(1,0.1,0.1)$
in practice. We implement the optimization via generalized coordinate
descent with soft-thresholding. See below for a detailed derivation of this
soft-thresholding routine. We use the same techniques as Chernozhukov,
Newey, and Singh (2018) to improve numerical stability in high dimensional
settings. We use $\hat{D}_{\ell}+0.2I$ instead of $\hat{D}_{\ell}$, and we
cap the maximum number of iterations at 10. We also use warm start: in a
given iteration, the optimization to determine $\hat{\rho}_{\ell}$ is
initialized as the value of $\hat{\rho}_{\ell}$ in the previous iteration.

\subsubsection{Cross-Validation Procedure}

In the theoretical tuning procedure, the hyperparameters $%
(c_{1},c_{2},c_{3}) $ are chosen by the analyst. The hyperparameter $c_{1}$
is of particular importance because it scales $r_{L}$. We now present a
procedure to determine $c_{1}\in\left\{ 5/4,1,3/4,1/2\right\} $ by cross
validation.

In the theoretical tuning procedure, denote by $r_{L}(c_{1})$ the value of
the regularization parameter and denote by $\hat{\rho}_{\ell}(c_{1})$ the
estimated Rr coefficient that are obtained using hyperparameter value $c_{1}$
and observations in $I_{\ell}^{c}$. We define the cross-validated loss for
the hyperparameter $c_{1}$ by 
\begin{equation*}
CV(c_{1})=\sum_{\ell=1}^{L}\sum_{i\in I_{\ell}}[-2m(W_{i},b)^{\prime}\hat {%
\rho}_{\ell}(c_{1})+\{b(X_{i})^{\prime}\hat{\rho}_{\ell}(c_{1})\}^{2}].
\end{equation*}
To determine $c_{1}$ by cross-validation, we solve the optimization problem 
\begin{equation*}
c_{1}^{\ast}=\arg\min_{c_{1}\in\{5/4,1,3/4,1/2\}}CV(c_{1}).
\end{equation*}

\subsubsection{Justification}

The iterative tuning procedure is analogous to Algorithm A.1 and therefore
justified by an argument analogous to Theorem 1 of Belloni et al. (2012).

The analogy is as follows. The normalization $\hat{D}_{\ell}$ is the square
root of the empirical second moment of the dictionary times the regression
residual, just as $\hat{\Gamma}_{\ell}$ in Belloni et al. (2012). The
formula for the regularization parameter is the same, after accounting for
the fact that the objective in the present work uses $r_{L}$ whereas the
objective in eq. 2.4 of Belloni et al. (2012) uses $\frac{\lambda}{n}$.

\subsection{Optimization}

\label{sec:optim}

\subsubsection{Procedure}

The tuning procedure, an elaboration of estimating equation~(\ref{RRLasso}),
involves the minimization of a generalized Lasso objective. We generalize
the coordinate descent approach for Lasso (Fu 1998, Daubechies et al. 2004,
Friedman et al. 2007, Friedman et al. 2010) to the minimum distance Lasso
objective used in the present work. Specifically, we use the following
coordinate-wise soft-thresholding update.

To lighten notation, we abstract from sample splitting, estimation of the
moments and normalization, and special treatment of the intercept. We also
scale the objective by $1/2$: 
\begin{equation*}
\hat{\rho}=\arg\min_{\rho} \frac{1}{2}\rho^{\prime}G\rho-\rho^{%
\prime}M+r_{L}\|D\rho\|_{1}
\end{equation*}
We denote the $j$-th element of a generic vector $V$ by $V_{j}$. We denote
the $(j,k)$-entry of the matrix $G$ by $G_{jk}$.

For $j=1:p$

\begin{enumerate}
\item Calculate loadings that do not depend on $\rho_{j}$ 
\begin{align*}
z_{j} & =G_{jj} \\
\pi_{j} & =M_{j}-\sum_{k\neq j} \rho_{k} G_{jk}
\end{align*}

\item Update coordinate $\rho_{j}$ 
\begin{align*}
\rho_{j} & = \frac{\pi_{j}+D_{j}r_{L}}{z_{j}},\quad\text{if }
\pi_{j}<-D_{j}r_{L} \\
& = 0,\quad\quad\quad\quad\quad\, \text{if } \pi_{j}\in[%
-D_{j}r_{L},D_{j}r_{L}] \\
& = \frac{\pi_{j}-D_{j}r_{L}}{z_{j}},\quad\text{if } \pi_{j}>D_{j}r_{L}
\end{align*}
\end{enumerate}

\subsubsection{Justification}

In this Section, we derive the coordinate-wise soft-thresholding update and
argue that the procedure converges to the minimizer.

Observe that 
\begin{equation*}
\frac{\partial}{\partial\rho_{j}}\left[ \frac{1}{2}\rho^{\prime}G\rho
-\rho^{\prime}M\right] =-\pi_{j}+\rho_{j} z_{j}
\end{equation*}
and the loadings $(z_{j},\pi_{j})$ do not depend on $\rho_{j}$.

The subgradient of the penalty term is 
\begin{align*}
\frac{\partial}{\partial\rho_{j}} r_{L}\|D\rho\|_{1} & =
-D_{j}r_{L}\quad\quad\quad\text{ if } \rho_{j}<0 \\
& = [-D_{j}r_{L},D_{j}r_{L}]\text{ if } \rho_{j}=0 \\
& = D_{j}r_{L} \quad\quad\quad\quad\text{ if }\rho_{j}>0
\end{align*}

In summary, the subgradient of the objective is 
\begin{align*}
\frac{\partial}{\partial\rho_{j}}\left[ \frac{1}{2}\rho^{\prime}G\rho
-\rho^{\prime}M +r_{L}\|D\rho\|_{1}\right] & = -\pi_{j}+\rho_{j}
z_{j}-D_{j}r_{L}\quad\quad\quad\text{ if } \rho_{j}<0 \\
& = [-\pi_{j}-D_{j}r_{L},-\pi_{j}+D_{j}r_{L}]\text{ if } \rho_{j}=0 \\
& = -\pi_{j}+\rho_{j} z_{j}+D_{j}r_{L}\quad\quad\quad\; \text{ if }\rho_{j}>0
\end{align*}
Rearranging yields the component-wise update.

In our minimum distance Lasso procedure, the objective is of the form of eq.
21 of Friedman et al. (2007). 
\begin{align*}
Q(\theta) & =g(\theta)+\sum_{k} h^{k}(\theta^{k}) \\
g(\theta) & =\frac{1}{2}\theta^{\prime}G\theta-M^{\prime}\theta \\
h^{k}(\theta^{k}) & =|\theta^{k}|
\end{align*}
where $g$ is differentiable and convex and $\{h^{k}\}$ are convex. Therefore
coordinate descent converges to the minimizer of the objective (Tseng, 2001).

\subsection{Minimum Distance Lasso Using Simulated Data}

We first validate the minimum distance Lasso estimator for $\hat{\rho}$ on a
design in which the truth is known. We compare our implementation to the
Lasso implementation \verb|LassoShooting.fit| in the \verb|hdm| package at
each point of departure: minimum distance Lasso formulation, theoretical $%
r_{L}$, normalization $\hat{D}$, iteration, and stabilization. Altogether,
this exercise confirms the validity of each technique introduced in the
tuning procedure.

In this design, the ground truth is $\rho_{0}=(1,1,1,0,0,...)$ where $%
dim(\rho_{0})=101$. The data generating process is 
\begin{equation*}
Y=X^{\prime}\rho_{0}+\epsilon
\end{equation*}
where $X=(1,X_{1},...,X_{100})^{\prime}$, $X_{j}\overset{i.i.d.}{\sim }%
\mathcal{N}(0,1)$, and $\epsilon\sim\mathcal{N}(0,1)$. Recall that the
regression coefficient $\rho_{0}$ can be recovered by using the functional $%
m(w,\gamma)=y\gamma(x)$ in the minimum distance Lasso formulation.

\begin{table}[ptb]
\centering
\begin{tabular}{c|c|c}
\hline\hline
algorithm & MSE & $R^{2}$ \\ \hline
Lasso & 0.0060 & 0.17 \\ 
generalized Lasso & 0.0060 & 0.17 \\ 
theoretical $r_{L}$ & 0.0014 & 0.48 \\
normalization $\hat{D}$ & 0.0016 & 0.56 \\ 
iteration: cold start & 0.0014 & 0.50 \\
iteration: warm start & 0.0014 & 0.50 \\
max iteration & 0.0014 & 0.50 \\ 
$\hat{D}+0.2 I$ & 0.0014 & 0.46 \\ 
\hline\hline
\end{tabular}
\vspace{8pt}
\caption{100 simulations}
\label{tab:sim}
\end{table}

In Table~\ref{tab:sim}, we report MSE defined as $|\hat{\rho}%
-\rho_{0}|^{2}_{2}$ of various implementations. Table~\ref{tab:sim} is
cumulative in the sense that each row implements one additional technique
relative to the preceding row. Before using theoretical $r_{L}$, we use $%
r_{L}=0.5$. We use the estimator reported in the final row in the empirical
examples of Sections 6 and 7; it is precisely the estimator defined in the
tuning procedure.

\bigskip

\section{Proofs of Results}

In this Appendix, we give the proofs of the results of the paper, partly
based on useful Lemmas that are stated and proved in this Appendix. We first
give a series of Lemmas like those in Bradic et al. (2021) except that $%
\varepsilon _{n}$ is allowed to be larger than $\sqrt{\ln(p)/n}$ in order to
allow $m(w,\gamma)$ to be nonlinear in $\gamma.$ These Lemmas are used to
prove Theorem 1. Let $\varepsilon_{n}$ be as given in Assumptions 2 and 6
and $s_{0}\geq C\varepsilon_{n}^{-2/(2\xi+1)}$. By Assumption 2 we can
define $J_{0}$ as indices of a sparse approximation with $\left\vert
J_{0}\right\vert =s_{0}$ and coefficients $\tilde{\rho}_{j}$ for $j\in J_{0}$
such that for $\tilde{\alpha}(x)=\sum_{j\in J_{0}}\tilde{\rho}_{j}b_{j}(X),$ 
\begin{equation*}
\mathrm{E}[\{\bar{\alpha}(X)-\tilde{\alpha}(X)\}^{2}]\leq Cs_{0}^{-2\xi}.
\end{equation*}
Define $\rho$ to be the coefficients of a linear projection of $%
\alpha_{0}(X) $ on $b(X)$ so that $\breve{\alpha}(X)=b(X)^{\prime}\rho$
satisfies%
\begin{equation*}
\mathrm{E}[b(X)\{\bar{\alpha}(X)-\breve{\alpha}(X)\}]=0.
\end{equation*}
Also define $\rho_{\ast}$ as 
\begin{equation}
\rho_{\ast}\in\arg\min_{v}\;(\rho-v)\prime
G(\rho-v)+2\varepsilon_{n}\sum_{j\in J_{0}^{c}}|v_{j}|\text{.}
\label{pistar}
\end{equation}

\bigskip

\textsc{Lemma A1:} $\Vert G(\rho_{\ast}-\rho)\Vert_{\infty}\leq\varepsilon
_{n}.$

\bigskip

Proof: Let $e_{j}\in\mathbb{R}^{p}$ denote the $j$-th column of $I_{p}$. The
first-order condition for $\rho^{\ast}$ imply that for $j\in J_{0}$, we have 
$e_{j}\prime G(\rho_{\ast}-\rho)=0$; for $j\in J_{0}^{c}$, we have that $%
e_{j}\prime G(\rho_{\ast}-\rho)+\varepsilon_{n}z_{j}=0$, where $z_{j}=\text{%
sign}(\rho_{\ast,j})$ if $\rho_{\ast,j}\neq0$ and $z_{j}\in \lbrack-1,1]$ if 
$\rho_{\ast,j}=0$. Therefore, for any $j$, we have that $|e_{j}\prime
G(\rho_{\ast}-\rho)|\leq\varepsilon_{n}$. Hence, $\Vert
G(\rho_{\ast}-\rho)\Vert_{\infty}\leq\varepsilon_{n}$. $\square.$

\bigskip

\textsc{Lemma A2:} $(\rho-\rho_{\ast})^{\prime}G(\rho-\rho_{\ast})\leq
C\varepsilon_{n}^{4\xi/(2\xi+1)}$ \textit{and }$\left\Vert \bar{\alpha }%
-b^{\prime}\rho^{\ast}\right\Vert =O(\varepsilon_{n}^{2\xi/(2\xi+1)}).$

\bigskip

Proof: Define $\tilde{\rho}=(\tilde{\rho}_{1},\ldots,\tilde{\rho}_{p})\prime$
as%
\begin{equation*}
\tilde{\rho}_{j}=%
\begin{cases}
\tilde{\rho}_{j} & \text{if}\;j\in J_{0} \\ 
0 & \text{otherwise}\text{.}%
\end{cases}%
\end{equation*}
By the definition of $\rho_{\ast}$, we have that 
\begin{equation}
(\rho-\rho_{\ast})\prime G(\rho-\rho_{\ast})+2\varepsilon_{n}\sum_{j\in
J_{0}^{c}}|\rho_{\ast,j}|\leq(\rho-\tilde{\rho})\prime G(\rho-\tilde{\rho }%
)+2\varepsilon_{n}\sum_{j\in J_{0}^{c}}|\tilde{\rho}_{j}|=(\rho-\tilde{\rho }%
)\prime G(\rho-\tilde{\rho})\text{.}  \label{starobj}
\end{equation}
Note that $\tilde{\alpha}(x)=b(x)^{\prime}\tilde{\rho}$, so by the defintion
of $\breve{\alpha}(x)=b(x)^{\prime}\rho$ we have%
\begin{align*}
(\rho-\tilde{\rho})\prime G(\rho-\tilde{\rho}) & =\mathrm{E}%
[\{b(X)^{\prime}(\rho-\tilde{\rho})\}^{2}]=\Vert\breve{\alpha}-\tilde{\alpha}%
\Vert_{2}^{2}=\left\Vert \bar{\alpha}-\tilde{\alpha}-(\bar{\alpha}-\breve{%
\alpha })\right\Vert _{2}^{2}\leq2(\left\Vert \bar{\alpha}-\tilde{\alpha}%
\right\Vert _{2}^{2}+\left\Vert \bar{\alpha}-\breve{\alpha}\right\Vert
_{2}^{2}) \\
& \leq4\left\Vert \bar{\alpha}-\breve{\alpha}\right\Vert ^{2}\leq
4Cs_{0}^{-2\xi}\leq C\varepsilon_{n}^{4\xi/(2\xi+1)}.
\end{align*}
where the last inequality follows by by $s_{0}\geq
C\varepsilon_{n}^{-2/(2\xi+1)}.$ The result then follows eq. (\ref{starobj})
by and $\varepsilon_{n}\sum_{j\in J_{0}^{c}}|\rho_{\ast,j}|\geq0.$ $\square.$

\bigskip

Define $J$ to be the vector of indices of nonzero elements of $\rho_{\ast}$
and $\left\vert A\right\vert $ be be the number non zero elements of any
finite set $A.$

\bigskip

\textsc{Lemma A3:} $\left\vert J\right\vert \leq
C\varepsilon_{n}^{-2/(2\xi+1)}$.

\bigskip

Proof: For all $j\in J\backslash J_{0}$ the first order conditions to
equation (\ref{pistar}) imply $|e_{j}^{\prime}G(\rho_{\ast}-\rho)|=%
\varepsilon_{n}$. Therefore, It follows that 
\begin{equation*}
\sum_{j\in J\backslash J_{0}}\left( e_{j}^{\prime}G(\rho_{\ast}-\rho)\right)
^{2}=\varepsilon_{n}^{2}|J\backslash J_{0}|\text{.}
\end{equation*}
In addition, 
\begin{align*}
\sum_{j\in J\backslash J_{0}}\left( e_{j}^{\prime}G(\rho_{\ast}-\rho)\right)
^{2} & \leq\sum_{j=1}^{p}\left( e_{j}^{\prime}G(\rho_{\ast}-\rho)\right)
^{2}=(\rho_{\ast}-\rho)^{\prime}G\left( \sum_{j=1}^{p}e_{j}e_{j}^{\prime
}\right) G(\rho_{\ast}-\rho) \\
& =(\rho_{\ast}-\rho)^{\prime}G^{2}(\rho_{\ast}-\rho)\leq\lambda_{\max
}(G)\{(\rho-\rho_{\ast})\prime G(\rho-\rho_{\ast})\}\leq C\varepsilon
_{n}^{4\xi/(2\xi+1)}\text{,}
\end{align*}
where the last inequality follows by Lemma A2 and $\lambda_{\max}(G)\leq C.$
Combining the above two displays, we obtain 
\begin{equation*}
\varepsilon_{n}^{2}|J\backslash J_{0}|\leq C\varepsilon_{n}^{4\xi/(2\xi +1)}%
\text{.}
\end{equation*}
Dividing through by $\varepsilon_{n}^{2}$ gives $|J\backslash J_{0}|\leq
C\varepsilon_{n}^{-2/(2\xi+1)}$. Thus by $s_{0}\leq
C\varepsilon_{n}^{-2/(2\xi+1)},$%
\begin{equation*}
|J|=|J_{0}|+|J\backslash J_{0}|=s_{0}+|J\backslash J_{0}|\leq
s_{0}+C\varepsilon_{n}^{-2/(2\xi+1)}\leq C\varepsilon_{n}^{-2/(2\xi+1)}.%
\text{ }\square.
\end{equation*}

\bigskip

\textsc{Lemma A4:} $\Vert\hat{G}\rho_{\ast}-G\rho_{\ast}\Vert_{\infty}=O_{p}(%
\sqrt{\ln(p)/n}).$

\bigskip

Proof: By $(\rho-\rho_{\ast})^{\prime}G(\rho-\rho_{\ast})\longrightarrow0$
and $\rho^{\prime}G\rho\leq \mathrm{E}[\bar{\alpha}(X)^{2}]$ it follows that 
$\mathrm{E}[\left\{ b(X)^{\prime}\rho_{\ast}\right\}
^{2}]=\rho_{\ast}^{\prime}G\rho_{\ast}\leq C.$ The conclusion then follows
by Assumption 4 and Lemma B2 of Bradic et al. (2021). $\square.$

\bigskip

\textsc{Lemma A5:} \textit{For }$\Delta=\hat{\rho}-\rho^{\ast}$\textit{\ and
any }$\hat{J}$\textit{\ such that }$(\rho^{\ast})_{\hat{J}^{c}}=0,$\textit{\
with probability one then with probability approaching one,}%
\begin{equation*}
\Delta^{\prime}\hat{G}\Delta\leq3r\Vert\Delta\Vert_{1}\text{, }\Vert \Delta_{%
\hat{J}^{c}}\Vert_{1}\leq3\Vert\Delta_{\hat{J}}\Vert_{1}.
\end{equation*}

\bigskip

Proof: By the definition of the estimator $\hat{\rho}$, we have 
\begin{equation*}
\hat{\rho}\prime\hat{G}\hat{\rho}-2\hat{M}^{\prime}\hat{\rho}+2r\Vert\hat {%
\rho}\Vert_{1}\leq\rho_{\ast}^{\prime}\hat{G}\rho_{\ast}-2\hat{M}^{\prime
}\rho_{\ast}+2r\Vert\rho_{\ast}\Vert_{1}\text{.}
\end{equation*}
Plugging $\hat{\rho}=\rho_{\ast}+\Delta$ into the above equation and
rearranging the terms gives%
\begin{equation}
\Delta\prime\hat{G}\Delta+2r\Vert\rho_{\ast}+\Delta\Vert_{1}\leq2r\Vert
\rho_{\ast}\Vert_{1}+2(\hat{M}-\hat{G}\rho_{\ast})\prime\Delta.
\label{key ineq}
\end{equation}
By the definition of $\rho$ and $M=\mathrm{E}[b(X)\bar{\alpha}(X)]$ we have $%
G\rho-M=0.$ Then by Assumption 6, Lemma A1, Lemma A4, and the triangle
inequality 
\begin{align*}
\Vert\hat{G}\rho_{\ast}-\hat{M}\Vert_{\infty} & \leq\Vert\hat{G}\rho_{\ast
}-G\rho_{\ast}\Vert_{\infty}+\Vert M-\hat{M}\Vert_{\infty}+\Vert G\rho_{\ast
}-M\Vert_{\infty} \\
& \leq O_{p}(\varepsilon_{n})+\Vert G\rho-M\Vert_{\infty}+\Vert G(\rho_{\ast
}-\rho)\Vert_{\infty}=O_{p}(\varepsilon_{n})\text{.}
\end{align*}
Therefore, by the Holder inequality we have $\left\vert (\hat{M}-\hat{G}%
\rho_{\ast})\prime\Delta\right\vert \leq\Vert\hat{M}-\hat{G}%
\rho_{\ast}\Vert_{\infty}\Vert\Delta\Vert_{1}=O_{p}(\varepsilon_{n})\Vert%
\Delta\Vert _{1},$ so that by $\varepsilon_{n}=o(r),$%
\begin{equation}
\Delta^{\prime}\hat{G}\Delta+2r\Vert\rho_{\ast}+\Delta\Vert_{1}\leq2r\Vert
\rho_{\ast}\Vert_{1}+O_{p}(\varepsilon_{n})\Vert\Delta\Vert_{1}\leq2r\Vert
\rho_{\ast}\Vert_{1}+r\Vert\Delta\Vert_{1},  \notag  \label{eq: thm bnd pi 2}
\end{equation}
with probability approaching one. Then the triangle inequality $\Vert
\rho_{\ast}\Vert_{1}=\Vert\rho_{\ast}+\Delta-\Delta\Vert_{1}\leq\Vert
\rho_{\ast}+\Delta\Vert_{1}+\Vert\Delta\Vert_{1}$ and subtracting $2r\Vert
\rho_{\ast}+\Delta\Vert_{1}$ from both sides gives the first conclusion.

Next, since $\Delta^{\prime}\hat{G}\Delta\geq0$ it also follows from
equation (\ref{key ineq}) that $2r\Vert\rho_{\ast}+\Delta\Vert_{1}\leq2r%
\Vert\rho _{\ast}\Vert_{1}+r\Vert\Delta\Vert_{1}$, so dividing through by $r$
gives%
\begin{equation*}
2\Vert\rho_{\ast}+\Delta\Vert_{1}\leq2\Vert\rho_{\ast}\Vert_{1}+\Vert
\Delta\Vert_{1}\text{.}
\end{equation*}
It follows by $(\rho_{\ast})_{\hat{J}^{c}}=0$ that $\Vert\rho_{\ast}+\Delta%
\Vert_{1}=\Vert(\rho_{\ast})_{\hat{J}}+\Delta_{\hat{J}}\Vert_{1}+\Vert%
\Delta_{\hat{J}^{c}}\Vert_{1}$ and $\Vert\rho_{\ast}\Vert_{1}=\Vert
(\rho_{\ast})_{\hat{J}}\Vert_{1}$. Substituting in the previous display then
gives%
\begin{align*}
2\Vert(\rho_{\ast})_{\hat{J}}+\Delta_{\hat{J}}\Vert_{1}+2\Vert\Delta_{\hat {J%
}^{c}}\Vert_{1} & \leq2\Vert(\rho_{\ast})_{\hat{J}}\Vert_{1}+\Vert
\Delta\Vert_{1}=2\Vert(\rho_{\ast})_{\hat{J}}\Vert_{1}+\Vert\Delta_{\hat{J}%
}\Vert_{1}+\Vert\Delta_{\hat{J}^{c}}\Vert_{1} \\
& \leq2\left( \Vert(\rho_{\ast})_{\hat{J}}+\Delta_{\hat{J}%
}\Vert_{1}+\Vert\Delta_{\hat{J}}\Vert_{1}\right) +\Vert\Delta_{\hat{J}%
}\Vert_{1}+\Vert\Delta_{\hat{J}^{c}}\Vert_{1} \\
& =2\Vert(\rho_{\ast})_{\hat{J}}+\Delta_{\hat{J}}\Vert_{1}+3\Vert\Delta _{%
\hat{J}}\Vert_{1}+\Vert\Delta_{\hat{J}^{c}}\Vert_{1}\text{.}
\end{align*}
Subtracting $2\Vert(\rho_{\ast})_{\hat{J}}+\Delta_{\hat{J}%
}\Vert_{1}+\Vert\Delta_{\hat{J}^{c}}\Vert_{1}$ from both sides gives the
second conclusion. $\square.$

\bigskip

\textsc{Lemma A6:} $\left\Vert \Delta\right\Vert _{2}=O_{p}((r/\varepsilon
_{n})\varepsilon_{n}^{2\xi/(2\xi+1)}).$

\bigskip

Proof: For $\hat{J}=J$ it follows from Assumption 3 and Lemma A5 that with
probability approaching one,%
\begin{align*}
\left\Vert \Delta_{J}\right\Vert ^{2} & \leq C\Delta^{\prime}\hat{G}%
\Delta\leq Cr\left\Vert \Delta\right\Vert _{1}=Cr(\left\Vert \Delta
_{J}\right\Vert _{1}+\left\Vert \Delta_{J}^{c}\right\Vert _{1})\leq
Cr\left\Vert \Delta_{J}\right\Vert _{1} \\
& \leq Cr\sqrt{\left\vert J\right\vert }\left\Vert \Delta_{J}\right\Vert
_{2}\leq Cr\varepsilon_{n}^{-1/(2\xi+1)}\left\Vert \Delta_{J}\right\Vert
_{2}=C((r/\varepsilon_{n})\varepsilon_{n}^{2\xi/(2\xi+1)}\left\Vert \Delta
_{J}\right\Vert _{2}.
\end{align*}
Dividing through by $\left\Vert \Delta_{J}\right\Vert _{2}$ then gives with
probability approaching one,%
\begin{equation*}
\left\Vert \Delta_{J}\right\Vert _{2}\leq C(r/\varepsilon_{n})\varepsilon
_{n}^{2\xi/(2\xi+1)}
\end{equation*}
Let $N$ denote the indices corresponding to the largest $\left\vert
J\right\vert $ entries in $\Delta_{J^{c}}$, so that $N\subset J^{c}$, $%
|N|=\left\vert J\right\vert $ and $|\Delta_{j}|\geq|\Delta_{k}|$ for any $%
j\in J^{c}\cap N$ and $k\in J^{c}\backslash N$. By Lemma A5 for $\hat {J}%
=J\cup N$ it follows exactly as in second previous display that 
\begin{equation*}
\left\Vert \Delta_{\hat{J}}\right\Vert _{2}\leq
C(r/\varepsilon_{n})\varepsilon_{n}^{2\xi/(2\xi+1)}.
\end{equation*}
By Lemma 6.9 of van de Geer and Buhlmann (2011) and Lemma A5, 
\begin{equation*}
\Vert\Delta_{\hat{J}^{c}}\Vert_{2}\leq(\left\vert J\right\vert
)^{-1/2}\Vert\Delta_{\hat{J}^{c}}\Vert_{1}\leq(\left\vert J\right\vert
)^{-1/2}3\Vert\Delta_{\hat{J}}\Vert_{1}\leq3(\left\vert J\right\vert )^{-1/2}%
\sqrt{\left\vert J\right\vert }\Vert\Delta_{J}\Vert_{2}\leq C(r/\varepsilon
_{n})\varepsilon_{n}^{2\xi/(2\xi+1)}.
\end{equation*}
Therefore, by the triangle inequality with probability approaching one, 
\begin{equation*}
\Vert\Delta\Vert_{2}\leq\Vert\Delta_{\hat{J}}\Vert_{2}+\Vert\Delta_{\hat {J}%
^{c}}\Vert_{2}\leq C(r/\varepsilon_{n})\varepsilon_{n}^{2\xi/(2\xi +1)}.%
\text{ }\square.
\end{equation*}

\bigskip

\textsc{Proof of Theorem 1:} By Lemma A6,%
\begin{equation*}
\Delta^{\prime}G\Delta\leq\lambda_{\max}(G)\left\Vert \Delta\right\Vert
_{2}^{2}=O_{p}((r/\varepsilon_{n})^{2}\varepsilon_{n}^{4\xi/(2\xi+1)}).
\end{equation*}
Then by Lemma A2, the triangle inequality, and Assumption 5, for any $c>0,$%
\begin{align*}
\left\Vert \bar{\alpha}-\hat{\alpha}\right\Vert ^{2} & \leq2\left\Vert \bar{%
\alpha}-b^{\prime}\rho^{\ast}\right\Vert ^{2}+2\left\Vert
b^{\prime}(\rho^{\ast}-\hat{\rho})\right\Vert
^{2}=O(\varepsilon_{n}^{4\xi/(2\xi +1)})+\Delta^{\prime}G\Delta \\
&
=O_{p}((r/\varepsilon_{n})^{2}\varepsilon_{n}^{4\xi/(2\xi+1)})=o_{p}(n^{2c}%
\varepsilon_{n}^{4\xi/(2\xi+1)}).
\end{align*}
Taking square roots of both sides gives the conclusion. $\square.$

\bigskip

Next we give a series of Lemmas that are used to prove Theorem 2.

\bigskip

\textsc{Lemma A7:}\ \textit{If Assumption 7 is satisfied then Assumption 2
is satisfied with }$\xi=1/2.$

\textit{\bigskip}

Proof: Let $J_{s}$ denote the indices of the $s$ largest coefficients in
absolute value and $j_{s}\in J_{s}$ be such that $\left\vert
\rho_{0j_{s}}\right\vert \leq\left\vert \rho_{0j}\right\vert $ for all $j\in
J_{s}.$ Then 
\begin{equation}
s\left\vert \rho_{0j_{s}}\right\vert \leq\sum_{j\in J_{s}}\left\vert \rho
_{j0}\right\vert \leq\sum_{j=1}^{\infty}\left\vert \rho_{j0}\right\vert =C.
\label{L1bound}
\end{equation}
By Assumption 7 $J_{s}\subset\{1,...,p\}.$ Define%
\begin{equation*}
\alpha_{p}(X):=\sum_{j=1}^{p}\rho_{0j}b_{j}(X),\text{ }\alpha_{s}(X):=%
\sum_{j\in J_{s}}\rho_{0j}b_{j}(X).
\end{equation*}
Let $\rho^{p}=(\rho_{01},...,\rho_{0p})$ and $\rho^{s}$ be the vector with $%
\rho_{j}^{s}=\rho_{0j}$ if $j\in J_{s}\,$and $\rho_{j}^{s}=0$ otherwise.
Then by $\left\vert \rho_{0j}\right\vert \leq\left\vert
\rho_{0j_{s}}\right\vert $ for all $j\notin J_{s},$%
\begin{align*}
\left\Vert \alpha_{p}-\alpha_{s}\right\Vert ^{2} & =(\rho^{p}-\rho
^{s})^{\prime}G(\rho^{p}-\rho^{s})\leq C\left\Vert \rho^{p}-\rho
^{s}\right\Vert ^{2}=C\sum_{j\notin J_{s}}\rho_{0j}^{2}\leq C\left\vert
\rho_{0j_{s}}\right\vert \sum_{j\notin J_{s}}\left\vert \rho_{0j}\right\vert
\\
& \leq C\left\vert \rho_{0j_{s}}\right\vert \sum_{j=1}^{\infty}\left\vert
\rho_{0j}\right\vert \leq C\left\vert \rho_{0j_{s}}\right\vert \leq C/s.
\end{align*}
It then follows by Assumption 7 and the triangle and Cauch-Scwartz
inequalities that 
\begin{equation*}
\left\Vert \bar{\alpha}-\alpha_{s}\right\Vert ^{2}\leq2\left\Vert \bar{%
\alpha }-\alpha_{p}\right\Vert ^{2}+2\left\Vert
\alpha_{p}-\alpha_{s}\right\Vert ^{2}\leq C/s.\text{ }\square.
\end{equation*}

\bigskip

Define $\rho_{\ast}\in\arg\min_{\rho}\left\{ \left\Vert \bar{\alpha }%
-b^{\prime}\rho\right\Vert ^{2}+2\varepsilon_{n}\left\vert \rho\right\vert
_{1}\right\} $.

\bigskip

\textsc{Lemma A8: }\textit{If Assumption 7 is satisfied then}%
\begin{equation*}
\left\Vert \bar{\alpha}-b^{\prime}\rho_{\ast}\right\Vert ^{2}\leq
C\varepsilon_{n},\text{ }\left\vert \rho_{\ast}\right\vert _{1}\leq C.
\end{equation*}

\bigskip

Proof: Note that by $\xi=1/2$ as in Lemma A7 we have $s=\varepsilon
_{n}^{-2/(2\xi+1)}=\varepsilon_{n}^{-1}$. By Lemma A7 and the definition of $%
\rho_{\ast},$%
\begin{equation*}
\left\Vert \bar{\alpha}-b^{\prime}\rho_{\ast}\right\Vert ^{2}+2\varepsilon
_{n}\left\vert \rho_{\ast}\right\vert _{1}\leq\left\Vert \bar{\alpha }%
-b^{\prime}\rho_{s}\right\Vert ^{2}+2\varepsilon_{n}\left\vert \rho
_{s}\right\vert _{1}\leq C\varepsilon_{n}.
\end{equation*}
The conclusion follows from the terms on the left-hand side both being
positive. $\square.$

\bigskip

\textsc{Lemma A9:} \textit{If }$\left\Vert \hat{G}-G\right\Vert _{\infty
}=O_{p}(\varepsilon_{n})$ and $\varepsilon_{n}=o(r)$\textit{\ then} $%
\left\Vert \hat{\rho}\right\Vert _{1}=O_{p}(1).$

\bigskip

Proof: For $\Delta=\hat{\rho}-\rho_{\ast}$ equation (\ref{key ineq}) can be
written as%
\begin{equation}
\Delta^{\prime}\hat{G}\Delta+2r\Vert\hat{\rho}\Vert_{1}\leq2r\Vert\rho_{\ast
}\Vert_{1}+2(\hat{M}-\hat{G}\rho_{\ast})^{\prime}\Delta.  \label{key ineq 2}
\end{equation}
By Lemma A8 $\left\Vert \bar{\alpha}-b^{\prime}\rho_{\ast}\right\Vert
^{2}\longrightarrow0$ so that $\mathrm{E}[(b(X)^{\prime}\rho^{\ast})^{2}]%
\leq C.$ Then by Assumption 7, Lemma A8, and the Holder inequality it
follows that 
\begin{equation*}
\left\Vert (\hat{G}-G)\rho_{\ast}\right\Vert _{\infty}\leq\left\Vert \hat {G}%
-G\right\Vert _{\infty}\left\Vert \rho_{\ast1}\right\Vert
=O_{p}(\varepsilon_{n})O_{p}(1)=O_{p}(\varepsilon_{n}).
\end{equation*}
Note that the first order conditions for the minimization of 
\begin{align*}
\left\Vert \bar{\alpha}-b^{\prime}\rho\right\Vert ^{2}+2\varepsilon
_{n}\left\vert \rho\right\vert _{1} & =C+\rho^{\prime}G\rho-2\mathrm{E}[\bar{%
\alpha }(X)b(X)]^{\prime}\rho+2\varepsilon_{n}\left\vert \rho\right\vert _{1}
\\
& =C+\rho^{\prime}G\rho-2M^{\prime}\rho+2\varepsilon_{n}\left\vert
\rho\right\vert _{1}
\end{align*}
imply that $\left\Vert G\rho_{\ast}-M\right\Vert _{\infty}=O(\varepsilon
_{n}),$ similarly to Lemma A1. Then by the triangle inequality,%
\begin{equation*}
\left\Vert \hat{M}-\hat{G}\rho_{\ast}\right\Vert _{\infty}\leq\left\Vert 
\hat{M}-M\right\Vert _{\infty}+\left\Vert (\hat{G}-G)\rho_{\ast}\right\Vert
_{\infty}+\left\Vert M-G\rho_{\ast}\right\Vert _{\infty}=O_{p}(\varepsilon
_{n}).
\end{equation*}
Then by the $\Delta\prime\hat{G}\Delta\geq0$, the Holder and triangle
inequalities, and dividing equation (\ref{key ineq 2}) by $2r$ we have%
\begin{equation*}
\Vert\hat{\rho}\Vert_{1}\leq\Vert\rho_{\ast}\Vert_{1}+\left\Vert \hat{M}-%
\hat{G}\rho_{\ast}\right\Vert _{\infty}\left\Vert \Delta\right\Vert
_{1}/r\leq C+O_{p}(\varepsilon_{n}/r)(\Vert\hat{\rho}\Vert_{1}+\Vert\rho
_{\ast}\Vert_{1})=C+o_{p}(1)\Vert\hat{\rho}\Vert_{1}.
\end{equation*}
Then noting that $o_{p}(1)\Vert\hat{\rho}\Vert_{1}\leq(1/2)\Vert\hat{\rho }%
\Vert_{1}$ with probability approaching one we have%
\begin{equation*}
\Vert\hat{\rho}\Vert_{1}\leq C.\text{ }\square.
\end{equation*}

\bigskip

\textsc{Proof of Theorem 2:} It follows by Lemma A9 that $\left\Vert (G-\hat{%
G})\hat{\rho}\right\Vert _{\infty}\leq\left\Vert G-\hat{G}\right\Vert
_{\infty}\left\Vert \hat{\rho}\right\Vert
_{1}=O_{p}(\varepsilon_{n})O_{p}(1)=O_{p}(\varepsilon_{n}).$ Also, the first
order conditions for Lasso imply $\left\Vert -\hat{G}\hat{\rho}+\hat{M}%
\right\Vert _{\infty}\leq r.$ Also $\left\Vert \hat{M}-M\right\Vert
_{\infty}=O_{p}(\varepsilon_{n})$ and $\left\Vert -G\rho_{\ast}+M\right\Vert
_{\infty}\leq\varepsilon_{n}$ by the first order conditions for $\rho_{\ast}$%
. Then by the triangle inequality%
\begin{equation*}
\left\Vert G(\hat{\rho}-\rho_{\ast})\right\Vert _{\infty}\leq\left\Vert (G-%
\hat{G})\hat{\rho}\right\Vert _{\infty}+\left\Vert -\hat{G}\hat{\rho}+\hat{M}%
\right\Vert _{\infty}+\left\Vert \hat{M}-M\right\Vert _{\infty }+\left\Vert
-G\rho^{\ast}+M\right\Vert _{\infty}=O_{p}(r).
\end{equation*}
Then by Lemma A8%
\begin{equation*}
(\hat{\rho}-\rho_{\ast})^{\prime}G(\hat{\rho}-\rho_{\ast})\leq\left\Vert 
\hat{\rho}-\rho_{\ast}\right\Vert _{1}\left\Vert G(\hat{\rho}-\rho_{\ast
})\right\Vert _{\infty}\leq(\left\Vert \hat{\rho}\right\Vert _{1}+\left\Vert
\rho_{\ast}\right\Vert _{1})O_{p}(r)=O_{p}(r).
\end{equation*}
Then we have%
\begin{equation*}
\left\Vert \bar{\alpha}-\hat{\alpha}\right\Vert ^{2}\leq2\left\Vert \bar{%
\alpha}-b^{\prime}\rho_{\ast}\right\Vert ^{2}+2\left\Vert b^{\prime}(\hat{%
\rho}-\rho_{\ast})\right\Vert ^{2}=O_{p}(\varepsilon_{n})+2(\hat{\rho }%
-\rho_{\ast})^{\prime}G(\hat{\rho}-\rho_{\ast})=O_{p}(r)=o_{p}(n^{2c}%
\varepsilon_{n}),
\end{equation*}
for any $c>0$. Taking square roots of both sides of the inequality gives the
conclusion. $\square.$

\bigskip

\textsc{Lemma A10:} \textit{If Assumption 4 is satisfied then} $\left\Vert 
\hat{G}-G\right\Vert _{\infty}=O_{p}(\sqrt{\ln(p)/n}).$

\bigskip

Proof: Define $\varepsilon_{n}^{\ast}=\sqrt{\ln(p)/n}$ and%
\begin{equation*}
T_{ijk}=b_{j}(X_{i})b_{k}(X_{i})-\mathrm{E}[b_{j}(X_{i})b_{k}(X_{i})],\text{ 
}U_{jk}=\dfrac{1}{n}\sum_{i=1}^{n}T_{ijk}.
\end{equation*}
For any constant $C,$ 
\begin{equation*}
\Pr(|\hat{G}-G|_{\infty}\geq
C\varepsilon_{n}^{\ast})\leq\sum_{j,k=1}^{p}\Pr(|U_{jk}|>C\varepsilon_{n}^{%
\ast})\leq p^{2}\max_{j,k}\Pr(|U_{jk}|>C\varepsilon_{n}^{\ast})
\end{equation*}
Note that $\mathrm{E}[T_{ijk}]=0$ and 
\begin{equation*}
|T_{ijk}|\leq|b_{j}(X_{i})|\cdot|b_{k}(X_{i})|+\mathrm{E}[|b_{j}(X_{i})|%
\cdot |b_{k}(X_{i})|]\leq2C_{b}^{2}.
\end{equation*}
Define $K=2C_{b}^{2}/\sqrt{\ln2}\geq\Vert T_{ijk}\Vert_{_{\Psi_{2}}}.$ By
Hoeffding's inequality (Vershynin, 2018) there is a constant $c$ such that 
\begin{align*}
p^{2}\max_{j,k}\Pr(|U_{jk}|>C\varepsilon_{n}^{\ast}) & \leq2p^{2}\exp\left( -%
\dfrac{c(nC\varepsilon_{n}^{\ast})^{2}}{nK^{2}}\right) \\
& =2p^{2}\exp\left( -\dfrac{cC^{2}\ln(p)}{K^{2}}\right) \\
& \leq2\exp\left( \ln(p)[2-\dfrac{cC^{2}}{K^{2}}]\right) \longrightarrow0
\end{align*}
for any $C>K\sqrt{2/c}.$ Thus for large enough $C$, $\Pr(|\hat{G}-G|_{\infty
}\geq C\sqrt{\ln(p)/n})\longrightarrow0$, implying the conclusion. $\square.$

\bigskip

\textsc{Proof of Theorem 3:} The proof proceeds verifying Assumptions 1-3 of
Chernozhukov et al. (2020, LR). Assumption 1 i) of LR is implied by
Assumption 10. Let $\phi(w,\gamma,\alpha)=\alpha(x)[y-\gamma(x)].$ Note that
by Assumption 9, 
\begin{align*}
\int\{\phi(w,\hat{\alpha}_{\ell},\bar{\gamma})-\phi(w,\bar{\alpha},\bar {%
\gamma})\}^{2}F(dw) & =\int\{\hat{\alpha}_{\ell}(x)-\bar{\alpha}(x)\}^{2}[y-%
\bar{\gamma}(x)]^{2}F(dw) \\
& \leq C\left\Vert \hat{\alpha}_{\ell}-\bar{\alpha}\right\Vert ^{2}\overset{p%
}{\longrightarrow}0, \\
\int\{\phi(w,\bar{\alpha},\hat{\gamma}_{\ell})-\phi(w,\bar{\alpha},\bar {%
\gamma})\}^{2}F(dw) & =\int\bar{\alpha}(x)^{2}[\hat{\gamma}_{\ell}(x)-\bar{%
\gamma}(x)]^{2}F(dx) \\
& \leq C\left\Vert \hat{\gamma}_{\ell}-\bar{\gamma}\right\Vert ^{2}\overset{p%
}{\longrightarrow}0,
\end{align*}
giving Assumptions 1 ii) and 1 iii) of LR.

To verify Assumption 2 of LR, note that by Assumption 8 it follows similarly
to Lemma A10 that Assumption 6 is satisfied for 
\begin{equation*}
\varepsilon _{n}=\sqrt{\ln (p)/n}.
\end{equation*}%
Consider first the first case of Assumption 11 where Assumptions 2 and 3 are
satisfied. By Theorem 1, for any $c>0$ we have 
\begin{equation*}
\left\Vert \hat{\alpha}_{\ell }-\bar{\alpha}\right\Vert =o_{p}(n^{c}[\ln
(n)/n]^{\xi /(2\xi +1)}).
\end{equation*}%
Choose $c=[d_{\gamma }+\xi /(2\xi +1)-1/2]/2>0.$ Then by Assumption 11,%
\begin{equation*}
\sqrt{n}\left\Vert \hat{\alpha}_{\ell }-\bar{\alpha}\right\Vert \left\Vert 
\hat{\gamma}_{\ell }-\bar{\gamma}\right\Vert =o_{p}(n^{c}[\ln (n)]^{\xi
/(2\xi +1)}n^{1/2-\xi /(2\xi +1)-d_{\gamma }})=o_{p}(n^{-c}[\ln (n)]^{\xi
/(2\xi +1)})=o_{p}(1).
\end{equation*}%
Consider now the second case of Assumption 11 where Assumption 7 is
satisfied. Then for $c=(1/4+d_{\gamma }-1/2)/2,$ the conclusion of Theorem 2
gives%
\begin{equation*}
\sqrt{n}\left\Vert \hat{\alpha}_{\ell }-\bar{\alpha}\right\Vert \left\Vert 
\hat{\gamma}_{\ell }-\bar{\gamma}\right\Vert =o_{p}(n^{c}[\ln
(n)]^{1/4}n^{-(1/4)-d_{\gamma }+1/2})=o_{p}(n^{-c}[\ln (n)]^{1/4})=o_{p}(1).
\end{equation*}%
Then by the Cauchy-Schwartz and conditional Markov inequalities we have%
\begin{align*}
& \left\vert \frac{1}{\sqrt{n}}\sum_{i\in I_{\ell }}\{\hat{\alpha}_{\ell
}(X_{i})-\bar{\alpha}(X_{i})\}\{\hat{\gamma}_{\ell }(X_{i})-\bar{\gamma}%
(X_{i})\}\right\vert  \\
& \leq \sqrt{n}\sqrt{\sum_{i\in I_{\ell }}\frac{\{\hat{\alpha}_{\ell
}(X_{i})-\bar{\alpha}(X_{i})\}^{2}}{n}}\sqrt{\sum_{i\in I_{\ell }}\frac{\{%
\hat{\gamma}_{\ell }(X_{i})-\bar{\gamma}(X_{i})\}^{2}}{n}} \\
& \leq \sqrt{n}\left\Vert \hat{\alpha}_{\ell }-\bar{\alpha}\right\Vert
\left\Vert \hat{\gamma}_{\ell }-\bar{\gamma}\right\Vert =o_{p}(1),
\end{align*}%
so that Assumption 2 iii) of LR is satisfied.

To verify Assumption 3 of LR, note that by Assumption 1 $\hat{\alpha}_{\ell
}(x)=b(x)^{\prime}\hat{\rho}_{\ell}\in\Gamma,$ so that%
\begin{equation*}
\int\phi(w,\bar{\gamma},\hat{\alpha}_{\ell})F_{0}(dw)=\int\hat{\alpha}_{\ell
}(x)[y-\bar{\gamma}(x)]F(dw)=0
\end{equation*}
and $\mathrm{E}[m(W,\gamma)-\bar{\theta}+\bar{\alpha}(X)\{Y-\gamma(X)\}]$ is
affine in $\gamma,$ giving Assumption 3 of LR. It then follow by Lemma 15 of
LR that%
\begin{equation*}
\sqrt{n}(\hat{\theta}-\bar{\theta})=\frac{1}{\sqrt{n}}\sum_{\ell=1}^{L}%
\sum_{i\in I_{\ell}}\psi(W_{i},\hat{\gamma}_{\ell},\hat{\alpha}_{\ell},\bar{%
\theta})=\frac{1}{\sqrt{n}}\sum_{i=1}^{n}\psi(W_{i},\bar{\gamma},\bar{\alpha}%
,\bar{\theta})+o_{p}(1).
\end{equation*}
The first conclusion then follows by the central limit theorem.

To show the second conclusion, let $\psi _{i}=\psi _{0}(W_{i})$. Then for $%
i\in I_{\ell },$%
\begin{align*}
(\hat{\psi}_{i\ell }-\psi _{i})^{2}& \leq C\left(
\sum_{j=1}^{3}R_{ij}+R\right) ,\text{ }R_{i1}=[m(W_{i},\hat{\gamma}_{\ell
})-m(W_{i},\gamma _{0})]^{2},\text{ }R_{i2}=\hat{\alpha}_{\ell }(X_{i})^{2}\{%
\hat{\gamma}_{\ell }(X_{i})-\bar{\gamma}(X_{i})\}^{2}, \\
R_{i3}& =\{\hat{\alpha}_{\ell }(X_{i})-\bar{\alpha}(X_{i})\}^{2}\{Y_{i}-\bar{%
\gamma}(X_{i})\}^{2},\text{ }R=(\hat{\theta}-\bar{\theta})^{2}.
\end{align*}%
The first conclusion implies $R\overset{p}{\longrightarrow }0$. Let $%
\mathcal{W}_{-\ell }$ denote the observations not in $I_{\ell }.$ By
Assumption 10,%
\begin{equation*}
\mathrm{E}[R_{i1}|\mathcal{W}_{-\ell }]=\int [m(w,\hat{\gamma}_{\ell })-m(w,%
\bar{\gamma})]^{2}F_{W}(dw)=o_{p}(1).
\end{equation*}%
By Assumption 4 and Lemma A9, uniformly in $x$%
\begin{equation*}
\left\vert \hat{\alpha}_{\ell }(x)\right\vert \leq
\sum_{j=1}^{p}|b_{j}(x)|\left\vert \hat{\rho}_{\ell j}\right\vert \leq
C\left\Vert \hat{\rho}_{\ell }\right\Vert _{1}=O_{p}(1).
\end{equation*}%
Then by Assumption 11, 
\begin{equation*}
\mathrm{E}[R_{i2}|\mathcal{W}_{-\ell }]\leq C\left\Vert \hat{\rho}_{\ell
}\right\Vert _{1}^{2}\int \{\hat{\gamma}_{\ell }(x)-\bar{\gamma}%
(x)\}^{2}F_{W}(dw)=C\left\Vert \hat{\rho}_{\ell }\right\Vert
_{1}^{2}\left\Vert \hat{\gamma}_{\ell }-\bar{\gamma}\right\Vert ^{2}\leq
O_{p}(1)o_{p}(1)=o_{p}(1).
\end{equation*}%
Also by Assumption 9 and iterated expectations%
\begin{align*}
\mathrm{E}[R_{i3}|\mathcal{W}_{-\ell }]& \leq \int \{\hat{\alpha}_{\ell }(x)-%
\bar{\alpha}(x)\}^{2}\mathrm{E}[(Y-\bar{\gamma}(x))^{2}|X=x]F_{X}(dx) \\
& \leq C\int \{\hat{\alpha}_{\ell }(x)-\bar{\alpha}(x)\}^{2}F_{X}(dx)=C\left%
\Vert \hat{\alpha}_{\ell }-\bar{\alpha}\right\Vert ^{2}=o_{p}(1).
\end{align*}%
Then by the triangle inequality,%
\begin{equation*}
\mathrm{E}[\frac{1}{n}\sum_{i\in I_{\ell }}\sum_{j=1}^{3}R_{ij}|\mathcal{W}%
_{-\ell }]\leq \mathrm{E}[R_{i1}|\mathcal{W}_{-\ell }]+\mathrm{E}[R_{i3}|%
\mathcal{W}_{-\ell }]+\mathrm{E}[R_{i3}|\mathcal{W}_{-\ell }]=o_{p}(1).
\end{equation*}%
It then follows by the conditional Markov inequality that $\sum_{i\in
I_{\ell }}\sum_{j=1}^{3}R_{ij}/n=o_{p}(1).$ The triangle inequality and
adding up over $\ell $ then gives $(\hat{\psi}_{i\ell }-\psi _{i})^{2}$%
\begin{equation*}
\frac{1}{n}\sum_{\ell =1}^{L}\sum_{i\in I_{\ell }}(\hat{\psi}_{i\ell }-\psi
_{i})^{2}=o_{p}(1).
\end{equation*}%
Note also that by Assumptions 9 and 10,%
\begin{equation*}
\mathrm{E}[\psi _{i}^{2}]\leq C(1+\mathrm{E}[m(W,\bar{\gamma})^{2}]+\mathrm{E%
}[\bar{\alpha}(X)^{2}\{Y-\bar{\gamma}(X)\}^{2}])<\infty .
\end{equation*}%
Then 
\begin{equation*}
\hat{V}=\frac{1}{n}\sum_{\ell =1}^{L}\sum_{i\in I_{\ell }}\hat{\psi}_{i\ell
}^{2}=\frac{1}{n}\sum_{\ell =1}^{L}\sum_{i\in I_{\ell }}(\hat{\psi}_{i\ell
}-\psi _{i}+\psi _{i})^{2}=\frac{1}{n}\sum_{\ell =1}^{L}\sum_{i\in I_{\ell
}}(\hat{\psi}_{i\ell }-\psi _{i})^{2}+2\frac{1}{n}\sum_{\ell
=1}^{L}\sum_{i\in I_{\ell }}(\hat{\psi}_{i\ell }-\psi _{i})\psi _{i}+\frac{1%
}{n}\sum_{i=1}^{n}\psi _{i}^{2};
\end{equation*}%
Furthermore by the Cauchy-Schwartz and Markov inequalities we have 
\begin{equation*}
\left\vert \frac{1}{n}\sum_{i=1}^{n}(\hat{\psi}_{i}-\psi _{i})\psi
_{i}\right\vert \leq \sqrt{\frac{1}{n}\sum_{i=1}^{n}(\hat{\psi}_{i}-\psi
_{i})^{2}}\sqrt{\frac{1}{n}\sum_{i=1}^{n}\psi _{i}^{2}}\overset{p}{%
\longrightarrow }0.
\end{equation*}%
Then $\hat{V}\overset{p}{\longrightarrow }V$ follows by the triangle
inequality and the law of large numbers. $\square $

\bigskip 

\textsc{Proof of Corollary 4:} Define $\hat{\zeta}_{\ell }=[\sum_{i\in
I_{\ell }}m(W_{i},\hat{\alpha}_{\ell })]/\sum_{i\in I_{\ell }}\hat{\alpha}%
_{\ell }(X_{i})^{2}.$ It follows by $m(W,\gamma )$ linear in $\gamma $ that 
\begin{eqnarray}
\tilde{\theta} &=&\frac{1}{n}\sum_{\ell =1}^{L}\sum_{i\in I_{\ell
}}\{m(W_{i},\hat{\gamma}_{\ell })+\frac{\sum_{i\in I_{\ell }}\hat{\alpha}%
_{\ell }(X_{i})[Y_{i}-\hat{\gamma}_{\ell }(X_{i})]}{\sum_{i\in I_{\ell }}%
\hat{\alpha}_{\ell }(X_{i})^{2}}m(W_{i},\hat{\alpha}_{\ell })\}
\label{ATML rem} \\
&=&\frac{1}{n}\sum_{\ell =1}^{L}\sum_{i\in I_{\ell }}\{m(W_{i},\hat{\gamma}%
_{\ell })+\hat{\alpha}_{\ell }(X_{i})[Y_{i}-\hat{\gamma}_{\ell }(X_{i})]\hat{%
\zeta}_{\ell }\}  \notag \\
&=&\hat{\theta}+\sum_{\ell =1}^{L}\{\hat{\zeta}_{\ell }-1\}\frac{1}{n}%
\sum_{i\in I_{\ell }}\{\hat{\alpha}_{\ell }(X_{i})[Y_{i}-\hat{\gamma}_{\ell
}(X_{i})]\}.  \notag
\end{eqnarray}%
It follows from Assumption 11 similarly to the proof of Theorem 3 that there
are $\nu _{\gamma n}\longrightarrow 0$ and $\nu _{\alpha n}\longrightarrow 0$
such that $\left\Vert \hat{\gamma}_{\ell }-\gamma _{0}\right\Vert
=O_{p}(\nu_{\gamma n})$, $\left\Vert \hat{\alpha}_{\ell }-\bar{\alpha}%
\right\Vert =O_{p}(\nu _{\alpha n})$, and $\sqrt{n}\nu_{\gamma n}\nu_{\alpha
n}\longrightarrow 0.$ We also have%
\begin{eqnarray*}
\frac{1}{n}\sum_{i\in I_{\ell }}\{m(W_{i},\hat{\alpha}_{\ell })-\hat{\alpha}%
_{\ell }(X_{i})^{2}\} &=&T_{1}+T_{2},\text{ }T_{1}=\frac{1}{n}\sum_{i\in
I_{\ell }}\{m(W_{i},\hat{\alpha}_{\ell })-\bar{\alpha}(X_{i})\hat{\alpha}%
_{\ell }(X_{i})\}, \\
T_{2} &=&\frac{1}{n}\sum_{i\in I_{\ell }}\hat{\alpha}(X_{i})\{\bar{\alpha}%
(X_{i})-\hat{\alpha}_{\ell }(X_{i})\}.
\end{eqnarray*}%
By $\hat{\alpha}\in \Gamma $ we have $\mathrm{E}[m(W_{i},\hat{\alpha}_{\ell })-\bar{%
\alpha}(X_{i})\hat{\alpha}_{\ell }(X_{i})|\mathcal{W}_{-\ell }]=0.$ Also by $%
\bar{\alpha}(X)$ bounded and $\mathrm{E}[m(W,\gamma )^{2}]\leq C\left\Vert \gamma
\right\Vert^{2} ,$ 
\begin{eqnarray*}
\mathrm{E}[\{m(W_{i},\hat{\alpha}_{\ell })-\bar{\alpha}(X_{i})\hat{\alpha}_{\ell
}(X_{i})\}^{2}|\mathcal{W}_{-\ell }] &\leq &2\mathrm{E}[m(W_{i},\hat{\alpha}_{\ell
})^{2}|\mathcal{W}_{-\ell }]+2\mathrm{E}[\bar{\alpha}(X_{i})^{2}\hat{\alpha}_{\ell
}(X_{i})^{2}|\mathcal{W}_{-\ell }] \\
&\leq &C\left\Vert \hat{\alpha}_{\ell }\right\Vert ^{2}=O_{p}(1).
\end{eqnarray*}%
Then by the triangle and conditional Markov inequalities $T_{1}=O_{p}(1/%
\sqrt{n})=O_{p}(\nu_{\alpha n}).$ Also by the Cauchy-Schwartz and conditional
Markov inequalities, $\left\Vert \hat{\alpha}_{\ell }\right\Vert^{2}
=O_{p}(1),$ and $\left\Vert \hat{\alpha}_{\ell }-\bar{\alpha}\right\Vert
=O_{p}(\nu_{\alpha n})$ we have%
\begin{equation*}
\left\vert T_{2}\right\vert \leq \{\frac{1}{n}\sum_{i\in I_{\ell }}\hat{%
\alpha}(X_{i})^{2}\}^{1/2}\{\frac{1}{n}\sum_{i\in I_{\ell }}[\bar{\alpha}%
(X_{i})-\hat{\alpha}_{\ell }(X_{i})]^{2}\}^{1/2}=O_{p}(\nu_{\alpha n}).
\end{equation*}%
Note also that $\mathrm{E}[\bar{\alpha}(X)^{2}]>0$ by $\bar{\alpha}(X)\neq 0$ and by
similar arguments to those previous it follows that $\sum_{i\in I_{\ell }}%
\hat{\alpha}_{\ell }(X_{i})^{2}/n_{1}\overset{p}{\longrightarrow }\mathrm{E}[\bar{%
\alpha}(X)^{2}]>0.$ Then%
\begin{equation}
\left\vert \hat{\zeta}_{\ell }-1\right\vert =\frac{\left\vert\frac{1}{n}\sum_{i\in
I_{\ell }}\{m(W_{i},\hat{\alpha}_{\ell })-\hat{\alpha}_{\ell }(X_{i})^{2}\}\right\vert}{%
\frac{1}{n}\sum_{i\in I_{\ell }}\hat{\alpha}_{\ell }(X_{i})^{2}}\leq \frac{%
\left\vert T_{1}\right\vert +\left\vert T_{2}\right\vert }{\frac{1}{n}%
\sum_{i\in I_{\ell }}\hat{\alpha}_{\ell }(X_{i})^{2}}=O_{p}(\nu_{\alpha n}).
\label{ATML remA}
\end{equation}

We also have%
\begin{eqnarray*}
\frac{1}{n}\sum_{i\in I_{\ell }}\hat{\alpha}_{\ell }(X_{i})[Y_{i}-\hat{\gamma%
}_{\ell }(X_{i})] &=&T_{1}+T_{2},\text{ }T_{1}=\frac{1}{n}\sum_{i\in I_{\ell
}}\hat{\alpha}_{\ell }(X_{i})[Y_{i}-\gamma _{0}(X_{i})],\text{ } \\
T_{2} &=&\frac{1}{n}\sum_{i\in I_{\ell }}\hat{\alpha}_{\ell }(X_{i})[\gamma
_{0}(X_{i})-\hat{\gamma}_{\ell }(X_{i})].
\end{eqnarray*}%
Similar to previous arguments we have $\left\vert T_{1}\right\vert =O_{p}(1/%
\sqrt{n})=O_{p}(\nu_{\gamma n})$ and $\left\vert T_{2}\right\vert
=O_{p}(\nu_{\gamma n}),$ so by the triangle inequality $\left\vert \frac{1}{n}%
\sum_{i\in I_{\ell }}\hat{\alpha}_{\ell }(X_{i})[Y_{i}-\hat{\gamma}_{\ell
}(X_{i})]\right\vert =O_{p}(\nu_{\gamma n}).$ It now follows by equation (\ref%
{ATML remA}) that%
\begin{equation*}
\{\hat{\zeta}_{\ell }-1\}\frac{1}{n}\sum_{i\in I_{\ell }}\{\hat{\alpha}%
_{\ell }(X_{i})[Y_{i}-\hat{\gamma}_{\ell }(X_{i})]\}=O_{p}(\nu _{\alpha
n})O_{p}(\nu _{\gamma n})=O_{p}(\nu _{\alpha n}\nu _{\gamma n}),
\end{equation*}%
so that by equation (\ref{ATML rem}) and the triangle inequality,%
\begin{equation*}
\sqrt{n}(\tilde{\theta}-\bar{\theta})=\sqrt{n}(\hat{\theta}-\bar{\theta})+%
\sqrt{n}(\tilde{\theta}-\hat{\theta})=\sqrt{n}(\hat{\theta}-\bar{\theta})+%
\sqrt{n}O_{p}(\nu _{\alpha n}\nu _{\gamma n})=\sqrt{n}(\hat{\theta}-\bar{%
\theta})+o_{p}(1).
\end{equation*}%
The conclusion then follows by Theorem 3 and the Slutzky Theorem. $Q.E.D.$

\bigskip 

\textsc{Proof of Corollary 5:} Note that by Assumption 4,%
\begin{equation*}
\left\vert m(W,b_{j})\right\vert \leq \int \left\vert b_{j}(x)\right\vert
[f_{1}(x)+f_{2}(x)]dx\leq C
\end{equation*}%
so that Assumption 8 is satisfied. Also, $\left\vert \alpha
_{0}(X)\right\vert =\left\vert [f_{1}(x)-f_{0}(x)]/f(x)\right\vert \leq C$
by hypothesis, so by the Cauchy-Schwartz inequality, 
\begin{equation*}
\mathrm{E}[m(W,\gamma )^{2}]=\left\vert \mathrm{E}[\gamma (X)\bar{\alpha}%
(X)]\right\vert ^{2}\leq C\left\vert \mathrm{E}[\left\vert \gamma
(X)\right\vert ]\right\vert ^{2}\leq C\mathrm{E}[\gamma (X)^{2}],
\end{equation*}%
implying Assumption 10. The conclusion then follows by Theorem 3. $Q,E.D.$

\bigskip

\textsc{Proof of Corollary 6}: Integration by parts and Assumption 4 give%
\begin{equation*}
\left\vert m(W,b_{j})\right\vert =\left\vert S(U)b(U,Z)dD\right\vert \leq
C\left\vert b(U,Z)\right\vert \leq C,
\end{equation*}%
so Assumption 8 is satisfied. Also%
\begin{equation*}
\mathrm{E}[m(W,\gamma )^{2}]=\mathrm{E}[\{S(U)\gamma (U,Z)\}^{2}]\leq C%
\mathrm{E}[\gamma (U,Z)^{2}]=\mathrm{E}[f(D|Z)^{-1}\omega (D)\gamma
(X)^{2}]\leq C\mathrm{E}[\gamma (X)^{2}],
\end{equation*}%
so Assumption 10 is satisfied. The conclusion then follows by Theorem 3. $%
Q,E.D.$

\bigskip

\textsc{Proof of Corollary 7:} By Assumption 4 and $m(w,\gamma )=\gamma
(1,z)-\gamma (0,z)$ so by the triangle inequality%
\begin{equation*}
\left\vert m(W,b_{j})\right\vert =\left\vert
b_{j}(1,Z)-b_{j}(0,Z)\right\vert \leq C,
\end{equation*}%
and Assumption 8 is satisfied. Also%
\begin{align*}
\mathrm{E}[m(W,\gamma )^{2}]& \leq C\mathrm{E}[\gamma (1,Z)^{2}]+C\mathrm{E}%
[\gamma (0,Z)^{2}]=C\mathrm{E}[\frac{D}{\pi _{0}(Z)}\gamma (1,Z)^{2}]+C%
\mathrm{E}[\frac{1-D}{1-\pi _{0}(Z)}\gamma (0,Z)^{2}] \\
& =C\mathrm{E}[\left\{ \frac{D}{\pi _{0}(Z)}+\frac{1-D}{1-\pi _{0}(Z)}%
\right\} \gamma (X)^{2}]\leq C\mathrm{E}[\gamma (X)^{2}],
\end{align*}%
so Assumption 10 is satisfied. The conclusion then follows by Theorem 3. $%
Q,E.D.$

\bigskip

\textsc{Proof of Lemma 8:} Define%
\begin{equation*}
\bar{M}_{k}(\gamma )=(\bar{M}_{k1}(\gamma ),...,\bar{M}_{kp}(\gamma
))^{\prime },\text{ }\bar{M}_{kj}(\gamma )=\int D_{k}(W,b_{kj},\gamma )F(dW).
\end{equation*}%
For notational convenience we henceforth suppress the $k$ superscript. Let $%
\mathcal{\Gamma }_{\ell ,\ell ^{\prime }}$ be the event that $\Vert \hat{%
\gamma}_{\ell ,\ell ^{\prime }}-\bar{\gamma}\Vert \leq \varepsilon $ and
note that $\Pr (\mathcal{\Gamma }_{\ell ,\ell ^{\prime }})\longrightarrow 1$
for each $\ell $ and $\ell ^{\prime }$. When $\mathcal{\Gamma }_{\ell ,\ell
^{\prime }}$ occurs, 
\begin{equation*}
\int A(W,\hat{\gamma}_{\ell ,\ell ^{\prime }})^{2}F(dW)\leq C,
\end{equation*}%
by Assumption 11. Define\textbf{\ }%
\begin{equation*}
T_{ij}(\gamma )=D(W_{i},b_{j},\gamma )-\bar{M}_{j}(\gamma ),\text{ }(i\in
I_{\ell ^{\prime }}),\text{ }U_{\ell ^{\prime }j}(\gamma )=\dfrac{1}{n_{\ell
^{\prime }}}\sum_{i\in I_{\ell ^{\prime }}}T_{ij}(\gamma ).
\end{equation*}%
Note that for any constant $C^{\prime }$ and the event $\mathcal{A=\{}%
\max_{j}|U_{\ell ^{\prime }j}(\hat{\gamma}_{\ell ,\ell ^{\prime }})|\geq
C^{\prime }\varepsilon _{n}^{\ast }\}$ where $\varepsilon _{n}^{\ast }=\sqrt{%
\ln (p)/n}$ 
\begin{align}
\Pr (\mathcal{A})& =\Pr (\mathcal{A}|\Gamma _{\ell ,\ell ^{\prime }})\Pr
(\Gamma _{\ell ,\ell ^{\prime }})+\Pr (\mathcal{A}|\Gamma _{\ell ,\ell
^{\prime }}^{c})[1-\Pr (\Gamma _{\ell ,\ell ^{\prime }})]  \notag \\
& \leq \Pr (\max_{j}|U_{\ell ^{\prime }j}(\hat{\gamma}_{\ell ,\ell ^{\prime
}})|\geq C^{\prime }\varepsilon _{n}^{\ast }|\Gamma _{\ell ,\ell ^{\prime
}})+1-\Pr (\Gamma _{\ell ,\ell ^{\prime }}).
\end{align}%
By Lemma B2 of Bradic et al. (2021) there is $C^{\prime }$ large enough that
for any $\delta >0$ with probability approaching one,%
\begin{equation*}
\Pr (\max_{j}|U_{\ell ^{\prime }j}(\hat{\gamma}_{\ell ,\ell ^{\prime
}})|\geq C^{\prime }\varepsilon _{n}^{\ast }|\Gamma _{\ell ,\ell ^{\prime
}})<\delta /2.
\end{equation*}%
Also $1-\Pr (\Gamma _{\ell ,\ell ^{\prime }})\longrightarrow 0$, so that $%
\Pr (\mathcal{A})<\delta $ for all $n$ large enough. Therefore 
\begin{equation*}
\left\Vert U_{\ell ^{\prime }}(\hat{\gamma}_{\ell ,\ell ^{\prime
}})\right\Vert _{\infty }=\max_{j}|U_{\ell ^{\prime }j}(\hat{\gamma}_{\ell
,\ell ^{\prime }})|=O_{p}(\varepsilon _{n}^{\ast }).
\end{equation*}

Next, for each $\ell$ it follows that $n-n_{\ell}=\sum_{\ell^{\prime}\neq%
\ell }n_{\ell^{\prime}}$ and 
\begin{equation*}
\left\vert \hat{M}_{\ell}-\sum_{\ell^{\prime}\neq\ell}\frac{n_{\ell^{\prime}}%
}{n-n_{\ell}}\bar{M}(\hat{\gamma}_{\ell,\ell^{\prime}})\right\vert _{\infty
}=\left\vert \sum_{\ell^{\prime}\neq\ell}\frac{n_{\ell^{\prime}}}{n-n_{\ell}}%
U_{\ell^{\prime}}(\hat{\gamma}_{\ell,\ell^{\prime}})\right\vert
_{\infty}\leq\sum_{\ell^{\prime}\neq\ell}\frac{n_{\ell^{\prime}}}{n-n_{\ell}}%
\left\Vert U_{\ell^{\prime}}(\hat{\gamma}_{\ell,\ell^{\prime}})\right\Vert
_{\infty }=O_{p}(\varepsilon_{n}^{\ast}).
\end{equation*}
Also by Assumption 12 and $\Pr(\Gamma_{\ell,\ell^{\prime}})\longrightarrow1$
for each $\ell$ and $\ell^{\prime},$%
\begin{equation*}
\left\vert \sum_{\ell^{\prime}\neq\ell}\frac{n_{\ell^{\prime}}}{n-n_{\ell}}%
\bar{M}(\hat{\gamma}_{\ell,\ell^{\prime}})-M\right\vert _{\infty}=\left\vert
\sum_{\ell^{\prime}\neq\ell}\frac{n_{\ell^{\prime}}}{n-n_{\ell^{\prime}}}[%
\bar{M}(\hat{\gamma}_{\ell,\ell^{\prime}})-M]\right\vert _{\infty}\leq
C\sum_{\ell^{\prime}\neq\ell}\frac{n_{\ell^{\prime}}}{n-n_{\ell^{\prime}}}%
\Vert\hat{\gamma}_{\ell,\ell^{\prime}}-\gamma_{0}\Vert=O_{p}(n^{-d_{%
\gamma}}).
\end{equation*}
The conclusion then follows by the triangle inequality. $\square.$

\bigskip

\textsc{Proof of Theorem 9:} The proof proceeds verifying Assumptions 1-3 of
Chernozhukov et al. (2020, LR) similarly to the proof of Theorem 3. By
Assumption 14, if Assumptions 2 and 3 are satisfied it follows by Lemma 8
that Assumption 6 is satisfied with $\varepsilon _{n}=n^{-d\gamma },$ so by
Theorem 1,%
\begin{equation*}
\left\Vert \hat{\alpha}_{k}-\bar{\alpha}\right\Vert
=o_{p}(n^{c}n^{-d_{\gamma }2\xi /(2\xi +1)}).
\end{equation*}%
Then for $c=[d_{\gamma }(2\xi /(2\xi +1)+d_{\gamma }-1/2]/2=[d_{\gamma
}(4\xi +1)/(2\xi +1)-1/2]/2>0$ we have $\left\Vert \hat{\alpha}_{k}-\bar{%
\alpha}_{k}\right\Vert ^{2}=o_{p}(1)$ and 
\begin{equation*}
\sqrt{n}\left\Vert \hat{\alpha}_{k}-\bar{\alpha}_{k}\right\Vert \left\Vert 
\hat{\gamma}_{k}-\bar{\gamma}_{k}\right\Vert =\sqrt{n}o_{p}(n^{c}n^{-d_{%
\gamma }2\xi /(2\xi +1)})O_{p}(n^{-d\gamma })=o_{p}(n^{c}n^{-2c})=o_{p}(1),
\end{equation*}%
for each $k$. Similarly, by Assumption 14 if Assumption 7 is satisfied
(rather than Assumptions 2 and 3) then by Theorem 2 for any $c>0,$ 
\begin{equation*}
\left\Vert \hat{\alpha}_{k}-\bar{\alpha}_{k}\right\Vert
=o_{p}(n^{c}n^{-d_{\gamma }/2}).
\end{equation*}%
Then for $c=[d_{\gamma }/2+d_{\gamma }-1/2]/2>0=[3d_{\gamma }/2-1/2]/2$ we
have $\left\Vert \hat{\alpha}_{k}-\bar{\alpha}\right\Vert =o_{p}(1)$ and%
\begin{equation*}
\sqrt{n}\left\Vert \hat{\alpha}_{k}-\bar{\alpha}_{k}\right\Vert \left\Vert 
\hat{\gamma}_{k}-\bar{\gamma}_{k}\right\Vert =\sqrt{n}o_{p}(n^{c}n^{-d_{%
\gamma }/2})O_{p}(n^{-d\gamma })=o_{p}(n^{c}n^{-2c})=o_{p}(1),
\end{equation*}%
for each $k.$

Next, Assumption 1 i) of LR is implied by Assumption 10. Let $\phi
_{k}(w,\gamma_{k},\alpha_{k})=\alpha_{k}(x_{k})[y_{k}-\gamma_{k}(x_{k})]$
and 
\begin{equation*}
\phi(w,\gamma,\alpha)=\sum_{k=1}^{K}\phi_{k}(w,\gamma_{k},\alpha_{k})
\end{equation*}
Note that by $\mathrm{E}[\{Y_{k}-\bar{\gamma}(X_{k})\}^{2}|X_{k}]$ and $\bar{%
\alpha }_{k}(X_{k})$ bounded,%
\begin{align*}
\int\{\phi_{k}(w,\hat{\alpha}_{k\ell},\bar{\gamma}_{k})-\phi_{k}(w,\bar {%
\alpha}_{k},\bar{\gamma}_{k})\}^{2}F(dw) & =\int\{\hat{\alpha}_{k\ell
}(x_{k})-\bar{\alpha}_{k}(x_{k})\}^{2}[y_{k}-\bar{\gamma}%
_{k}(x_{k})]^{2}F(dw) \\
& \leq C\left\Vert \hat{\alpha}_{k\ell}-\bar{\alpha}_{k}\right\Vert ^{2}%
\overset{p}{\longrightarrow}0, \\
\int\{\phi_{k}(w,\bar{\alpha}_{k},\hat{\gamma}_{k\ell})-\phi_{k}(w,\bar {%
\alpha}_{k},\bar{\gamma}_{k})\}^{2}F(dw) & =\int\bar{\alpha}_{k}(x_{k})^{2}[%
\hat{\gamma}_{k\ell}(x_{k})-\bar{\gamma}(x_{k})]^{2}F(dx_{k}) \\
& \leq C\left\Vert \hat{\gamma}_{k\ell}-\bar{\gamma}_{k}\right\Vert ^{2}%
\overset{p}{\longrightarrow}0,
\end{align*}
so that Assumptions 1 ii) and 1 iii) of LR are satisfied by the triangle
inequality.

By the Cauchy-Schwartz and conditional Markov inequalities we have%
\begin{align*}
& \left\vert \frac{1}{\sqrt{n}}\sum_{i\in I_{\ell}}\{\hat{\alpha}_{k\ell
}(X_{ki})-\bar{\alpha}_{k}(X_{ki})\}\{\hat{\gamma}_{k\ell}(X_{ki})-\bar {%
\gamma}_{k}(X_{ki})\}\right\vert \\
& \leq\sqrt{n}\sqrt{\sum_{i\in I_{\ell}}\frac{\{\hat{\alpha}_{k\ell}(X_{i})-%
\bar{\alpha}_{k}(X_{i})\}^{2}}{n}}\sqrt{\sum_{i\in I_{\ell}}\frac{\{\hat{%
\gamma}_{k\ell}(X_{ki})-\bar{\gamma}_{k}(X_{ki})\}^{2}}{n}} \\
& =O_{p}(\sqrt{n}\left\Vert \hat{\alpha}_{k\ell}-\bar{\alpha}_{k}\right\Vert
\left\Vert \hat{\gamma}_{k\ell}-\bar{\gamma}_{k}\right\Vert )=o_{p}(1).
\end{align*}
Then by the triangle inequality Assumption 2 of LR is satisfied.

To verify Assumption 3 of LR, note that by Assumption 1 $\hat{\alpha}_{\ell
}(x)=b(x)^{\prime}\hat{\rho}_{\ell}\in\Gamma,$ so that%
\begin{equation*}
\int\phi_{k}(w,\bar{\gamma}_{k},\hat{\alpha}_{k\ell})F_{0}(dw)=\int\hat {%
\alpha}_{k\ell}(x_{k})[y_{k}-\bar{\gamma}_{k}(x_{k})]F(dw)=0.
\end{equation*}
Also note that for each $k,$ 
\begin{align*}
\mathrm{E}[\phi_{k}(W,\gamma_{k},\bar{\alpha}_{k})] & =\mathrm{E}[\bar{\alpha%
}_{k}(X_{k})\{Y_{k}-\gamma_{k}(X_{k})\}]=\mathrm{E}[\bar{\alpha}_{k}(X_{k})\{%
\bar{\gamma }_{k}(X_{k})-\gamma_{k}(X_{k})\}] \\
& =\mathrm{E}[D_{k}(W,\bar{\gamma}_{k},\bar{\gamma})]-\mathrm{E}%
[D_{k}(W,\gamma_{k},\bar {\gamma})]=-\mathrm{E}[D_{k}(W,\gamma_{k}-\bar{%
\gamma}_{k},\bar{\gamma})].
\end{align*}
Then by Assumption 13 for all $\gamma$ with $\left\Vert \gamma-\bar{\gamma }%
\right\Vert <\varepsilon,$%
\begin{align*}
\left\vert \mathrm{E}[\psi(W,\gamma,\bar{\alpha},\bar{\theta})]\right\vert &
=\left\vert \mathrm{E}[m(W,\gamma)-m(W,\bar{\gamma})+\sum_{k=1}^{K}%
\phi_{k}(W,\gamma_{k},\bar{\alpha}_{k})]\right\vert \\
& =\left\vert \mathrm{E}[m(W,\gamma)-m(W,\bar{\gamma})-%
\sum_{k=1}^{K}D_{k}(W,\gamma_{k}-\bar{\gamma}_{k},\bar{\gamma})]\right\vert
\leq C\left\Vert \gamma-\bar{\gamma}\right\Vert ^{2},
\end{align*}
giving Assumption 3 of LR.

It then follows by Lemma 15 of LR that%
\begin{equation*}
\sqrt{n}(\hat{\theta}-\bar{\theta})=\frac{1}{\sqrt{n}}\sum_{\ell=1}^{L}%
\sum_{i\in I_{\ell}}\psi(W_{i},\hat{\gamma}_{\ell},\hat{\alpha}_{\ell},\bar{%
\theta})=\frac{1}{\sqrt{n}}\sum_{i=1}^{n}\psi(W_{i},\bar{\gamma},\bar{\alpha}%
,\bar{\theta})+o_{p}(1).
\end{equation*}
The first conclusion then follows by the central limit theorem$.$

The second conclusion follows by the triangle inequality as in the proof of
Theorem 3 with $R_{i2}$ replaced by $\hat{\alpha}_{k\ell}(X_{ki})^{2}\{\hat{%
\gamma}_{k\ell}(X_{ki})-\bar{\gamma}_{k}(X_{ki})\}^{2}$ and $R_{i3}$ by $\{%
\hat{\alpha}_{k\ell}(X_{ki})-\bar{\alpha}_{k}(X_{ki})\}^{2}\{Y_{ki}-\bar{%
\gamma}_{k}(X_{ki})\}^{2}$ for each $k.$ $\square$

\bigskip

\textsc{Proof of Corollary 10:} The proof proceeds by showing that the
conditions of Theorem 9 are satisfied. By $\bar{\gamma}_{k}$ bounded for
each $k$ and the triangle inequality, $\mathrm{E}[m(W,\bar{\gamma}%
)^{2}]<\infty .$ Also, by the triangle inequality, 
\begin{align*}
\lbrack m(W,\hat{\gamma})-m(W,\bar{\gamma})]^{2}& \leq
C\sum_{k=1}^{K-1}\left\vert \hat{\gamma}_{K}(d,k,Z)\hat{\gamma}%
_{k}(d^{\prime },Z)-\bar{\gamma}_{K}(d,k,Z)\bar{\gamma}_{k}(d^{\prime
},Z)\right\vert ^{2} \\
& \leq C\sum_{k=1}^{K-1}\left\vert \hat{\gamma}_{K}(d,k,Z)-\bar{\gamma}%
_{K}(d,k,Z)|^{2}\hat{\gamma}_{k}(d^{\prime },Z)^{2}-|\bar{\gamma}%
_{K}(d,k,Z)|^{2}|\hat{\gamma}_{k}(d^{\prime },Z)-\bar{\gamma}_{k}(d^{\prime
},Z)\right\vert ^{2} \\
& \leq C\sum_{k=1}^{K-1}(\left\vert \hat{\gamma}_{K}(d,k,Z)-\bar{\gamma}%
_{K}(d,k,Z)|^{2}+|\hat{\gamma}_{k}(d^{\prime },Z)-\bar{\gamma}_{k}(d^{\prime
},Z)\right\vert ^{2}).
\end{align*}%
Therefore we have%
\begin{equation*}
\int [m(w,\hat{\gamma})-m(w,\bar{\gamma})]^{2}F(dw)\leq
\sum_{k=1}^{K-1}(\int |\hat{\gamma}_{K}(d,k,z)-\bar{\gamma}%
_{K}(d,k,z)|^{2}F_{Z}(dz)+\int |\hat{\gamma}_{k}(d^{\prime },z)-\bar{\gamma}%
_{k}(d^{\prime },z)|^{2}F_{Z}(dz)).
\end{equation*}%
By $\pi (d,k|Z)=\Pr (D=d,Q=k|Z)\geq C$ for each $d$ and $k$ we have for any $%
\gamma _{K}(d,k,Z)$ 
\begin{align*}
\int |\gamma _{K}(d,k,z)-\bar{\gamma}_{K}(d,k,z)|^{2}F_{Z}(dz)& =\mathrm{E}%
[|\gamma _{K}(d,k,Z)-\bar{\gamma}_{K}(d,k,Z)|^{2}] \\
& =\mathrm{E}[\frac{1(D=d,Q=k)}{\pi (d,k|Z)}|\gamma _{K}(d,k,Z)-\bar{\gamma}%
_{K}(d,k,Z)|^{2}] \\
& =\mathrm{E}[\frac{1(D=d,Q=k)}{\pi (d,k|Z)}|\gamma _{K}(D,Q,Z)-\bar{\gamma}%
_{K}(D,Q,Z)|^{2}] \\
& \leq C\mathrm{E}[|\gamma _{K}(D,Q,Z)-\bar{\gamma}_{K}(D,Q,Z)|^{2}]=C\left%
\Vert \gamma _{K}-\bar{\gamma}_{K}\right\Vert ^{2}.
\end{align*}%
Applying this calculation to $\gamma _{K}=\hat{\gamma}_{K}$ gives%
\begin{equation*}
\int |\hat{\gamma}_{K}(d,k,z)-\bar{\gamma}_{K}(d,k,z)|^{2}F_{Z}(dz)\leq
C\left\Vert \hat{\gamma}_{K}-\bar{\gamma}_{K}\right\Vert ^{2}.
\end{equation*}%
Also it follows by $\pi (d,k|Z)\geq C$ for each $k$ that $\pi (d|Z)=\Pr
(D=d|Z)\geq C.$ Then similarly to the previous inequality we have$)$ 
\begin{equation*}
|\gamma _{k}(d^{\prime },Z)-\bar{\gamma}_{k}(d^{\prime },Z)|\leq C\left\Vert
\gamma _{k}-\bar{\gamma}_{k}\right\Vert ^{2},\text{ }k=1,...,K-1.
\end{equation*}%
Then collecting terms we have 
\begin{equation*}
\int [m(w,\hat{\gamma})-m(w,\bar{\gamma})]^{2}F(dw)\leq
C\sum_{k=1}^{K}\left\Vert \hat{\gamma}_{k}-\bar{\gamma}_{k}\right\Vert
^{2}\leq C\left\Vert \hat{\gamma}-\bar{\gamma}\right\Vert ^{2},
\end{equation*}%
for $\left\Vert \gamma \right\Vert =\sum_{k=1}^{K}\left\Vert \gamma
_{k}\right\Vert .$ Thus Assumption 10 is satisfied.

Next, by the Gateaux derivative formula in the body of the paper for $%
(k=1,...,K-1)$ we have%
\begin{equation*}
D_{k}(W,b_{kj},\gamma )=a_{kj}(W)A_{k}(W,\gamma ),\text{ }%
a_{kj}(W)=b_{1j}(d^{\prime },Z),\text{ }A_{k}(W,\gamma )=\gamma _{K}(d,k,Z).
\end{equation*}%
It follows similarly to the verification of Assumption 10 and by Assumption
4 that%
\begin{equation*}
\max_{j\leq p}\left\vert a_{kj}(W)\right\vert \leq C,\text{ and }\mathrm{E}%
[A_{k}(W,\gamma )^{2}]\leq C\left\Vert \gamma \right\Vert ^{2},\text{ }%
(k=1,...,K-1).
\end{equation*}%
Also, we have 
\begin{equation*}
D_{K}(W,b_{Kj},\gamma )=\sum_{k=1}^{K-1}b_{Kj}(d,k,Z)\gamma _{k1}(d^{\prime
},Z),
\end{equation*}%
which also has the form like that Assumption 12 where the conclusion of
Lemma 8 will also be satisfied. The second part of Assumption 12 follows by
a similar argument, so that Assumption 12 is satisfied.

Turning now to Assumption 13, note that for $(k=1,...,K-1),$%
\begin{align*}
\mathrm{E}[D_{k}(W,\gamma_{k},\bar{\gamma})] & =\mathrm{E}[\bar{\gamma}%
_{K}(d,k,Z)\gamma _{k}(d^{\prime},Z)]=\mathrm{E}[\bar{\alpha}%
_{k}(X_{k})\gamma_{k}(X_{k})],\text{ }\bar{\alpha}_{k}(X_{k})=\frac{\bar{%
\gamma}_{K}(d,k,Z)1(D=d^{\prime})}{\pi(d^{\prime}|Z)}, \\
\mathrm{E}[D_{K}(W,\gamma_{K},\bar{\gamma})] & =\mathrm{E}%
[\sum_{k=1}^{K-1}\gamma _{K}(d,k,Z)\bar{\gamma}_{k}(d^{\prime},Z)]=\mathrm{E}%
[\bar{\alpha}_{K}(X_{K})\gamma _{K}(X_{K})],\text{ } \\
\bar{\alpha}_{K}(X_{K}) & =\sum_{k=1}^{K-1}\frac{1(D=d,Q=k)\bar{\gamma}%
_{k}(d,Z)}{\pi(d,k|Z)}.
\end{align*}
Each of $\bar{\alpha}_{k}(X_{k})$ is bounded by $\pi(d,k|Z)\geq C$ for $%
d\in\{0,1\}$ and $(k=1,...,K-1)$ and $\bar{\gamma}_{k}(X_{k})$ bounded for
each $k.$ Similarly these conditions imply that $\mathrm{E}[Y_{ki}-\bar{%
\gamma}($

To verify Assumption 13 iii) note that by algebra we have 
\begin{equation*}
m(W,\gamma )-m(W,\bar{\gamma})-\sum_{k=1}^{K}D_{k}(W,\gamma _{k}-\bar{\gamma}%
_{k},\bar{\gamma})=\sum_{k=1}^{K-1}\{\gamma _{K}(d,k,Z)-\bar{\gamma}%
_{K}(d,k,Z)\}\{\gamma _{k}(d^{\prime },Z)-\bar{\gamma}_{K}(d^{\prime },Z)\}.
\end{equation*}%
Therefore by the Cauchy Scwhartz, triangle, arithmetic mean-geometric mean
inequalities,%
\begin{align*}
& \left\vert \mathrm{E}[m(W,\gamma )-m(W,\bar{\gamma})-%
\sum_{k=1}^{K}D_{k}(W,\gamma _{k}-\bar{\gamma}_{k},\bar{\gamma})]\right\vert 
\\
& =\left\vert \mathrm{E}[\sum_{k=1}^{K-1}\{\gamma _{K}(d,k,Z)-\bar{\gamma}%
_{K}(d,k,Z)\}\{\gamma _{k}(d^{\prime },Z)-\bar{\gamma}_{k}(d^{\prime
},Z)\}]\right\vert  \\
& \leq \sum_{k=1}^{K-1}\left\{ \mathrm{E}[|\gamma _{K}(d,k,Z)-\bar{\gamma}%
_{K}(d,k,Z)|^{2}]+\mathrm{E}[\left\vert \gamma _{k}(d^{\prime },Z)-\bar{%
\gamma}_{k}(d^{\prime },Z)\right\vert ^{2}]\right\} \leq C\left\Vert \gamma -%
\bar{\gamma}\right\Vert ^{2},
\end{align*}%
where the last inequality follows similarly to previous results. The
conclusion now follows by Theorem 9. $\square .$

\bigskip

\textsc{Proof of Corollary 11}\textbf{: }Note first that for any $\gamma (X)$
it follows as in the proof of Corollary 7 that by $\Pr (D=1|Z)<1-C,$%
\begin{equation*}
\mathrm{E}[D\gamma (0,Z)^{2}]\leq \mathrm{E}[\gamma (0,Z)^{2}]=\mathrm{E}[%
\frac{1-D}{1-\pi _{0}(Z)}\gamma (0,Z)^{2}]=\mathrm{E}[\frac{1-D}{1-\pi
_{0}(Z)}\gamma (D,Z)^{2}]\leq C\mathrm{E}[\gamma (X)^{2}].
\end{equation*}%
Also note that 
\begin{equation*}
\mathrm{E}[D\gamma (0,Z)]=\mathrm{E}[\pi _{0}(Z)\gamma (0,Z)]=\mathrm{E}%
\left[ \pi _{0}(Z)\frac{1-D}{1-\pi _{0}(X)}\gamma (0,Z)\right] =\mathrm{E}%
[\alpha _{0}(X)\gamma (X)].
\end{equation*}%
The remainder of the proof follows analogously to the proof of Corollary 6. $%
\square .$

\bigskip

\subsection{Panel Average Derivative and Demand Elasticities}

\label{sec:delta}

Since own-price elasticity $\theta^{*}_{0}$ is a deterministic mapping of $%
\tilde{\theta}_{0}:=( \theta_{0}, \mathbb{E}[Y_{it}] )^{\prime}$, we obtain
the asymptotic variance $V^{*}$ of $\theta^{*}_{0}$ from the asymptotic
variance $\tilde{V}$ of $\tilde{\theta}_{0}$ using delta method.
Specifically, 
\begin{equation*}
V^{*}=H\tilde{V}H^{\prime}
\end{equation*}
where 
\begin{equation*}
H=\frac{\partial\theta^{*}_{0}}{\partial\tilde{\theta}_{0}}=%
\begin{bmatrix}
\frac{1}{\mathbb{E}[Y_{it}]} & \frac{-\theta_{0}}{[\mathbb{E}[Y_{it}]]^{2}}%
\end{bmatrix}%
\end{equation*}
and 
\begin{align*}
\tilde{V}=%
\begin{bmatrix}
\mathbb{E}[\psi_{0}(W_{it})]^{2} & \mathbb{E}[\psi_{0}(W_{it})Y_{it}] \\ 
\cdot & \mathbb{E}[Y_{it}]^{2}-\left\{ \mathbb{E}[Y_{it}]\right\} ^{2}%
\end{bmatrix}
.
\end{align*}

We estimate the asymptotic variance $V^{*}$ using the empirical analogue $%
\hat{V}^{*}$, where $\psi_{0}(W_{it})$ is replaced by 
\begin{align*}
\hat{\psi}_{it} & =\frac{\partial\hat{\gamma}_{\ell}(\tilde{X}_{i})}{%
\partial d}-\hat{\theta}+\hat{\alpha}^{*}_{\ell}(\tilde{X}_{i})[Y_{it}-\hat{%
\gamma}_{\ell}(\tilde{X}_{i})],\quad i\in I_{\ell}.
\end{align*}

The covariance estimator recognizes that household $i$'s observations form a
cluster $T_{i}$. For example, the estimator for the first component of $%
\tilde{V}$ is 
\begin{equation*}
\frac{1}{\sum_{i=1}^{n} T_{i}}\sum_{i=1}^{n} \sum_{t\in T_{i}}\sum_{s\in
T_{i}}\hat{\psi}_{it}\hat{\psi}_{is}.
\end{equation*}

More generally, we may consider estimating not only own price elasticity but
also income elasticity and cross price elasticity. The same arguments go
through with light modification.

Concatenate the derivatives as 
\begin{equation*}
\theta_{0}=%
\begin{bmatrix}
(\theta_{0})^{\text{income}} \\ 
(\theta_{0})^{\text{own}} \\ 
(\theta_{0})^{\text{cross}} \\ 
\end{bmatrix}
=%
\begin{bmatrix}
\mathrm{E} \left[ \frac{\partial\gamma_{0}(\tilde{X}_{i})}{\text{log income}}%
\right] \\ 
\mathrm{E} \left[ \frac{\partial\gamma_{0}(\tilde{X}_{i})}{\text{log own
price}}\right] \\ 
\mathrm{E} \left[ \frac{\partial\gamma_{0}(\tilde{X}_{i})}{\text{log cross
price}}\right] \\ 
\end{bmatrix}%
\end{equation*}
where the first and second components are scalars and the third component is
a vector.

The elasticities are a smooth transform thereof. By arguments in
Chernozhukov, Hausman, and Newey (2019) 
\begin{equation*}
\theta^{*}_{0}=%
\begin{bmatrix}
(\theta_{0}^{*})^{\text{income}} \\ 
(\theta_{0}^{*})^{\text{own}} \\ 
(\theta_{0}^{*})^{\text{cross}} \\ 
\end{bmatrix}
= 
\begin{bmatrix}
\frac{(\theta_{0})^{\text{income}}}{\mathbb{E}[Y_{it}]}-1 \\ 
\frac{(\theta_{0})^{\text{own}}}{\mathbb{E}[Y_{it}]}-1 \\ 
\frac{(\theta_{0})^{\text{cross}}}{\mathbb{E}[Y_{it}]}%
\end{bmatrix}
.
\end{equation*}

Likewise the delta method argument goes through. Elasticites $\theta_{0}$
are a deterministic mapping of $\tilde{\theta}_{0}=((\theta_{0}^{*})^{%
\prime}, \mathbb{E}[Y_{it}])^{\prime}$. We obtain the asymptotic variance $%
V^{*}$ of $\theta^{*}_{0}$ from the asymptotic variance $\tilde{V}$ of $%
\tilde{\theta }_{0}$ using delta method. Specifically, 
\begin{equation*}
V^{*}=H\tilde{V}H^{\prime}
\end{equation*}
where 
\begin{equation*}
H=\frac{\partial\theta^{*}_{0}}{\partial\tilde{\theta}_{0}}=%
\begin{bmatrix}
\frac{1}{\mathbb{E}[Y_{it}]}\cdot I & \frac{-\theta_{0}}{[\mathbb{E}%
[Y_{it}]]^{2}}%
\end{bmatrix}%
\end{equation*}
and $\tilde{V}$ is as before, where the influence function $\psi_{0}$ is
vector-valued, corresponding to the vector $\theta_{0}$.

As an aside, when using OLS, the empirical influence function used in
estimating off-diagonal terms is 
\begin{align*}
\psi_{0}(W_{it}) & =(\mathbb{E}[b_{it}b_{it}^{\prime-1}b_{it}\epsilon_{it} \\
\hat{\psi}_{it} & =\left( \frac{1}{\sum_{i=1}^{n} T_{i}}\sum_{i=1}^{n}%
\sum_{t\in T_{i}}\sum_{s\in T_{i}}b_{it}b_{it}^{\prime}\right) ^{-1}
b_{it}\epsilon_{it}
\end{align*}
where $\epsilon_{it}$ is the OLS residual for observation $W_{it}$. As
before, we use a variance estimator that recognizes clustering.

\section{Additional Empirical Results}

\label{sec:additional_empirics}

\subsection{Regression Decomposition and ATET}

We present ATET estimates from Auto-DML using cross validation rather than
theoretical iteration to tune the regularization. Our results are broadly
similar, with larger standard errors.

\begin{table}[ptb]
\centering
\begin{tabular}{c|c|c|c|c|c|c|c|c}
\hline \hline 
spec. & treated & untreated & Lasso ATET & Lasso SE & RF ATET & RF SE & NN
ATET & NN SE \\ \hline
1 & 185 & 172 & 4071.88 & 3390.11 & 4170.99 & 3277.92 & 1807.48 & 2656.05 \\ 
2 & 185 & 172 & 1618.74 & 500.49 & 2047.18 & 504.70 & 1754.79 & 531.10 \\ 
3 & 185 & 172 & 3379.15 & 1466.45 & 3589.10 & 1385.49 & 1175.21 & 1735.56 \\ 
\hline\hline
\end{tabular}
\caption{ATET using NSW treatment and NSW control, by cross validation}
\label{tab:ATT_nsw2_bias}
\end{table}

\begin{table}[ptb]
\centering
\begin{tabular}{c|c|c|c|c|c|c|c|c}
\hline\hline
spec. & treated & untreated & Lasso ATET & Lasso SE & RF ATET & RF SE & NN
ATET & NN SE \\ \hline
1 & 185 & 727 & 2194.07 & 1060.97 & 1986.25 & 1031.62 & 834.13 & 1004.58 \\ 
2 & 185 & 727 & 1686.64 & 1092.13 & 1422.68 & 1125.96 & 1909.87 & 1404.22 \\ 
3 & 185 & 727 & 2974.55 & 1108.72 & 2579.75 & 1042.94 & 3057.04 & 1454.26 \\ 
\hline\hline
\end{tabular}
\caption{ATET using NSW treatment and PSID comparison, by cross validation}
\label{tab:ATT_psid2_bias}
\end{table}

\begin{table}[ptb]
\centering
\begin{tabular}{c|c|c||c|c|c|c|c|c}
\hline\hline
spec. & treated & untreated & Lasso ATET & Lasso SE & RF ATET & RF SE & NN
ATET & NN SE \\ \hline
1 & 185 & 5904 & 1413.98 & 636.68 & 1813.82 & 662.06 & 2043.87 & 657.46 \\ 
2 & 185 & 5904 & 1405.09 & 644.10 & 1756.57 & 669.68 & 2025.50 & 653.32 \\ 
3 & 185 & 5904 & 1756.87 & 654.73 & 2013.84 & 676.72 & 1823.67 & 651.66 \\ 
\hline\hline
\end{tabular}
\caption{ATET using NSW treatment and CPS comparison, by cross validation}
\label{tab:ATT_cps2_bias}
\end{table}

\subsection{Panel Average Derivative and Demand Elasticities}

We present elasticity estimates from OLS with a simpler specification than
the specification used in the main text. We take as $b_{1}(X_{it})$ the
concatenation of the following variables: log expenditure, and log price of
each good. For $\tilde{H}_{i}$, we use the time averages of $b_{1}(X_{it})$.
Note that $K=dim(b_{1}(X_{it}))=16$ and $p=dim(b_{it})=288$. We calculate
clustered standard errors derived by delta method as explained in Appendix~%
\ref{sec:delta}. Tables~\ref{tab:reg_milk2} and~\ref{tab:reg_soda2}
summarize results.

\begin{table}[ptb]
\centering
\begin{tabular}{c|c|c}
\hline\hline
variable & elasticity & SE  \\ \hline
income & 0.42 & 0.05 \\ 
own-price & -0.68 & 0.05 \\ 
bread & -0.03 & 0.02 \\ 
butter & 0.00 & 0.02 \\ 
cereal & 0.00 & 0.02 \\ 
chips & 0.02 & 0.03 \\ 
coffee & 0.00 & 0.02 \\ 
cookies & 0.00 & 0.02 \\ 
eggs & -0.03 & 0.03 \\ 
ice cream & -0.03 & 0.03 \\ 
orange juice & -0.01 & 0.05 \\ 
salad & 0.02 & 0.02 \\ 
soda & -0.02 & 0.02 \\ 
soup & -0.03 & 0.02 \\ 
water & -0.01 & 0.02 \\ 
yogurt & 0.01 & 0.04 \\ 
\hline\hline
\end{tabular}
\vspace{8pt}
\caption{Milk elasticities, by OLS}
\label{tab:reg_milk2}
\end{table}

\begin{table}[ptb]
\centering
\begin{tabular}{c|c|c}
\hline\hline
variable & elasticity & SE  \\ \hline
income & 0.64 & 0.02 \\ 
own-price & -0.65 & 0.04 \\ 
bread & -0.01 & 0.03 \\ 
butter & -0.04 & 0.02 \\ 
cereal & 0.01 & 0.03 \\ 
chips & 0.02 & 0.03 \\ 
coffee & 0.01 & 0.02 \\ 
cookies & -0.03 & 0.02 \\ 
eggs & -0.03 & 0.03 \\ 
ice cream & 0.01 & 0.03 \\ 
milk & 0.02 & 0.04 \\ 
orange juice & -0.05 & 0.05 \\ 
salad & -0.03 & 0.02 \\ 
soup & 0.01 & 0.03 \\ 
water & 0.01 & 0.02 \\ 
yogurt & 0.05 & 0.04 \\ 
\hline\hline
\end{tabular}
\vspace{8pt}
\caption{Soda elasticities, by OLS}
\label{tab:reg_soda2}
\end{table}

%\bigskip

%Making the analogous replacement in the constraints of the Dantzig selector
%(Candes and Tao, 2007) gives a Dantzig estimator%

%\begin{equation}
%\hat{\rho}_{D}=\arg\min_{\rho}|\rho|_{1}\,s.t.|\hat{M}-\hat{G}\rho|_{\infty
%}\leq\lambda_{D}, \label{RRDantzig}%
%\end{equation}
%where $\lambda_{D}>0$ is the slackness size. These two minimization problems
%can be thought of as minimum distance versions of Lasso and Dantzig,
%respectively. Either $\hat{\rho}_{L}$ or $\hat{\rho}_{D}$ may be used in
%equation (\ref{Riesz est}) to form an estimator $\hat{\alpha}(x)=b(x)^{\prime
%}\hat{\rho}_{L}$ or $\hat{\alpha}(x)=b(x)^{\prime}\hat{\rho}_{D}$. This
%$\hat{\alpha}(x)$ may then be substituted in equation (\ref{Estimator}), along
%with a machine learner $\hat{\gamma}$ of the regression, to construct Auto-DML
%$\hat{\theta}$.

\pagebreak

\section*{References}

Ahn, H. and C.F. Manski (1993): \textquotedblleft Distribution Theory for
the Analysis of Binary Choice under Uncertainty with Nonparametric
Estimation of Expectations,\textquotedblright\ \textit{Journal of
Econometrics }56,\textit{\ }291--321.

Athey, S., G. Imbens, and S. Wager (2018): \textquotedblleft Approximate
Residual Balancing: Debiased Inference of Average Treatment Effects in High
Dimensions,\textquotedblright\ \textit{Journal of the Royal Statistical
Society, Series B }80,\textit{\ }597--623.

Avagyan, V. and S. Vansteelandt (2017): "Honest data-adaptive inference for
the average treatment effect under model misspecification using penalised
bias-reduced double-robust estimation," https://arxiv.org/abs/1708.03787.

Belloni, A., D. Chen, and V. Chernozhukov (2012): \textquotedblleft Sparse
Models and Methods for Optimal Instruments with an Application to Eminent
Domain,\textquotedblright\ \textit{Econometrica }80,\textit{\ }2369--429.

Belloni, A. and V. Chernozhukov (2013): "Least Squares After Model Selection
in High-dimensional Sparse Models," \textit{Bernoulli} 19, 521--547.

Belloni, A., V. Chernozhukov, and C. Hansen (2014a): \textquotedblleft
Inference on Treatment Effects after Selection among High-Dimensional
Controls,\textquotedblright\ \textit{Review of Economic Studies} 81,
608--650.

Belloni, A., V. Chernozhukov, L. Wang (2014b): \textquotedblleft Pivotal
Estimation via Square-Root Lasso in Nonparametric
Regression,\textquotedblright\ \textit{Annals of Statistics}\emph{\ }42,
757--788.

Belloni, A., V. Chernozhukov, K. Kato (2015): "Uniform Post-selection
Inference for Least Absolute Deviation Regression and Other Z-estimation
Problems," \textit{Biometrika} 102, 77--94.

Bickel, P.J. (1982): \textquotedblleft On Adaptive
Estimation,\textquotedblright\ \textit{Annals of Statistics} 10, 647--671.

Bickel, P.J. and Y. Ritov (1988): \textquotedblleft Estimating Integrated
Squared Density Derivatives: Sharp Best Order of Convergence
Estimates,\textquotedblright\ \textit{Sankhy\={a}: The Indian Journal of
Statistics, Series A }238, 381--393. \ 

Bickel, P.J., C.A.J. Klaassen, Y. Ritov and J.A. Wellner (1993): \textit{%
Efficient and Adaptive Estimation for Semiparametric Models}, Baltimore:
Johns Hopkins University Press.

Bickel, P.J., Y. Ritov, and A. Tsybakov (2009): \textquotedblleft
Simultaneous Analysis of {L}asso and {D}antzig Selector,\textquotedblright\ 
\textit{Annals of Statistics} 37, 1705--1732.

Blundell, R.W. and J.L. Powell (2004): "Endogeneity in Binary Response
Models," \textit{Review of Economic Studies} 71, 655-679.

Bradic, J. and M. Kolar (2017): \textquotedblleft Uniform Inference for
High-Dimensional Quantile Regression: Linear Functionals and Regression Rank
Scores,\textquotedblright\ arXiv:1702.06209.

Bradic, J., S. Wager, and Y. Zhu (2019): "Sparsity Double Robust Inference
of Average Treatment Effects," https://arxiv.org/pdf/1905.00744.pdf.

Bradic, J., V. Chernozhukov, W. Newey, and Y. Zhu (2019): "Minimax
Semiparametric Learning with Approximate Sparsity," arXiv.

Burda, M., M. Harding, J.A. Hausman (2008): \textquotedblleft A Bayesian
Mixed Logit Probit Model for Multinomial Choice,\textquotedblright\ \textit{%
Journal of Econometrics} 147, 232--46.

Burda, M., M. Harding, J.A. Hausman (2012): \textquotedblleft A Poisson
Mixture Model of Discrete Choice,\textquotedblright\ \textit{Journal of
Econometrics} 166, 184--203.

Cai, T.T. and Z. Guo (2017): ``Confidence Intervals for High-Dimensional
Linear Regression: Minimax Rates and Adaptivity," \textit{Annals of
Statistics }45, 615-646.

Candes, E. and T. Tao (2007): \textquotedblleft The Dantzig Selector:
Statistical Estimation when \textit{p} is much Larger than \textit{n}%
,\textquotedblright\ \textit{Annals of Statistics} 35, 2313--2351.

Cattaneo, M.D., M. Jansson, and W.K. Newey (2018): "Inference in Linear
Regression Models with Many Covariates and Heteroscedasticity," \textit{%
Journal of the American Statistical Association} 113, 1350-1361.

Chamberlain, G. (1982): \textquotedblleft Multivariate Regression Models for
Panel Data,\textquotedblright\ \textit{Journal of Econometrics }18,\textit{\ 
}5--46.

Chamberlain, G. (1982): \textquotedblleft Efficiency Bounds for
Semiparametric Regression,\textquotedblright\ \textit{Econometrica }60,%
\textit{\ }567--96.

Chamberlain, G. (1984): "Panel Data," \textit{Handbook of Econometrics Vol 2}%
, Z. Griliches and M. Intriligator, eds., 1247-1318.

Chatterjee, S. and J. Jafarov (2015): \textquotedblleft Prediction Error of
Cross-Validated Lasso,\textquotedblright\ arXiv:1502.06291.

Chen, X. and H. White (1999): "Improved Rates and Asymptotic Normality for
Nonparametric Neural Network Estimators," \textit{IEEE Transactions on
Information Theory} 45, 682-691.

Chernozhukov, V., D. Chetverikov, and K. Kato (2013a): \textquotedblleft
Gaussian Approximations and Multiplier Bootstrap for Maxima of Sums of
High-Dimensional Random Vectors,\textquotedblright\ \textit{Annals of
Statistics} 41, 2786--2819.

Chernozhukov, V., I. Fernandez-Val, J. Hahn, W. Newey (2013b):
\textquotedblleft Average and Quantile Effects in Nonseparable Panel
Models,\textquotedblright\ \textit{Econometrica }81,\textit{\ }535--80.

Chernozhkov, V., C. Hansen, and M. Spindler (2015): \textquotedblleft Valid
Post-Selection and Post-Regularization Inference: An Elementary, General
Approach,\textquotedblright\ \textit{Annual Review of Economics 7}, 649--688.

Chernozhukov, V., J. C. Escanciano, H. Ichimura, W.K. Newey, and J. Robins
(2016): \textquotedblleft Locally Robust Semiparametric
Estimation,\textquotedblright\ https://arxiv.org/abs/1608.00033v1.

Chernozhukov, V., D. Chetverikov, M. Demirer, E. Duflo, C. Hansen, W.K.
Newey (2017): "Double/Debiased/Neyman Machine Learning of Treatment
Effects," American Economic Review 107, 261-65.

Chernozhukov, V., D. Chetverikov, M. Demirer, E. Duflo, C. Hansen, W.K.
Newey, J.M. Robins (2018): "Double/debiased machine learning for treatment
and structural parameters," Econometrics Journal 21, C1-C68.

Chernozhukov, V., W.K. Newey, and J. Robins (2018): \textquotedblleft
Double/De-Biased Machine Learning Using Regularized Riesz Representers,"
https://arxiv.org/pdf/1802.08667v1.pdf.

Chernozhukov, V., W.K. Newey, and R. Singh (2018): "Learning L2-Continuous
Regression Functionals via Regularized Riesz Representers,"
https://arxiv.org/pdf/1809.05224v1.pdf.

Chernozhukov, V., W.K. Newey, and R. Singh (2019): "Double/De-Biased Machine
Learning of Global and Local Parameters Using Regularized Riesz
Representers," https://arxiv.org/abs/1802.08667v3.

Chernozhukov, V., J.A. Hausman, W.K. Newey (2021): "Demand Analysis with
Many Prices," NBER Working Paper 26424.

Chernozhukov, V., J. C. Escanciano, H. Ichimura, W.K. Newey, and J. Robins
(2020): \textquotedblleft Locally Robust Semiparametric
Estimation,\textquotedblright\ https://arxiv.org/abs/1608.00033v4.

Chiang, H.D., K. Kato, Y. Ma, Y. Sasaki (2019): \textquotedblleft Multiway
Cluster Robust Double/Debiased Machine Learning," arXiv:1909.03489.

Daubechies, I., M Defrise, and C. De Mol (2004): \textquotedblleft An
Iterative Thresholding Algorithm for Linear Inverse Problems with a Sparsity
Constraint,\textquotedblright\ \textit{Communications on Pure and Applied
Mathematics }57,\textit{\ }1413--57.

Dehejia, R.H. and S. Wahba (1999): \textquotedblleft Causal Effects in
Nonexperimental Studies: Reevaluating the Evaluation of Training
Programs,\textquotedblright\ \textit{Journal of the American Statistical
Association} 94 (448): 1053--62.

Farbmacher, M., M. Huber, L. Laff\'{e}rs, H. Langen, M. Spindler (2020):
"Causal Mediation Analysis with Double Machine Learning,"
https://arxiv.org/abs/2002.12710.

Farrell, M. (2015): \textquotedblleft Robust Inference on Average Treatment
Effects with Possibly More Covariates than Observations,\textquotedblright\ 
\textit{Journal of Econometrics} 189, 1--23.

Farrell, M., T. Liang, S. Misra (2021): "Deep Neural Networks for Estimation
and Inference," \textit{Econometrica} 89, 181--213.

Friedman, J., T. Hastie, H. H{\"{o}}fling, and R. Tibshirani (2007):
\textquotedblleft Pathwise Coordinate Optimization,\textquotedblright\ 
\textit{The Annals of Applied Statistics }1,\textit{\ }302--32.

Friedman, J., T. Hastie, and R. Tibshirani (2010): \textquotedblleft
Regularization Paths for Generalized Linear Models via Coordinate
Descent,\textquotedblright\ \textit{Journal of Statistical Software }33,%
\textit{\ }1-22.

Fu, W.J. (1998): \textquotedblleft Penalized Regressions: The Bridge versus
the Lasso,\textquotedblright\ \textit{Journal of Computational and Graphical
Statistics }7,\textit{\ }397--416.

Graham, B. and J.L. Powell (2012): \textquotedblleft Identification and
Estimation of Average Partial Effects in ``Irregular'' Correlated Random
Coefficient Panel Data Models,\textquotedblright\ \textit{Econometrica }80,%
\textit{\ }2105--52.

Hasminskii, R.Z. and I.A. Ibragimov (1979): \textquotedblleft On the
Nonparametric Estimation of Functionals,\textquotedblright\ in P. Mandl and
M. Huskova (eds.), \textit{Proceedings of the 2nd Prague Symposium on
Asymptotic Statistics, 21-25 August 1978}, Amsterdam: North-Holland, pp.
41-51.

Hausman, J.A. and W.K. Newey (2016): \textquotedblleft Individual
Heterogeneity and Average Welfare,\textquotedblright\ \textit{Econometrica}
84, 1225--1248.

Hirshberg, D.A. and S. Wager (2017): \textquotedblleft Balancing Out
Regression Error: Efficient Treatment Effect Estimation without Smooth
Propensities,\textquotedblright\ arXiv:1712.00038v1.

Hirshberg, D.A. and S. Wager (2020): \textquotedblleft Debiased Inference of
Average Partial Effects in Single-Index Models,\textquotedblright\ \textit{%
Journal of Business and Economic Statistics} 38, 19-24.

Hirshberg, D.A. and S. Wager (2018): \textquotedblleft Augmented minimax
linear estimation,\textquotedblright\ \newline
arXiv:1712.00038v5.

Huber, P. J.: "The Behavior of Maximum Likelihood Estimates Under
Nonstandard Conditions," in \textit{Proceedings of the Fifth Berkeley
Symposium in Mathematical Statistics and Probability}. Berkeley: University
of California Press, 1967.

Imai, K, L. Keele, and D. Tingley (2010): "A General Approach to Causal
Mediation Analysis," \textit{Psychological Methods }15, 309 --334.

Imbens, G.W. and W.K. Newey (2009): "Identification and Estimation of
Triangular Simultaneous Equations Models Without Additivity," \textit{%
Econometrica }77, 1481-1512.

Jankova, J. and S. Van De Geer (2015): \textquotedblleft Confidence
Intervals for High-Dimensional Inverse Covariance
Estimation,\textquotedblright\ \textit{Electronic Journal of Statistics} 90,
1205--1229.

Jankova, J. and S. Van De Geer (2016a): \textquotedblleft Semi-Parametric
Efficiency Bounds and Efficient Estimation for High-Dimensional
Models,\textquotedblright\ arXiv:1601.00815.

Jankova, J. and S. Van De Geer (2016b): \textquotedblleft Confidence Regions
for High-Dimensional Generalized Linear Models under
Sparsity,\textquotedblright\ arXiv:1610.01353.

Javanmard, A. and A. Montanari (2014a): \textquotedblleft Hypothesis Testing
in High-Dimensional Regression under the Gaussian Random Design Model:
Asymptotic Theory,\textquotedblright\ \textit{IEEE Transactions on
Information Theory} 60, 6522--6554.

Javanmard, A. and A. Montanari (2014b): \textquotedblleft Confidence
Intervals and Hypothesis Testing for High-Dimensional
Regression,\textquotedblright\ \textit{Journal of Machine Learning Research}
15: 2869--2909.

Javanmard, A. and A. Montanari (2015): \textquotedblleft De-Biasing the
Lasso: Optimal Sample Size for Gaussian Designs,\textquotedblright\
arXiv:1508.02757.

Jing, B.Y., Q.M. Shao, and Q. Wang (2003): \textquotedblleft Self-Normalized
Cram\'{e}r-Type Large Deviations for Independent Random
Variables,\textquotedblright\ \textit{Annals of Probability} 31, 2167--2215.

Kennedy, E.H. (2020): "Optimal Doubly Robust Estimation of Heterogeneous
Causal Effects," https://arxiv.org/abs/2004.14497.

Klaassen, C.A.J. (1987): "Consistent Estimation of the Influence Function of
Locally Asymptotically Linear Estimators," \textit{Annals ot Statistics} 15,
1548-1562.

LaLonde, R.J. (1986): ``Evaluating the Econometric Evaluations of Training
Programs with Experimental Data," \textit{The American Economic Review} 76,
604--20.

Leeb, H., and B.M. P\"{o}tscher (2008a): ``Recent Developments in Model
Selection and Related Areas," \textit{Econometric Theory} 24, 319--22.

Leeb H., and B.M. P\"{o}tscher (2008b): ``Sparse Estimators and the Oracle
Property, or the Return of Hodges' Estimator," \textit{Journal of
Econometrics} 142, 201--211.

Luo, Ye and M. Spindler (2016): "High-Dimenstional L2 Boosting: Rate of
Convergence," https://arxiv.org/pdf/1602.08927.pdf.

Luedtke, A. R. and M. J. van der Laan (2016): \textquotedblleft Optimal
Individualized Treatments in Resource-limited Settings," \textit{The
International Journal of Biostatistics} 12, 283-303.

Newey, W.K. (1994): \textquotedblleft The Asymptotic Variance of
Semiparametric Estimators,\textquotedblright\ \textit{Econometrica} 62,
1349--1382.

Newey, W.K., F. Hsieh, and J.M. Robins (1998): \textquotedblleft
Undersmoothing and Bias Corrected Functional Estimation,\textquotedblright\
MIT Dept. of Economics working paper\ 98-17.

Newey, W.K., F. Hsieh, and J.M. Robins (2004): \textquotedblleft Twicing
Kernels and a Small Bias Property of Semiparametric
Estimators,\textquotedblright\ \textit{Econometrica} 72, 947--962.

Newey, W.K. and J.M. Robins (2017): \textquotedblleft Cross Fitting and Fast
Remainder Rates for Semiparametric Estimation,\textquotedblright\
arXiv:1801.09138.

Neykov, M., Y. Ning, J.S. Liu, and H. Liu (2015): \textquotedblleft A
Unified Theory of Confidence Regions and Testing for High Dimensional
Estimating Equations,\textquotedblright\ arXiv:1510.08986.

Ning, Y. and H. Liu (2017): \textquotedblleft A General Theory of Hypothesis
Tests and Confidence Regions for Sparse High Dimensional
Models,\textquotedblright\ \textit{Annals of Statistics} 45, 158-195.

Powell, J.L., J.H. Stock, and T.M. Stoker (1989): "Semiparametric Estimation
of Index Coefficients," \textit{Econometrica }57, 1403-1430.

Ren, Z., T. Sun, C.H. Zhang, and H. Zhou (2015): \textquotedblleft
Asymptotic Normality and Optimalities in Estimation of Large Gaussian
Graphical Models,\textquotedblright\ \textit{Annals of Statistics} 43,
991--1026.

Robins, J.M. and A. Rotnitzky (1995): \textquotedblleft Semiparametric
Efficiency in Multivariate Regression Models with Missing
Data,\textquotedblright\ \textit{Journal of the American Statistical
Association} 90 (429): 122--129.

Robins, J.M., A. Rotnitzky, and L.P. Zhao (1995): \textquotedblleft Analysis
of Semiparametric Regression Models for Repeated Outcomes in the Presence of
Missing Data,\textquotedblright\ \textit{Journal of the American Statistical
Association }90, 106--121.

Robins, J., P. Zhang, R. Ayyagari, R. Logan, E. Tchetgen, L. Li, A. Lumley,
and A. van der Vaart (2013): \textquotedblleft New Statistical Approaches to
Semiparametric Regression with Application to Air Pollution Research,"
Research Report Health E Inst.

Rosenbaum, P.R. and D. B. Rubin (1983): \textquotedblleft The Central Role
of the Propensity Score in Observational Studies for Causal
Effects,\textquotedblright\ \textit{Biometrika }70: 41--55.

Rothenh{\"{a}}usler, D. and B. Yu (2019): \textquotedblleft Incremental
Causal Effects,\textquotedblright\ arXiv:1907.13258.

Rudelson, M. and S. Zhou (2013): "Reconstruction From Anisotropic Random
Measurements," \textit{IEEE Transactions on Informating Theory} 59,
3434--3447.

Scharfstein, D.O., A. Rotnitzky, and J.M. Robins (1999): "Rejoinder to
Adjusting for Nonignorable Drop-out Using Semiparametric Nonresponse Models," 
\textit{Journal of the American Statistical Association} 94, 1096-1146.

Schick, A. (1986): \textquotedblleft On Asymptotically Efficient Estimation
in Semiparametric Models,\textquotedblright\ \textit{Annals of Statistics}
14, 1139--1151.

Schmidt-Hieber, J. (2020): "Nonparametric Regression Using Deep Neural
Networks with RELU Activation Function," \textit{The Annals of Statistics}
48, 1875--1897.

Singh, R. and L. Sun (2019): \textquotedblleft De-biased Machine Learning
for Compliers,\textquotedblright\ arXiv:1909.05244.

Smucler, E., A. Rotnitzky, and J.R. Robins (2019): "A Unifying Approach for
Doubly-robust L1 Regularized Estimation of Causal Contrasts,"
https://arxiv.org/abs/1904.03737.

Stock, J.H. (1989): \textquotedblleft Nonparametric Policy
Analysis,\textquotedblright\ \textit{Journal of the American Statistical
Association} 84, 567--575.

Syrgkanis, V., and M. Zampetakis (2020): "Estimation and Inference with
Trees and Forests in High Dimensions," https://arxiv.org/abs/2007.03210.

Tchetgen Tchetgen, E.J. and I. Shipster (2012): "Semiparametric Theory for
Causal Mediation Analysis: Efficiency Bounds, Multiple Robustness and
Sensitivity Analysis," \textit{The Annals of Statistics} 40, 1816-1845.

Toth, B. and M. J. van der Laan (2016), \textquotedblleft TMLE for Marginal
Structural Models Based On An Instrument," U.C. Berkeley Division of
Biostatistics Working Paper Series, Working Paper 350.

Tseng, P. (2001): \textquotedblleft Convergence of a Block Coordinate
Descent Method for Nondifferentiable Minimization,\textquotedblright\ 
\textit{Journal of Optimization Theory and Applications} 109,\textit{\ }%
475--94.

Van De Geer, S., P. B{\"{u}}hlmann, Y. Ritov, and R. Dezeure (2014):
\textquotedblleft On Asymptotically Optimal Confidence Regions and Tests for
High-Dimensional Models,\textquotedblright\ \textit{Annals of Statistics},
42: 1166--1202.

Van der Laan, M. and D. Rubin (2006): \textquotedblleft Targeted Maximum
Likelihood Learning,\textquotedblright\ \textit{International Journal of
Biostatistics} 2.

Van der Laan, M. J. and S. Rose (2011): \textit{Targeted Learning: Causal
Inference for Observational and Experimental Data,} Springer.

Van der Vaart, A.W. (1991): \textquotedblleft On Differentiable
Functionals,\textquotedblright\ \textit{Annals of Statistics}, 19: 178--204.

Van der Vaart, A.W. (1998): \textit{Asymptotic Statistics}. New York:
Cambridge University Press.

Vermeulen, K. and S. Vansteelandt (2015): "Bias-Reduced Doubly Robust
Estimation," \textit{Journal of the American Statistical Association} 110,
1024-1036.

Vershynin, R. (2018): \textit{High-Dimensional Probability}, New York:
Cambridge University Press.

White, H. (1982): "Maximum Likelihood Estimation of Misspecified Models," 
\textit{Econometrica} 50, 1-25.

Wooldridge, J.M. (2002): \textit{Econometric Analysis of Cross-Section and
Panel Data}, Cambridge, MIT Press.

Wooldridge, J.M. (2019): \textquotedblleft Correlated Random Effects Models
with Unbalanced Panels,\textquotedblright\ \textit{Journal of Econometrics}
211, 137--50.

Wooldridge, J.M. and Y. Zhu (2020): "Inference in Approximately Sparse
Correlated Random Effects Probit Models With Panel Data," \textit{Journal of
Business and Economic Statistics} 38, 1-18.

Zhang, C. and S. Zhang (2014): \textquotedblleft Confidence Intervals for
Low-Dimensional Parameters in High-Dimensional Linear
Models,\textquotedblright\ \textit{Journal of the Royal Statistical Society,
Series B }76, 217--242.

Zheng, W., Z. Luo, and M. J. van der Laan (2016), ``Marginal Structural
Models with Counterfactual Effect Modifiers," U.C. Berkeley Division of
Biostatistics Working Paper Series, Working Paper 348.

Zhu, Y. and J. Bradic (2017a): \textquotedblleft Linear Hypothesis Testing
in Dense High-Dimensional Linear Models,\textquotedblright\ \textit{Journal
of the American Statistical Association} 112.

Zhu, Y. and J. Bradic (2017b): \textquotedblleft Breaking the Curse of
Dimensionality in Regression,\textquotedblright\ arXiv: 1708.00430.

Zubizarreta, J.R. (2015): \textquotedblleft Stable Weights that Balance
Covariates for Estimation with Incomplete Outcome Data,\textquotedblright\ 
\textit{Journal of the American Statistical Association} 90 (429): 122--129.

\end{document}